\documentclass[10pt,a4paper]{article}
\usepackage{amsmath,amssymb,amsthm,graphics,graphicx,epic,eepic,multicol,ascmac}
\usepackage[mathscr]{euscript}
\usepackage[all]{xy}

\numberwithin{equation}{section}
\theoremstyle{theorem}
\newtheorem{thm}{Theorem}[section]
\newtheorem{prop}[thm]{Proposition}
\newtheorem{lem}[thm]{Lemma}
\newtheorem{rem}[thm]{Remark}

\newtheorem{ex}[thm]{Example}

\theoremstyle{definition}
\newtheorem{defn}[thm]{Definition}

\def\al{\alpha}

\def\wht(#1){\widehat{\ #1\ }}

\newcommand{\cA}{{\mathcal A}}

\newcommand{\cF}{{\mathcal F}}

\newcommand{\frg}{\mathfrak g}
\newcommand{\frh}{\mathfrak h}

\newcommand{\frn}{\mathfrak n}

\newcommand{\bbC}{\mathbb C}

\newcommand{\ch}{\mathrm{ch}}

\newcommand{\lbr}{\begin{bmatrix}}
\newcommand{\rbr}{\end{bmatrix}}

\def\ge{\frg}

\def\al{\alpha}

\def\beneme{\begin{enumerate}}
\def\beq{\begin{equation}}
\def\beqn{\begin{eqnarray}}
\def\beqnn{\begin{eqnarray*}}

\def\bfii0{{\bf i_0}}

\def\bbra#1,#2,#3{\left\{\begin{array}{c}\hspace{-5pt}
#1;#2\\ \hspace{-5pt}#3\end{array}\hspace{-5pt}\right\}}

\def\ci(#1,#2){c_{#1}^{(#2)}}
\def\Ci(#1,#2){C_{#1}^{(#2)}}
\def\mpp(#1,#2,#3){#1^{(#2)}_{#3}}
\def\bCi(#1,#2){\ovl C_{#1}^{(#2)}}
\def\ch(#1,#2){c_{#2,#1}^{-h_{#1}}}
\def\cc(#1,#2){c_{#2,#1}}

\def\di(#1,#2){D_{#1}^{(#2)}}
\def\dbi(#1,#2){\ovl D_{#1}^{(#2)}}

\def\eneme{\end{enumerate}}

\def\eeq{\end{equation}}
\def\eeqn{\end{eqnarray}}
\def\eeqnn{\end{eqnarray*}}

\def\gau#1,#2{\left[\begin{array}{c}\hspace{-5pt}#1\\
\hspace{-5pt}#2\end{array}\hspace{-5pt}\right]}

\def\ji(#1,#2){j_{#1}^{(#2)}}

\def\lan{\langle}

\def\lm{\lambda}
\def\Lm{\Lambda}

\def\nd{\noindent}

\def\ovl{\overline}

\def\qed{\hfill\framebox[2mm]{}}

\def\ran{\rangle}

\def\TY(#1,#2,#3){#1^{(#2)}_{#3}}

\def\xxi(#1,#2,#3){\displaystyle {}^{#1}\Xi^{(#2)}_{#3}}
\def\xsi(#1,#2,#3){\displaystyle {}^{#1}\Sigma^{(#2)}_{#3}}
\def\xE(#1,#2,#3){\displaystyle {}^{#1}E_{#2}[#3]}
\def\xF(#1,#2){\displaystyle {}^{#1}F_{#2}}
\def\xx(#1,#2){\displaystyle {}^{#1}\Xi_{#2}}
\def\W1{W(\varpi_1)}

\def\m@th{\mathsurround=0pt}
\def\fsquare(#1,#2){
\hbox{\vrule$\hskip-0.4pt\vcenter to #1{\normalbaselines\m@th
\hrule\vfil\hbox to #1{\hfill$\scriptstyle #2$\hfill}\vfil\hrule}$\hskip-0.4pt
\vrule}}

\newcommand{\ba}{\begin{array}}
\newcommand{\ea}{\end{array}}

\newcommand{\eq}{\begin{eqnarray}}
\newcommand{\eneq}{\end{eqnarray}}

\title{\textbf{\large{Cluster algebras of finite type via a Coxeter element and Demazure Crystals of type A}}}
\author{\normalsize{YUKI KANAKUBO\thanks{Division of Mathematics, 
Sophia University, Kioicho 7-1, Chiyoda-ku, Tokyo 102-8554,
Japan: {j\_chi\_sen\_you\_ky@eagle.sophia.ac.jp}}
\ and\ 
TOSHIKI NAKASHIMA\thanks{Division of Mathematics, 
Sophia University, Kioicho 7-1, Chiyoda-ku, Tokyo 102-8554,
Japan: { toshiki@sophia.ac.jp}:
supported in part by JSPS Grants 
in Aid for Scientific Research $15K04794$.}
}}
\date{}

\begin{document}

\maketitle
\vspace{-10pt}

\begin{abstract}
Let $G$ be a simply connected simple algebraic group over $\mathbb{C}$, 
$B$ and $B_-$ be its two opposite Borel subgroups. 
For two elements $u$, $v$ of the Weyl group $W$, 
it is known that the coordinate ring 
${\mathbb C}[G^{u,v}]$ of the double Bruhat cell $G^{u,v}=BuB\cap B_-vB_-$ is 
isomorphic to a cluster algebra $\cA(\textbf{i})_{{\mathbb C}}$ \cite{A-F-Z,GY}. In the case $u=e$, $v=c^2$ ($c$ is a Coxeter element), the algebra ${\mathbb C}[G^{e,c^2}]$ has only finitely many cluster variables. In this article, for $G={\rm SL}_{r+1}(\mathbb{C})$, we obtain explicit forms of all the cluster variables in $\mathbb{C}[G^{e,c^2}]$ by considering its additive categorification via preprojective algebras,
and describe them in terms of monomial realizations of Demazure crystals.
\end{abstract}


\tableofcontents

\section{Introduction}

A cluster algebra is a commutative ring generated by so-called ``cluster 
variables'', which has been introduced in order to study certain combinatorial 
properties of dual (semi) canonical bases by Fomin and Zelevinsky(\cite{FZ2}).
Nowadays, it has influenced to remarkably wide areas of mathematics and 
physics.

In \cite{A-F-Z}, Berenstein et al constructed the upper cluster algebra
structures on the
coordinate algebra $\bbC[G^{u,v}]$ of double Bruhat cell $G^{u,v}$, where 
$G$ is a simply-connected simple algebraic group over $\mathbb{C}$ and $u,v$ are elements
of the associated Weyl group $W$. Recently, Goodearl and Yakimov 
showed that $\bbC[G^{u,v}]$ also has a cluster algebra structure (\cite{GY}).
In \cite{GLS} Geiss et al initiated categorification of cluster algebras by 
considering semi-canonical bases. 

Cluster algebras which have only finitely many cluster variables are called {\it finite type}. In \cite{FZ3}, cluster algebras of finite type are studied thoroughly, and they are classified by the set of Cartan matrices up to coefficients. For a fixed Cartan matrix, all the cluster variables are parametrized by 
the set of ``almost positive roots'', which is, a union of all positive roots and negative simple roots corresponding to the Cartan matrix. Thus, we define the type of such cluster algebra to be the type of the corresponding Cartan matrix. Let $c\in W$ be a Coxeter element whose length $l(c)$ satisfies 
$l(c^2)=2 l(c)=2$rank$(G)$. It is known that one can realize a cluster algebra of finite type on the coordinate ring $\bbC[G^{e,c^2}]$, whose type coincides with the Cartan-Killing type of $G$ \cite{A-F-Z}. 

In \cite{KaN,KaN2}, we showed that certain cluster variables of $\bbC[G^{u,e}]$ ($u\in W$)
are realized as a sum of monomials in Demazure crystals in the case $G$ is type A, B, C or D. Then we treated only a part of the cluster variables, and calculated them by using an inner product
on the fundamental representation of Lie(G). In this paper, we consider the case $G={\rm SL}_{r+1}(\mathbb{C})$ $(r\geq 3)$ and describe all the cluster variables in $\bbC[G^{e,c^2}]$ using direct sum of 
certain monomial realizations of Demazure crystals and give a new parametrization of these cluster variables different from those of \cite{FZ3}. For the proof, we use the {\it additive categorification} of the coordinate ring $\mathbb{C}[L^{e,c^2}]$ (\cite{GLS}). Each cluster in $\mathbb{C}[L^{e,c^2}]$ is associated with a {\it cluster-tilting module} of the {\it preprojective algebra} (see Sect.4), and each cluster variable
is associated with a direct summand of the corresponding cluster-tilting module. There
exists an explicit formula of such cluster variables (Proposition \ref{dualsemi}). With the help of this formula and the 
additive categorification, we shall prove that each initial cluster variable in $\mathbb{C}[G^{e,c^2}]$ is described as a sum of monomials in the Demazure crystal $B(\Lambda_{k})_{w_k}$
with some $k\in\{1,2,\cdots,r\}$ and $w_k\in W$. And other variables are described as sums of monomials in Demazure crystals in the forms
$B(\sum^b_{s=a}\Lambda_{s})_w\oplus \bigoplus^{p}_{t=1} B(\lambda_t)_{w_t}$ with some $w,w_t\in W$, $p,a,b\in \mathbb{Z}_{>0}$ and 
$\lambda_t\in \sum^b_{s=a}\Lambda_{s} - \sum_{i\in I}\mathbb{Z}_{\geq0} \al_i$. As a corollary of these results, we see that a natural correspondence 
$-\al_k\mapsto B(\Lambda_{k})_{w_k}$, $\sum^b_{s=a}\al_{s}\rightarrow B(\sum^b_{s=a}\Lambda_{s})_w\oplus \bigoplus^{p}_{t=1} B(\lambda_t)_{w_t}$
gives a parametrization of the cluster variables in $\mathbb{C}[G^{e,c^2}]$ by the set of almost positive roots.

For example, let us consider the case $G={\rm SL}_{4}(\mathbb{C})$ (type ${\rm A}_3$ algebraic group). For the monomial realization of the crystal $B(\Lambda_2)$ of type ${\rm A}_3$, its crystal graph is as follows:
\[
\begin{xy}
(-20,87)*{Y_{1,2}}="h",
(0,87)*{\frac{Y_{1,1}Y_{1,3}}{Y_{2,2}}}="1j",
(0,74)*{\frac{Y_{1,1}}{Y_{2,3}}}="1j-1",
(20,87)*{\frac{Y_{1,3}}{Y_{2,1}}}="1j+1",
(20,74)*{\frac{Y_{2,2}}{Y_{2,1}Y_{2,3}}}="1j+1-1",
(40,74)*{\frac{1}{Y_{3,2}}}="1j+1-10",
\ar@{->} "h";"1j"^{2}
\ar@{->} "1j";"1j+1"^{1}
\ar@{->} "1j";"1j-1"_{3}
\ar@{->} "1j+1";"1j+1-1"^{3}
\ar@{->} "1j-1";"1j+1-1"_{1}
\ar@{->} "1j+1-1";"1j+1-10"_{2}
\end{xy}
\]
On the other hand, taking a Coxeter element $c=s_2s_1s_3\in W$, specific initial cluster variables in 
$\mathbb{C}[G^{e,c^2}]$ are given by generalized minors $\Delta_{\Lambda_i,c^{-1}\Lambda_i}$ $(i=1,2,3)$ (see \ref{bilingen}). It is known that $\Delta_{\Lambda_1,c^{-1}\Lambda_1}$, $\Delta_{\Lambda_2,c^{-1}\Lambda_2}$ and $\Delta_{\Lambda_3,c^{-1}\Lambda_3}$ coincide with ordinary minors $D_{1,2}$, $D_{12,24}$ and $D_{123,124}$, respectively, where $D_{\{1,2,\cdots,k\},\{i_1,i_2,\cdots,i_k\}}$ denote the minor whose rows are labelled by $\{1,2,\cdots,k\}$, columns are labelled by $\{i_1,i_2,\cdots,i_k\}$.
Using the biregularly isomorphism $\ovl{x}^G_{\textbf{i}}:H\times (\mathbb{C}^{\times})^{6}\rightarrow G^{e,c^2}$ ($\textbf{i}:=(2,1,3,2,1,3)$) in Proposition
 \ref{gprime}, we have
\[ D_{12,24}\circ \ovl{x}^G_{\textbf{i}}(a;\textbf{Y})=
a_1a_2(Y_{1,2}+\frac{Y_{1,1}Y_{1,3}}{Y_{2,2}}+\frac{Y_{1,3}}{Y_{2,1}}+\frac{Y_{2,2}}{Y_{2,1}Y_{2,3}}+\frac{Y_{1,1}}{Y_{2,3}}),
\]
where we set $a:={\rm diag}(a_1,a_2,a_3,a_4)\in H$ and $\textbf{Y}:=(Y_{1,2},Y_{1,1},Y_{1,3},Y_{2,2},Y_{2,1},Y_{2,3})\in(\mathbb{C}^{\times})^{6}$. Comparing with the above crystal graph of $B(\Lambda_2)$,
we see that the set of terms $\{Y_{1,2},\frac{Y_{1,1}Y_{1,3}}{Y_{2,2}},\frac{Y_{1,3}}{Y_{2,1}},\frac{Y_{2,2}}{Y_{2,1}Y_{2,3}},\frac{Y_{1,1}}{Y_{2,3}}\}$ in $D_{12,24}\circ \ovl{x}^G_{\textbf{i}}$ coincides with the monomial realization of the Demazure crystal $B(\Lambda_2)_{s_3s_1s_2}$ (see \ref{Demcrysub}). Similarly, we get
\[
D_{1,2}\circ \ovl{x}^G_{\textbf{i}}(a;\textbf{Y})=
a_1(Y_{1,1}+\frac{Y_{2,2}}{Y_{2,1}}),\quad
D_{123,124}\circ \ovl{x}^G_{\textbf{i}}(a;\textbf{Y})=
a_1a_2a_3(Y_{1,3}+\frac{Y_{2,2}}{Y_{2,3}}),
\]
which coincide with the total sums of monomials in Demazure crystals $B(\Lambda_1)_{s_3s_1s_2}$, $B(\Lambda_3)_{s_3s_1s_2}$ respectively. All other cluster variables in $\mathbb{C}[G^{e,c^2}]$ are
\[ (D_{12,12}D_{13,34})\circ \ovl{x}^G_{\textbf{i}}=a_1^2a_2a_3 Y_{2,2},\quad D_{12,14}\circ \ovl{x}^G_{\textbf{i}}=a_1a_2(Y_{1,2}Y_{2,1}+\frac{Y_{1,1}Y_{1,3}Y_{2,1}}{Y_{2,2}}+\frac{Y_{1,1}Y_{2,1}}{Y_{2,3}}), \]
\[ D_{12,23}\circ \ovl{x}^G_{\textbf{i}}=a_1a_2(Y_{1,2}Y_{2,3}+
\frac{Y_{1,1}Y_{1,3}Y_{2,3}}{Y_{2,2}}+\frac{Y_{1,3}Y_{2,3}}{Y_{2,1}}),\quad
D_{1,3}\circ \ovl{x}^G_{\textbf{i}}=a_1Y_{2,3},
 \]
\[ D_{123,134}\circ \ovl{x}^G_{\textbf{i}}=a_1a_2a_3Y_{2,1},\quad
D_{12,13}\circ \ovl{x}^G_{\textbf{i}}=a_1a_2(Y_{1,2}Y_{2,1}Y_{2,3}+
\frac{Y_{1,1}Y_{1,3}Y_{2,1}Y_{2,3}}{Y_{2,2}}),
 \]
which coincide with the total sums of monomials in Demazure crystals $B(\Lambda_2)_e$, $B(\Lambda_1+\Lambda_2)_{s_3s_2}$, $B(\Lambda_2+\Lambda_3)_{s_1s_2}$, $B(\Lambda_3)_e$, $B(\Lambda_1)_e$ and $B(\Lambda_1+\Lambda_2+\Lambda_3)_{s_2}$ respectively. These results imply the statement of Theorem \ref{thm1} for $r=3$. Thus, the correspondence $-\al_i\mapsto B(\Lambda_i)_{s_3s_1s_2}$, $\al_i\mapsto B(\Lambda_i)_{e}$ $(i=1,2,3)$,
$\al_{1}+\al_2\mapsto B(\Lambda_1+\Lambda_2)_{s_3s_2}$, $\al_{2}+\al_3\mapsto B(\Lambda_2+\Lambda_3)_{s_1s_2}$, and $\al_{1}+\al_2+\al_3\mapsto B(\Lambda_1+\Lambda_2+\Lambda_3)_{s_2}$ yields an alternative parametrization of all cluster variables in $\mathbb{C}[G^{e,c^2}]$ by the set of almost positive roots, which differs from the one in \cite{FZ3}. This correspondence means the claim in Theorem \ref{maincor} for $r=3$.

The article is organized as follows. In section 2, we recall properties of (reduced) double Bruhat cells $G^{u,v}$ and $L^{u,v}$. In section 3, 
after a concise reminder on cluster algebras, we review an isomorphism between the coordinate ring of a double Bruhat cell 
$G^{e,v}$ and a cluster algebra $\mathcal{A}(\textbf{i})$. In section 4, we recall the cluster algebra structure of $\mathbb{C}[L^{e,v}]$ and basic notions of preprojective algebras. We also review the additive categorifications of the cluster algebras $\mathbb{C}[L^{e,v}]$ following \cite{LD,GLS}. In section 5, we shortly review 
the definition of monomial realizations of crystal bases. Section 6 is devoted to present our main results, which provide a relation between all cluster variables in $\mathbb{C}[G^{e,c^2}]$ and monomial 
realizations of Demazure crystals. In Section 7, we complete the proof of the main theorems.

\section{Factorization theorem}\label{DBCs}

In this section, we shall introduce (reduced) double Bruhat cells $G^{u,v}$, $L^{u,v}$, and their properties\cite{B-Z, F-Z}. For $l\in \mathbb{Z}_{>0}$, we set $[1,l]:=\{1,2,\cdots,l\}$.

\subsection{Double Bruhat cells}\label{factpro}

Let $G$ be a simple complex algebraic group of classical type, $B$ and $B_-$ be two opposite Borel subgroups in $G$, $N\subset B$ and $N_-\subset B_-$ be their unipotent radicals, 
$H:=B\cap B_-$ a maximal torus. We set $\frg:={\rm Lie}(G)$ with the triangular decomposition $\frg=\frn_-\oplus \frh \oplus \frn$. Let $e_i$, $f_i$ $(i\in[1,r])$ be the generators of $\frn$, $\frn_-$. For $i\in[1,r]$ and $t \in \mathbb{C}$, we set
\begin{equation}\label{xiyidef} 
x_i(t):={\rm exp}(te_i),\ \ \ y_{i}(t):={\rm exp}(tf_i).
\end{equation}
Let $W:=\lan s_i |i=1,\cdots,r \ran$ be the Weyl group of $\frg$, where
$\{s_i\}$ are the simple reflections. We identify the Weyl group $W$ with ${\rm Norm}_G(H)/H$. An element 
\begin{equation}\label{smpl}
\ovl{s_i}:=x_i(-1)y_i(1)x_i(-1)
\end{equation}
is in ${\rm Norm}_G(H)$, which is a representative of $s_i\in W={\rm Norm}_G(H)/H$ \cite{N1}. For $u\in W$, let $u=s_{i_1}\cdots s_{i_n}$ be its reduced expression. Then we write $\ovl{u}=\ovl{s_{i_1}}\cdots \ovl{s_{i_n}}$, call $l(u):=n$ the length of $u$. We have two kinds of Bruhat decompositions of $G$ as follows:
\[ G=\displaystyle\coprod_{u \in W}B\ovl{u}B=\displaystyle\coprod_{u \in W}B_-\ovl{u}B_- .\]
Then, for $u$, $v\in W$, 
we define the {\it double Bruhat cell} $G^{u,v}$ as follows:
\[ G^{u,v}:=B\ovl{u}B \cap B_-\ovl{v}B_-. \]
We also define the {\it reduced double Bruhat cell} $L^{u,v}$ as follows:
\[ L^{u,v}:=NuN \cap B_-vB_- \subset G^{u,v}. \] 

\begin{defn}\label{redworddef}
Let $v=s_{j_n}\cdots s_{j_1}$ be a reduced expression of $v\in W$ $(j_n,\cdots,j_1\in [1,r])$. Then the finite sequence $\textbf{i}:=(j_n,\cdots,j_1)$ is called a {\it reduced word} for $v$.
\end{defn}

For example, the sequence $(2,1,3,2,1,3)$ is a reduced word of the longest element $s_2s_1s_3s_2s_1s_3$ of the Weyl group of type ${\rm A}_3$. In this paper, we mainly treat (reduced) double Bruhat cells of the form $G^{e,v}:=B \cap B_-\ovl{v} B_-$, $L^{e,v}:=N \cap B_-v B_-$.

\subsection{Factorization theorem}\label{fuctorisec}

In this subsection, we shall introduce the isomorphisms between double Bruhat cell $G^{e,v}$ and $H\times (\mathbb{C}^{\times})^{l(v)}$, and between $L^{e,v}$ and $(\mathbb{C}^{\times})^{l(v)}$. For $i \in [1,r]$ and $t\in \mathbb{C}^{\times}$, we set $\alpha_i^{\vee}(t):=t^{h_i}$. 

For a reduced word $\textbf{i}=(i_1, \cdots ,i_n)$ 
($i_1,\cdots,i_n\in[1,r]$), 
we define a map $x^G_{\textbf{i}}:H\times \mathbb{C}^n \rightarrow G$ as 
\begin{equation}\label{xgdef}
x^G_{\textbf{i}}(a; t_1, \cdots, t_n):=a\cdot x_{i_1}(t_1)\cdots x_{i_n}(t_n).
\end{equation}

\begin{thm}\label{fp}${\cite{B-Z,F-Z}}$ For $v\in W$ and its reduced word ${\rm \bf{i}}$, the map $x^G_{{\rm \bf{i}}}$ is a biregular isomorphism from $H\times (\mathbb{C}^{\times})^{l(v)}$ to a Zariski open subset of $G^{e,v}$. The map $(\mathbb{C}^{\times})^{l(v)}\rightarrow L^{e,v}$, $(t_1, \cdots, t_n)\mapsto x^G_{{\rm \bf{i}}}(1;t_1,\cdots,t_n)$ is a biregular isomorphism to a Zariski open subset of $L^{e,v}$.
\end{thm}

For $\textbf{i}=(i_1, \cdots ,i_n)$
($i_1,\cdots,i_n\in[1,r]$), we define a map
$\ovl{x}^G_{\textbf{i}}:H\times(\mathbb{C}^{\times})^{n}\rightarrow
G^{e,v}$ as
\[ \ovl{x}^G_{\textbf{i}}(a;t_1,\cdots,t_n)
=a x_{i_1}(t_1)\alpha_{i_1}^{\vee}(t_1)x_{i_2}(t_2)\alpha_{i_2}^{\vee}(t_2)\cdots x_{i_n}(t_n)\alpha_{i_n}^{\vee}(t_n) , \]
where $a\in H$ and $(t_1,\cdots,t_n)\in (\mathbb{C}^{\times})^{n}$. 

Now, let $G={\rm SL}_{r+1}(\mathbb{C})$ and $c\in W$ be a Coxeter element such that a reduced word ${\rm \bf{i}}$ of $c^2$ can be written as
\begin{equation}\label{redwords2}
{\rm \bf{i}}=
\begin{cases}
(2,4,6,\cdots,r,1,3,5,\cdots,r-1,2,4,6,\cdots,r,1,3,5,\cdots,r-1) & {\rm if}\ r\ {\rm is\ even},\\ 
(2,4,6,\cdots,r-1,1,3,5,\cdots,r,2,4,6,\cdots,r-1,1,3,5,\cdots,r) & {\rm if}\ r\ {\rm is\ odd}.
\end{cases} 
\end{equation}

\begin{rem}\label{importantrem}
In the rest of the paper, we use double indexed variables $Y_{s,j}$ $(s\in \mathbb{Z}$, $j\in [1,r])$. If we see the variables $Y_{s,0}$, $Y_{s,j}$ 
$(r+1\leq j)$ then 
we understand $Y_{s,0}=Y_{s,j}=1$. For example, if $l=1$ then $Y_{s,l-1}=1$.
\end{rem}

\begin{prop}\label{gprime} 
In the above setting, the map $\ovl{x}^G_{{\rm \bf{i}}}$ is a biregular isomorphism between $H\times(\mathbb{C}^{\times})^{2r}$ and a Zariski open subset of $G^{e,c^2}$.
\end{prop}
\nd
{\sl [Proof.]}

Let $j_k$ be the $k$ th index of \textbf{i} in (\ref{redwords2}) from the right, which means that
$\textbf{i}=(j_{2r},\cdots,j_{r+1},j_r,\cdots,j_2,j_1)$. Note that $j_{i+r}=j_i$ $(1\leq i\leq r)$.
In this proof, we use the notation
\[ \textbf{Y}:=(Y_{1,j_{r}},\cdots,Y_{1,j_{1}},Y_{2,j_r},\cdots,Y_{2,j_2},Y_{2,j_1}), \]
for variables instead of $(t_1,\cdots,t_{2r})\in (\mathbb{C}^{\times})^{2r}$.

We define a map
$\phi:H\times(\mathbb{C}^{\times})^{2r}\rightarrow
H\times(\mathbb{C}^{\times})^{2r}$,
\begin{equation*}
\phi(a;\textbf{Y})=(\Phi_H(a;\textbf{Y});\Phi_{1,j_{r}}(\textbf{Y}),\cdots,\Phi_{1,j_{1}}(\textbf{Y}),
\Phi_{2,j_r}(\textbf{Y}),\cdots,\Phi_{2,j_2}(\textbf{Y}),\Phi_{2,j_1}(\textbf{Y})),
\end{equation*} 
as
\begin{equation}\label{mbasea} 
 \Phi_H(a;\textbf{Y}):=a\cdot \prod^{r}_{i=1}\prod^{2}_{j=1}\al_i^{\vee}(Y_{j,i}), 
\end{equation}
and for $l\in\{1,2,\cdots,r \}$, 

\begin{equation}\label{mbase0} 
\Phi_{1,l}(\textbf{Y}):=
\begin{cases}
\frac{(Y_{1,l-1}Y_{2,l-1})(Y_{1,l+1}Y_{2,l+1})}{Y_{1,l}Y_{2,l}^{2}} & {\rm if}\ l\ {\rm is\ even}, \\
\frac{(Y_{2,l-1})(Y_{2,l+1})}{Y_{1,l}Y_{2,l}^{2}} & {\rm if}\ l\ {\rm is\ odd},
\end{cases}
\end{equation}
\begin{equation}\label{mbase01} 
\Phi_{2,l}(\textbf{Y}):=
\begin{cases}
\frac{(Y_{2,l-1})(Y_{2,l+1})}{Y_{2,l}} & {\rm if}\ l\ {\rm is\ even}, \\
\frac{1}{Y_{2,l}} & {\rm if}\ l\ {\rm is\ odd}.
\end{cases}
\end{equation}

Note that $\phi$ is a biregular isomorphism since we can construct the inverse map $\psi:H\times(\mathbb{C}^{\times})^{2r}\rightarrow
H\times(\mathbb{C}^{\times})^{2r}$, 
\[ \psi(a;\textbf{Y})=(\Psi_H(a;\textbf{Y});\Psi_{1,j_{r}}(\textbf{Y}),\cdots,\Psi_{1,j_{1}}(\textbf{Y}),\Psi_{2,j_r}(\textbf{Y}),\cdots,\Psi_{2,j_1}(\textbf{Y}))
\]
of $\phi$ as follows:
\[
\Psi_{1,l}(\textbf{Y}):=
\begin{cases}
(Y_{1,l-1}Y_{1,l}Y_{1,l+1}Y_{2,l-3}Y_{2,l-2}Y_{2,l+2}Y_{2,l+3})^{-1} & {\rm if}\ l\ {\rm is\ even}, \\
(Y_{1,l}Y_{2,l-2}Y_{2,l-1}Y_{2,l+1}Y_{2,l+2})^{-1} & {\rm if}\ l\ {\rm is\ odd},
\end{cases}
\]
\[
\Psi_{2,l}(\textbf{Y}):=
\begin{cases}
(Y_{2,l-1}Y_{2,l}Y_{2,l+1})^{-1} & {\rm if}\ l\ {\rm is\ even}, \\
\frac{1}{Y_{2,l}} & {\rm if}\ l\ {\rm is\ odd},
\end{cases}
\]
\[ \Psi_{H}(a;\textbf{Y}):=
a\cdot(\prod^{r}_{i=1}\prod^{2}_{j=1}\al_i^{\vee}(\Psi_{j,i}(\textbf{Y})))^{-1}.
\]
Then, the map $\psi$ is the inverse map of $\phi$.

Let us prove
\[ \ovl{x}^G_{\textbf{i}}(a;\textbf{Y})=(x^G_{\textbf{i}}\circ\phi)(a;\textbf{Y}), \]
which implies that $\ovl{x}^G_{\textbf{i}}:H\times(\mathbb{C}^{\times})^{2r}\rightarrow G^{e,c^2}$ is a biregular isomorphism by Theorem \ref{fp}. First, it is known that for $1\leq i,\ j\leq r$ and $s,\ t\in \mathbb{C}^{\times}$,
\begin{equation}\label{base2}
\al_j^{\vee}(s)x_{i}(t)=\begin{cases}
	x_{i}(s^{2}t)\al_i^{\vee}(s) & {\rm if}\ i=j, \\
	x_{i}(s^{-1}t)\al_j^{\vee}(s) & {\rm if}\ |i-j|=1, \\
	x_{i}(t)\al_j^{\vee}(s) & {\rm otherwise}.
\end{cases}
\end{equation}
On the other hand, it follows from the definition (\ref{xgdef}) of $x^G_{\textbf{i}}$ and $(\ref{mbasea})$ that
\begin{multline*}(x^G_{\textbf{i}}\circ\phi)(a;\textbf{Y})
=a\cdot \left(\prod^{r}_{i=1}\prod^{2}_{s=1}\al_i^{\vee}(Y_{s,i})\right)
 \times x_{j_{r}}(\Phi_{1,j_{r}}(\textbf{Y}))
\cdots x_{j_{1}}(\Phi_{1,j_{1}}(\textbf{Y}))\\
\times x_{j_r}(\Phi_{2,j_r}(\textbf{Y}))
\cdots x_{j_2}(\Phi_{2,j_2}(\textbf{Y}))x_{j_1}(\Phi_{2,j_1}(\textbf{Y})).
\end{multline*}

For each even $l$ $(1\leq l\leq r)$, we can move 
\[
\al_{1}^{\vee}(Y_{1,1})\al_{3}^{\vee}(Y_{1,3})\cdots\al_{l-3}^{\vee}(Y_{1,l-3})\al_{l-1}^{\vee}(Y_{1,l-1})\prod^{r}_{i=l}\al_{i}^{\vee}(Y_{1,i})
\prod^{r}_{i=1}\al_{i}^{\vee}(Y_{2,i})
\]
to the right of
$x_{l}(\Phi_{1,l}(\textbf{Y}))$ by using the relations (\ref{base2}):
\begin{eqnarray*}
& &\hspace{-20pt}\al_{1}^{\vee}(Y_{1,1})\al_{3}^{\vee}(Y_{1,3})\cdots\al_{l-3}^{\vee}(Y_{1,l-3})\al_{l-1}^{\vee}(Y_{1,l-1})\left(\prod^{r}_{i=l}\al_{i}^{\vee}(Y_{1,i})
\prod^{r}_{i=1}\al_{i}^{\vee}(Y_{2,i})\right)
x_{l}(\Phi_{1,l}(\textbf{Y}))\\
&=&x_{l}(\Phi_{1,l}(\textbf{Y})\frac{Y_{1,l}^2Y_{2,l}^2}{Y_{1,l-1}Y_{2,l-1}Y_{1,l+1}Y_{2,l+1}})
\al_{1}^{\vee}(Y_{1,1})\cdots\al_{l-1}^{\vee}(Y_{1,l-1})\prod^{r}_{i=l}\al_{i}^{\vee}(Y_{1,i})
\prod^{r}_{i=1}\al_{i}^{\vee}(Y_{2,i})\\
&=&x_{l}(Y_{1,l})
\al_{1}^{\vee}(Y_{1,1})\al_{3}^{\vee}(Y_{1,3})\cdots\al_{l-3}^{\vee}(Y_{1,l-3})\al_{l-1}^{\vee}(Y_{1,l-1})\prod^{r}_{i=l}\al_{i}^{\vee}(Y_{1,i})
\prod^{r}_{i=1}\al_{i}^{\vee}(Y_{2,i}).
\end{eqnarray*}
Similarly, we can also move $\al_{1}^{\vee}(Y_{2,1})\al_{3}^{\vee}(Y_{2,3})\cdots\al_{l-1}^{\vee}(Y_{2,l-1})\prod^{r}_{i=l}\al_{i}^{\vee}(Y_{2,i})$ to the right of $x_{l}(\Phi_{2,l}(\textbf{Y}))$:
\begin{multline*}
\al_{1}^{\vee}(Y_{2,1})\al_{3}^{\vee}(Y_{2,3})\cdots\al_{l-3}^{\vee}(Y_{2,l-3})\al_{l-1}^{\vee}(Y_{2,l-1})\left(\prod^{r}_{i=l}\al_{i}^{\vee}(Y_{2,i})\right)x_{l}(\Phi_{2,l}(\textbf{Y}))\\
=x_{l}(Y_{2,l})\al_{1}^{\vee}(Y_{2,1})\al_{3}^{\vee}(Y_{2,3})\cdots\al_{l-3}^{\vee}(Y_{2,l-3})\al_{l-1}^{\vee}(Y_{2,l-1})\prod^{r}_{i=l}\al_{i}^{\vee}(Y_{2,i}).
\end{multline*}
For odd $l$, we obtain
\begin{multline*}
\al_{l}^{\vee}(Y_{1,l})\al_{l+2}^{\vee}(Y_{1,l+2})\cdots \al_{j_1}^{\vee}(Y_{1,j_1})\al_{1}^{\vee}(Y_{2,1})\al_{2}^{\vee}(Y_{2,2})\cdots \al_{r}^{\vee}(Y_{2,r})x_{l}(\Phi_{1,l}(\textbf{Y}))\\
=x_{l}(Y_{1,l})\al_{l}^{\vee}(Y_{1,l})\al_{l+2}^{\vee}(Y_{1,l+2})\cdots \al_{j_1}^{\vee}(Y_{1,j_1})\al_{1}^{\vee}(Y_{2,1})\al_{2}^{\vee}(Y_{2,2})\cdots \al_{r}^{\vee}(Y_{2,r}),
\end{multline*}
\[\al_{l}^{\vee}(Y_{2,l})\al_{l+2}^{\vee}(Y_{2,l+2})\cdots \al_{j_1}^{\vee}(Y_{2,j_1})x_{l}(\Phi_{2,l}(\textbf{Y}))
=x_{l}(Y_{2,l})\al_{l}^{\vee}(Y_{2,l})\al_{l+2}^{\vee}(Y_{2,l+2})\cdots \al_{j_1}^{\vee}(Y_{2,j_1}).\]

Thus, we get
\begin{multline*} 
(x^G_{\textbf{i}}\circ\phi)(a;\textbf{Y})=
a\cdot x_{j_r}(Y_{1,j_{r}})\al^{\vee}_{j_r}(Y_{1,j_{r}})\cdots x_{j_{1}}(Y_{1,j_{1}})\al^{\vee}_{j_{1}}(Y_{1,j_{1}})\\
 x_{j_r}(Y_{2,j_r})\al^{\vee}_{j_r}(Y_{2,j_r})\cdots
x_{j_2}(Y_{2,j_2})\al^{\vee}_{j_2}(Y_{2,j_2}) x_{j_1}(Y_{2,j_1})\al^{\vee}_{j_1}(Y_{2,j_1})=\ovl{x}^G_{\textbf{i}}(a;\textbf{Y}).
\end{multline*}
\qed

\section{Cluster algebras and generalized minors}\label{CluSect}
Following \cite{A-F-Z,F-Z,FZ2,M-M-A}, we review the definitions of cluster algebras and their generators called cluster variables. It is known that the coordinate rings of double Bruhat cells have cluster algebra structures, and generalized minors are their cluster variables \cite{GY}. We will refer to a relation between cluster variables on double Bruhat cells and crystal bases in Sect.\ref{gmc}.

We set $[1,l]:=\{1,2,\cdots,l\}$ and $[-1,-l]:=\{-1,-2,\cdots,-l\}$ for $l\in \mathbb{Z}_{>0}$. For $n,m\in \mathbb{Z}_{>0}$, let $x_1, \cdots ,x_n,x_{n+1}, \cdots
,x_{n+m}$ be commuting variables and $\cF:=\mathbb{C}(x_{1}, \cdots ,x_{n},x_{n+1},\cdots,x_{n+m})$ 
be the field of rational functions.

\subsection{Cluster algebras of geometric type}

In this subsection, we recall the definitions of cluster algebras. Let $\tilde{B}=(b_{ij})_{1\leq i\leq n+m,\ 1\leq j \leq n}$ be an $(n+m)\times
n$ integer matrix. The {\it principal part} $B$ of $\tilde{B}$ is obtained from $\tilde{B}$ by deleting the last $m$ rows. For $\tilde{B}$ and $k\in [1,n]$, the new $(n+m)\times n$ integer matrix $\mu_k(\tilde{B})=(b'_{ij})$ is defined by
\[b_{ij}':=
\begin{cases}
	-b_{ij} & {\rm if}\ i=k\ {\rm or}\ j=k, \\
	b_{ij}+\frac{|b_{ik}|b_{kj}+b_{ik}|b_{kj}|}{2} & {\rm otherwise}.
\end{cases}
\]
One calls $\mu_k(\tilde{B})$ the {\it matrix mutation} in direction $k$ of $\tilde{B}$. If there exists a positive 
integer diagonal matrix $D$ such that $DB$ is skew symmetric, we say $B$ is {\it skew symmetrizable}. Then we also say $\tilde{B}$ is skew symmetrizable. It is easily verified that if $\tilde{B}$ is skew symmetrizable then $\mu_k(\tilde{B})$ is also skew symmetrizable${\cite[Proposition3.6]{M-M-A}}$. We can also verify that $\mu_k\mu_k(\tilde{B})=\tilde{B}$. Define $\textbf{x}:=(x_1,\cdots,x_{n+m})$ and we call the pair $(\textbf{x}, \tilde{B})$ {\it initial seed}. For $1\leq k\leq n$, a new cluster variable $x_k'$ is defined by the following {\it exchange relation}.
\begin{equation}\label{exrel} x_k x_k' = 
\prod_{1\leq i \leq n+m,\ b_{ik}>0} x_i^{b_{ik}}
+\prod_{1\leq i \leq n+m,\ b_{ik}<0} x_i^{-b_{ik}}. \end{equation}
Let $\mu_k(\textbf{x})$ be the set of variables obtained from $\textbf{x}$ by replacing $x_k$ by $x'_k$. Ones call the pair $(\mu_k(\textbf{x}), \mu_k(\tilde{B}))$ the {\it mutation} in direction $k$ of the seed $(\textbf{x}, \tilde{B})$ and denote by $\mu_k((\textbf{x}, \tilde{B}))$.

Now, we can repeat this process of mutation and obtain a set of seeds inductively. Hence, each seed consists of an $(n+m)$-tuple of variables and a matrix. Ones call this $(n+m)$-tuple and matrix {\it cluster} and {\it exchange matrix} respectively. Variables in cluster is called {\it cluster variables}. In particular, the variables $x_{n+1},\cdots,x_{n+m}$ are called {\it frozen cluster variables}.

\begin{defn}${\cite{F-Z, M-M-A}}$\label{clusterdef}
Let $\tilde{B}$ be an integer matrix whose principal part is skew symmetrizable and $\Sigma=(\textbf{x},\tilde{B})$ a seed. We set ${\mathbb A}:={\mathbb Z}[x_{n+1}^{\pm1}, \cdots ,x_{n+m}^{\pm1}]$. The cluster algebra (of geometric type) $\cA=\cA(\Sigma)$ over $\mathbb A$ associated with
seed $\Sigma$ is defined as the ${\mathbb A}$-subalgebra of $\cF$ generated by all cluster variables in all seeds which can be obtained from $\Sigma$ by sequences of mutations.
\end{defn}

\subsection{Cluster algebra $\cA({\rm \bf{i}})$}\label{cAi}

In the rest of this section, let $G={\rm SL}_{r+1}(\mathbb{C})$ be the complex simple algebraic group of type ${\rm A}_r$. Let $\ge:={\rm Lie}(G)$ and $A=(a_{i,j})$ be its Cartan matrix. In Definition \ref{redworddef}, we define a reduced word ${\rm \bf{i}}=(j_{n},\cdots,j_2,j_{1})$ for an element $v$ of Weyl group $W$. In this subsection, we define a cluster algebra $\cA({\rm \bf{i}})$, which obtained from ${\rm \bf{i}}$. It satisfies that $\cA({\rm \bf{i}})\otimes \mathbb{C}$ is isomorphic to the coordinate ring $\mathbb{C}[G^{e,v}]$ of the double Bruhat cell \cite{A-F-Z}. Let $j_k$ $(k\in[1,n])$ be the $k$-th index of ${\rm \bf{i}}$ from the right. Let us add $r$ additional entries $j_{-r},\cdots,j_{-1}$ at the beginning of ${\rm \bf{i}}$. by setting $j_{-t}=-t$ $(t\in [1,r])$.

For $l\in[1,n]$, we denote by $l^-$ the largest index $k\in[1,n]$ such that $k<l$ and $j_l=j_k$. If $l\in[-1,-r]$, let $l^-$ be the largest index $k\in[1,n]$ such that $|j_l|=|j_k|$. For example, if $[-1,-3]\cup{\rm \bf{i}}=(-3,-2,-1,2,1,3,2,1,3)$ then, $(-1)^-=5$, $(-2)^-=6$, $(-3)^-=4$, $4^-=1$, $5^-=2$, $6^-=3$, and $3^-$, $2^-$, $1^-$ are not defined. We define a set e({\rm \bf{i}}) as 
\[ e({\rm \bf{i}}):=\{k\in[1,n]| k^-\ \hbox{is well-defined}\}.  \]
Following \cite{A-F-Z}, we define a quiver $\Gamma_{{\rm \bf{i}}}$ as follows. The vertices of $\Gamma_{{\rm \bf{i}}}$ are the variables $x_k$ ($k\in[-1,-r]\cup[1,n]$). For two vertices $x_k$ $(k\in [-1,-r]\cup[1,n])$ and $x_l$ $(l\in e({\rm \bf{i}}))$ with either $l<k$ or $k\in[-1,-r]$, there exists an arrow $x_k\rightarrow x_l$ (resp. $x_l\rightarrow x_k$) if and only if $l=k^-$ (resp. $l^-<k^-<l$ and $a_{|j_k|,|j_l|}<0$). Next, let us define a matrix $\tilde{B}=\tilde{B}({\rm \bf{i}})$. 

\begin{defn}
Let $\tilde{B}({\rm \bf{i}})$ be an integer matrix with rows labelled by all the indices in $[-1,-r]\cup [1,n]$ and columns labelled by all the indices in $e({\rm \bf{i}})$. For $k\in[-1,-r]\cup [1,n]$ and $l\in e({\rm \bf{i}})$, an entry $b_{kl}$ of $\tilde{B}({\rm \bf{i}})$ is determined as follows: If there exists an arrow $x_k\rightarrow x_l$ (resp. $x_l\rightarrow x_k$) in $\Gamma_{{\rm \bf{i}}}$, then
\[
b_{kl}:=\begin{cases}
		1\ ({\rm resp.}\ -1)& {\rm if}\ |j_k|=|j_l|, \\
		-a_{|j_k||j_l|}\ ({\rm resp.}\ a_{|j_k||j_l|})& {\rm if}\ |j_k|\neq|j_l|.
	\end{cases}
\]
Unless there exist any arrow between $k$ and $l$, we set $b_{kl}=0$. The principal part $B({\rm \bf{i}})$ of $\tilde{B}({\rm \bf{i}})$ is the submatrix $(b_{i,j})_{i,j\in e({\rm \bf{i}})}$. We also define $\Sigma_{{\rm \bf{i}}}:=(\textbf{x},\tilde{B}({\rm \bf{i}}))$.
\end{defn}

\begin{prop}\label{propss}${\cite{A-F-Z}}$
$\tilde{B}({\rm \bf{i}})$ is skew symmetrizable. 
\end{prop}

In general, for an $(m+l)$-tuple of variables $\textbf{y}=(y_i)^{m+l}_{i=1}$ and an $(m+l)\times l$-skew symmetrizable matrix $\tilde{B}=(b_{i,j})$ with $m,\ l>0$, let $\Gamma((\textbf{y},\tilde{B}))$ be a quiver whose vertices are $y_1,\cdots,y_{m+l}$, and whose arrows are determined as follows: For $i\in[1,m+l]$ and $j\in[1,l]$, there exists an arrow $y_i\rightarrow y_j$ (resp. $y_j\rightarrow y_i$) if $b_{i,j}>0$
(resp. $b_{i,j}<0$). We can easily check that $\Gamma((\textbf{x},\tilde{B}({\rm \bf{i}})))=\Gamma_{{\rm \bf{i}}}$.
\begin{lem}$\cite{M-M-A}$\label{mutgamlem}
Let $(\textbf{y},\tilde{B})$ be a seed, where $\textbf{y}=(y_i)^{m+l}_{i=1}$ and $\tilde{B}=(b_{i,j})$ is an $(m+l)\times l$-skew symmetrizable matrix with $|b_{i,j}|=1$ or $0$. For $k\in[1,l]$, the quiver $\Gamma((\mu_k(\textbf{y}),\mu_k(\tilde{B})))$ has vertices $y_1,\cdots,y_k',\cdots,y_{m+l}$ and arrows determined as follows:
\begin{enumerate}
\item[$(1)$]If $y_i\rightarrow y_k$ (resp. $y_k\rightarrow y_i$) in $\Gamma((\textbf{y},\tilde{B}))$ then $y_k'\rightarrow y_i$ (resp. $y_i\rightarrow y_k'$) in $\Gamma((\mu_k(\textbf{y}),\mu_k(\tilde{B})))$.
\item[$(2)$]We suppose that there exist arrows $y_i\rightarrow y_k$ and $y_k\rightarrow y_j$ in $\Gamma((\textbf{y},\tilde{B}))$ with either $i\in[1,l]$ or $j\in[1,l]$. If there exists an arrow $y_j\rightarrow y_i$ (resp. $y_i$ and $y_j$ are not connected) in $\Gamma((\textbf{y},\tilde{B}))$, then $y_i$ and $y_j$ are not connected (resp. $y_i\rightarrow y_j$) in $\Gamma((\mu_k(\textbf{y}),\mu_k(\tilde{B})))$.
\item[$(3)$]The rest of the arrows are the same as the one of $\Gamma((\textbf{y},\tilde{B}))$.
\end{enumerate}
\end{lem}

All the skew symmetrizable matrices $B=(b_{i,j})$ appearing in this article satisfy $|b_{i,j}|=1$ or $0$. Thus, we can use the above Lemma repeatedly.

\begin{ex}\label{gammaex}
Let us consider the case $G={\rm SL}_{5}(\mathbb{C})$ and ${\rm \bf{i}}=(2,4,1,3,2,4,1,3)$. The quiver $\Gamma_{{\rm \bf{i}}}$ is described as
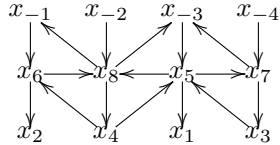
\[
\begin{xy}
(80,98) *{x_{-4}}="7",
(80,90)*{x_7}="3",
(80,82)*{x_3}="-4",
(70,98) *{x_{-3}}="5",
(70,90)*{x_5}="1",
(70,82)*{x_1}="-3",
(60,98) *{x_{-2}}="8",
(60,90)*{x_8}="4",
(60,82)*{x_4}="-2",
(50,98) *{x_{-1}}="6",
(50,90)*{x_6}="2",
(50,82)*{x_2}="-1",
(35,98) *{}="12",
\ar@{->} "7";"3"
\ar@{->} "5";"1"
\ar@{->} "8";"4"
\ar@{->} "6";"2"
\ar@{->} "3";"-4"
\ar@{->} "1";"-3"
\ar@{->} "4";"-2"
\ar@{->} "2";"-1"
\ar@{->} "1";"3"
\ar@{->} "1";"4"
\ar@{->} "2";"4"
\ar@{->} "3";"5"
\ar@{->} "4";"5"
\ar@{->} "4";"6"
\ar@{->} "-4";"1"
\ar@{->} "-2";"1"
\ar@{->} "-2";"2"
\end{xy} \] 
In general, let us consider the case $G={\rm SL}_{r+1}(\mathbb{C})$ and the sequence ${\rm \bf{i}}$ in $(\ref{redwords2})$. Let $\lfloor \ \rfloor $ denote the Gaussian symbol and $j_k$ be the $k$-th index of ${\rm \bf{i}}$ from the right:
${\rm \bf{i}}=(j_{2r},\cdots,j_{r+1},j_r\cdots,j_1)$. For $k$ $(1\leq k\leq \lfloor \frac{r+1}{2} \rfloor)$, vertices and arrows around the vertex $x_k$ in the quiver $\Gamma_{{\rm \bf{i}}}$ are described as

\begin{equation}\label{gammaex-1}
\begin{xy}
(95,90)*{\cdots}="emp1",
(85,98)*{x_{-j_k-1}}="7",
(85,90) *{x_{r+\lfloor \frac{r}{2} \rfloor + k}}="3",
(85,82)*{x_{\lfloor \frac{r}{2} \rfloor +k}}="-4",
(60,98)*{x_{-j_k}}="r+k",
(60,90) *{x_{r+k}}="k",
(60,82)*{x_k}="-j_k",
(35,98)*{x_{-j_k+1}}="8",
(35,90) *{x_{r+\lfloor \frac{r}{2} \rfloor +k+1}}="4",
(35,82)*{x_{\lfloor \frac{r}{2} \rfloor +k+1}}="-2",
(10,98)*{x_{-j_k+2}}="6",
(10,90) *{x_{r+k+1}}="2",
(10,82)*{x_{k+1}}="-1",
(0,90)*{\cdots}="emp",
\ar@{->} "7";"3"
\ar@{->} "r+k";"k"
\ar@{->} "8";"4"
\ar@{->} "6";"2"
\ar@{->} "3";"-4"
\ar@{->} "k";"-j_k"
\ar@{->} "4";"-2"
\ar@{->} "2";"-1"
\ar@{->} "k";"3"
\ar@{->} "k";"4"
\ar@{->} "2";"4"
\ar@{->} "3";"r+k"
\ar@{->} "4";"r+k"
\ar@{->} "4";"6"
\ar@{->} "-4";"k"
\ar@{->} "-2";"k"
\ar@{->} "-2";"2"
\end{xy} \end{equation}
For $k$ $(\lfloor \frac{r+1}{2} \rfloor< k \leq r)$, it is described as

\begin{equation}\label{gammaex-2}
\begin{xy}
(95,90)*{\cdots}="emp1",
(85,98) *{x_{-j_k-2}}="7",
(85,90)*{x_{r+k-1}}="3",
(85,82)*{x_{k-1}}="-4",
(60,98) *{x_{-j_k-1}}="r+k",
(60,90)*{x_{r+k-\lfloor \frac{r}{2} \rfloor-1}}="k",
(60,82)*{x_{k-\lfloor \frac{r}{2} \rfloor-1}}="-j_k",
(35,98) *{x_{-j_k}}="8",
(35,90)*{x_{r+k}}="4",
(35,82)*{x_k}="-2",
(10,98) *{x_{-j_k+1}}="6",
(10,90)*{x_{r+k-\lfloor \frac{r}{2} \rfloor}}="2",
(10,82)*{x_{k-\lfloor \frac{r}{2} \rfloor}}="-1",
(0,90)*{\cdots}="emp",
\ar@{->} "7";"3"
\ar@{->} "r+k";"k"
\ar@{->} "8";"4"
\ar@{->} "6";"2"
\ar@{->} "3";"-4"
\ar@{->} "k";"-j_k"
\ar@{->} "4";"-2"
\ar@{->} "2";"-1"
\ar@{->} "k";"3"
\ar@{->} "k";"4"
\ar@{->} "2";"4"
\ar@{->} "3";"r+k"
\ar@{->} "4";"r+k"
\ar@{->} "4";"6"
\ar@{->} "-4";"k"
\ar@{->} "-2";"k"
\ar@{->} "-2";"2"
\end{xy} \end{equation} 
\end{ex}

\begin{defn}$\cite{A-F-Z}$
By Definition $\ref{clusterdef}$ and Proposition \ref{propss}, we can construct the cluster algebra.
We denote this cluster algebra by $\cA({\rm \bf{i}})$.
\end{defn}

\subsection{Generalized minors}\label{bilingen}

Set $\cA({\rm \bf{i}})_{\mathbb{C}}:=\cA({\rm \bf{i}})\otimes \mathbb{C}$. It is known that the coordinate ring $\mathbb{C}[G^{e,v}]$ of the double Bruhat cell is isomorphic to $\cA({\rm \bf{i}})_{\mathbb{C}}$ (Theorem \ref{clmainthm}). To describe this isomorphism explicitly, we need generalized minors.  

We set $G_0:=N_-HN$, and let $x=[x]_-[x]_0[x]_+$ with $[x]_-\in N_-$, $[x]_0\in H$, $[x]_+\in N$ be the corresponding decomposition. 

\begin{defn}
For $i\in[1,r]$ and $w\in W$, the {\it generalized minor} $\Delta_{\Lambda_i,w\Lambda_i}$ is a regular function on $G$ whose restriction to the open set $G_0\ovl{w}^{-1}$ is given by $\Delta_{\Lambda_i,w\Lambda_i}(x)=([x \ovl{w} ]_0)^{\Lambda_i}$. Here, $\Lambda_i$ is the $i$-th  fundamental weight. 
\end{defn}
The generalized minor $\Delta_{\Lambda_i,w\Lambda_i}$ depends on $w\Lambda_i$ and does not depend on $w$. By definition, for $a\in H$, $x\in G$, $w\in W$, $i,j\in I$ and $t\in\mathbb{C}$,
\begin{equation}\label{genbasic}
\Delta_{\Lambda_i,w\Lambda_i}(ax)=a^{\Lambda_i}\Delta_{\Lambda_i,w\Lambda_i}(x),\ 
\Delta_{\Lambda_i,\Lambda_i}(x x_{j}(t))=\Delta_{\Lambda_i,\Lambda_i}(x),
\end{equation}
where $x_{j}(t)\in N$ is the one in (\ref{xiyidef}).

\subsection{Cluster algebras on double Bruhat cells}

For a reduced expression $v=s_{j_n}s_{j_{n-1}}\cdots s_{j_1}\in W$ and $k\in [1,n]$, we set
\begin{equation}\label{inc}
v_{> k}=v_{> k}({\rm \bf{i}}):=s_{j_{1}}s_{j_{2}}\cdots s_{j_{n-k}}.
\end{equation}
For $k \in [1,n]$, we define $\Delta(k;{\rm \bf{i}})(x):=\Delta_{\Lambda_{j_k},v_{>n-k+1}\Lambda_{j_k}}(x)$, and for $k \in [-1,-r]$, $\Delta(k;{\rm \bf{i}})(x):=\Delta_{\Lambda_{|k|},v^{-1}\Lambda_{|k|}}(x)$. 

Finally, we set $F({\rm \bf{i}}):=\{ \Delta(k;{\rm \bf{i}})(x)|k \in [-1,-r]\cup[1,n] \}$. It is known that the set $F({\rm \bf{i}})$ is an algebraically independent generating set for the field of rational functions $\mathbb{C}(G^{e,v})$ \cite[Theorem 1.12]{F-Z}. Then, we have the following.
\begin{thm}\label{clmainthm}${\cite{A-F-Z, GLS, GY}}$
The isomorphism of fields $\varphi :\cF \rightarrow \mathbb{C}(G^{e,v})$ defined by $\varphi (x_k)=\Delta(k;{\rm \bf{i}})\ (k \in [-1,-r]\cup [1,l(v)] )$ restricts to an isomorphism of algebras $\cA({\rm \bf{i}})_{\mathbb{C}}\rightarrow \mathbb{C}[G^{e,v}]$.
\end{thm}

\begin{ex}\label{clmainthmex}
Let $v=c^2$ be the square of Coxeter element such that whose reduced word ${\rm \bf{i}}=(j_{2r},\cdots,j_{r+1},j_r\cdots,j_1)$ is written as in $(\ref{redwords2})$. Then for $k\in[1,r]$, the correspondence of the initial cluster variables are as follows:
\begin{eqnarray*}
 x_{-j_k}&\mapsto& \Delta_{\Lambda_{j_k},c^{-2}\Lambda_{j_k}}=
\Delta_{\Lambda_{j_k},s_{j_1}s_{j_2}\cdots s_{j_{2r}}\Lambda_{j_k}}
=\Delta_{\Lambda_{j_k},s_{j_1}s_{j_2}\cdots s_{j_{r+k}}\Lambda_{j_k}}
=\Delta_{\Lambda_{j_k},c^2_{>r-k}\Lambda_{j_k}}
,\\
 x_{r+k}&\mapsto & \Delta_{\Lambda_{j_k},c^2_{>r-k+1}\Lambda_{j_k}}=
\Delta_{\Lambda_{j_k},s_{j_1}s_{j_2}\cdots j_{r+k-1}\Lambda_{j_k}}\\
& &=\Delta_{\Lambda_{j_k},s_{j_1}s_{j_2}\cdots s_{j_k}\Lambda_{j_k}}
=\Delta_{\Lambda_{j_k},c^2_{>2r-k}\Lambda_{j_k}},\\
 x_k &\mapsto& \Delta_{\Lambda_{j_k},c^2_{>2r-k+1}\Lambda_{j_k}}=
\Delta_{\Lambda_{j_k},s_{j_1}s_{j_2}\cdots s_{j_{k-1}}\Lambda_{j_k}}
=\Delta_{\Lambda_{j_k},\Lambda_{j_k}}.
\end{eqnarray*}

\end{ex}

\subsection{Finite type}

Let $\mathcal{S}$ be the set of seeds of a cluster algebra $\mathcal{A}$. If $\mathcal{S}$ is finite, then $\mathcal{A}$ is said to be of {\it finite type}. In this subsection, we shall review cluster algebras of finite type \cite{FZ3}.

Let $B=(b_{ij})$ be an integer square matrix. The {\it Cartan counter part} of $B$ is a generalized Cartan matrix $A=A(B)=(a_{i,j})$ defined as follows:
\[ a_{i,j}=
\begin{cases}
2 & {\rm if}\ i=j,\\
-|b_{i,j}| & {\rm if}\ i\neq j.
\end{cases}
\]

\begin{thm}\label{finthm1}$\cite{FZ3}$
The cluster algebra $\cA$ is of finite type if and only if there exists a seed $\Sigma=(\textbf{y},\tilde{B})$ such that $\cA= \cA(\Sigma)$ and $A(B)$ is a Cartan matrix of finite type, where $B$ is the principal part of $\tilde{B}$.
\end{thm}

By this theorem, we can define the {\it type} of each cluster algebra of finite type mirroring the Cartan-Killing classification.

Let $\Phi$ be the root system associated with a Cartan matrix, with the set of simple roots $\Pi=\{\al_i|\ i\in I\}$ and the set of positive roots $\Phi_{>0}$. The set of {\it almost positive roots} is defined as $\Phi_{\geq -1}:=\Phi_{>0}\cup -\Pi$.

\begin{thm}\label{finthm2}$\cite{FZ3}$
\begin{enumerate}
\item For a cluster algebra $\cA$ of finite type, the number of the cluster variables included in $\cA$ is equal to $|\Phi_{\geq -1}|$, where $\Phi$ is the root system associated with the Cartan matrix of the same type as $\cA$.
\item Let $c\in W$ be a Coxeter element of $G$ whose length $l(c)$ satisfies 
$l(c^2)=2 l(c)=2${\rm rank}$(G)$. Then the coordinate ring $\mathbb{C}[G^{e,c^2}]$ has a structure of cluster algebra of finite type under the isomorphism in Theorem \ref{clmainthm}, and its type is the Cartan-Killing type of $G$.
\end{enumerate}
\end{thm}

\section{Additive categorifications of cluster algebras}
We fix an element $v\in W$ and set $n:=l(v)$. In this section, we set $G={\rm SL}_{r+1}(\mathbb{C})$ and review the additive categorifications of the coordinate rings $\mathbb{C}[L^{e,v}]$. We refer to \cite{ASS, LD, GLS}. 

\subsection{Preprojective algebras and Category $\mathcal{C}_v$}\label{ppalg-cat-sub}

Let $Q=(Q_0,Q_1,s,t)$ be a Dynkin quiver of type A and 
\[ \Lambda=\mathbb{C}\ovl{Q}/(\mathscr{C}) \]
the associated preprojective algebra. Here $\ovl{Q}$ is the double quiver of $Q$:
\[
\xymatrix{
1 \ar@<0.5ex>[r] & 2 \ar@<0.5ex>[r] \ar@<0.5ex>[l] & 3 \ar@<0.5ex>[r] \ar@<0.5ex>[l] & \cdots \ar@<0.5ex>[r] \ar@<0.5ex>[l]  & r \ar@<0.5ex>[l] 
}
\] 
and $\mathbb{C}\ovl{Q}$ is its path algebra, and $(\mathscr{C})$ is the ideal generated by
\[ \mathscr{C}=\sum_{a\in Q_1}(a^{*}a-a a^{*}). \]

Let $\widehat{I}_1,\cdots,\widehat{I}_r$ be the indecomposable injective $\Lambda$-modules which have the simple socle isomorphic to $S_1,\cdots,S_r$, respectively, where $S_i$ is the $1$-dimensional simple $\Lambda$-module which corresponds to the vertex $i$ in $Q$. The module $\widehat{I}_j$ is described as follows:
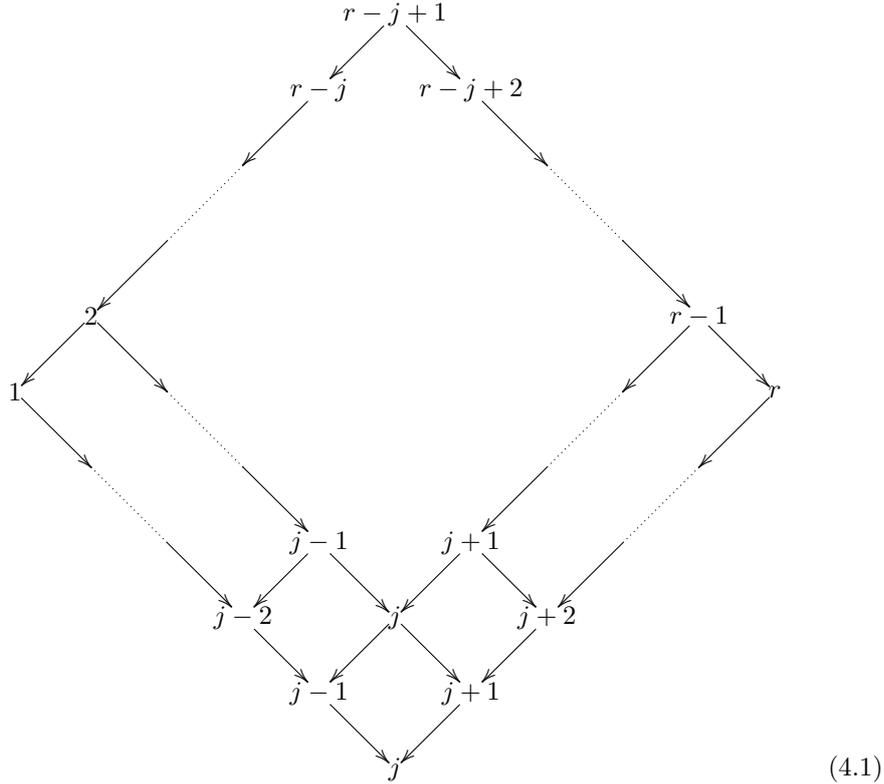
\begin{equation}\label{inj-mod}
\begin{xy}
(50,100)*{r-j+1}="r-j+1",
(40,90)*{r-j}="5",(60,90)*{r-j+2}="r-4",
(30,80)*{}="4",(70,80)*{}="r-3",
(20,70)*{}="33",(80,70)*{}="r-2r-2",
(10,60)*{2}="2",(90,60)*{r-1}="r-1",
(0,50)*{1}="1",(20,50)*{}="emp0002",(80,50)*{}="emp002",(100,50)*{r}="r",
(10,40)*{}="emp1",(30,40)*{}="emp0001",(70,40)*{}="emp001",(90,40)*{}="emp01",
(20,30)*{}="emp",(40,30)*{j-1}="emp000",(60,30)*{j+1}="emp00",(80,30)*{}="emp0",
(30,20)*{j-2}="j-2",(50,20)*{j}="jj",(70,20)*{j+2}="j+2",
(40,10)*{j-1}="j-1",(60,10)*{j+1}="j+1",
(50,0)*{j}="j",
\ar@{->} "r-j+1";"5"
\ar@{->} "r-j+1";"r-4"
\ar@{->} "5";"4"
\ar@{->} "r-4";"r-3"
\ar@{.} "4";"33"
\ar@{.} "r-3";"r-2r-2"
\ar@{->} "33";"2"
\ar@{->} "r-2r-2";"r-1"
\ar@{->} "2";"1"
\ar@{->} "r-1";"r"
\ar@{->} "1";"emp1"
\ar@{->} "r";"emp01"
\ar@{->} "r-1";"emp002"
\ar@{.} "emp002";"emp001"
\ar@{->} "emp001";"emp00"
\ar@{->} "2";"emp0002"%
\ar@{.} "emp0002";"emp0001"
\ar@{->} "emp0001";"emp000"%
\ar@{.} "emp01";"emp0"
\ar@{.} "emp1";"emp"
\ar@{->} "emp000";"jj"
\ar@{->} "emp000";"j-2"
\ar@{->} "emp00";"jj"
\ar@{->} "emp00";"j+2"
\ar@{->} "emp0";"j+2"
\ar@{->} "emp";"j-2"
\ar@{->} "j+2";"j+1"
\ar@{->} "j-2";"j-1"
\ar@{->} "jj";"j-1"
\ar@{->} "jj";"j+1"
\ar@{->} "j+1";"j"
\ar@{->} "j-1";"j"
\end{xy}
\end{equation}

In (\ref{inj-mod}), each vertex $k$ $(1\leq k\leq r)$ means a basis of $\widehat{I}_j$, and each arrow $k\rightarrow k+1$ (resp. $k\rightarrow k-1$) means the action of the edge $k\rightarrow k+1$ (resp. $k\rightarrow k-1$) $\in \Lambda$ on the basis $k$. The vertex $e_k\in \Lambda$ acts on each basis $k'$ as 
\[e_k.k'=
\begin{cases}
k' & {\rm if}\ k=k',\\
0 & {\rm if}\ k\neq k'. 
\end{cases}
\]
For example, the vertex $e_j\in\Lambda$ acts on the basis $j$ located at the bottom of (\ref{inj-mod}) identically, and all other paths act trivially. Thus, $1$-dimensional submodule generated by this basis $j$ is isomorphic to the simple module $S_j$.

Let mod$(\Lambda)$ be the category of finite dimensional $\Lambda$-modules. Note that though in \cite{GLS} the category ${\rm nil}(\Lambda)$ is treated, we consider the category ${\rm mod}(\Lambda)$ instead of ${\rm nil}(\Lambda)$ since ${\rm mod}(\Lambda)={\rm nil}(\Lambda)$ holds in our setting. For $j\in Q_0$ and $\Lambda$-module $X$ in mod$(\Lambda)$, let ${\rm soc}_{j}(X)$ be the sum of all submodules $U$ of $X$ with $U\cong S_j$. For a sequence $(i_1,\cdots,i_t)$ $(i_1,i_2,\cdots,i_t\in Q_0)$, there exists a unique chain
\[ 0=X_0\subset X_1\subset\cdots X_t\subset X \]
of submodules such that $X_p/X_{p-1}={\rm soc}_{i_p}(X/X_{p-1})$ $(p=1,2,\cdots,t)$. We define ${\rm soc}_{(i_1,\cdots,i_t)}(X):=X_t$.

Let $v\in W$ and $\textbf{i}=(j_n,\cdots,j_1)$ be its reduced word. Without loss of generality, we assume that for each $j\in[1,r]$, there exist some $k\in[1,n]$ such that $j_k=j$. The $\Lambda$-modules $V_{k}=V_{\textbf{i},k}$ $(k=1,2,\cdots,n)\in {\rm mod}(\Lambda)$ are defined as
\[ V_{k}:=V_{\textbf{i},k}={\rm soc}_{(j_k,\cdots,j_1)}(\widehat{I}_{j_k} ). \]
Let $V_{\bf{i}}:=\bigoplus^n_{k=1} V_{k}$ and $\mathcal{C}_{\bf{i}}$ be the full subcategory of mod$(\Lambda)$
whose objects are factor modules of direct sums of finitely many copies of $V_{\bf{i}}$. For $j\in [1,r]$, let $m_j:={\rm max}\{1\leq m\leq n | j_m=j\}$ and $I_{\textbf{i},j}:=V_{\textbf{i},m_j}$. We also set $I_{\textbf{i}}:=I_{\textbf{i},1}\oplus\cdots\oplus I_{\textbf{i},r}$. The category $\mathcal{C}_{\bf{i}}$ and $I_{\textbf{i}}$ depend on only $v$, and do not depend on the choice of reduced word $\bf{i}$. Thus, we define
\[ \mathcal{C}_{v}:=\mathcal{C}_{\bf{i}},\ \ I_{v}:=I_{\textbf{i}}. \]

A $\Lambda$-module $C$ in $\mathcal{C}_{v}$ is called $\mathcal{C}_{v}$-{\it projective} (resp. $\mathcal{C}_{v}$-{\it injective}) if ${\rm Ext}^1_{\Lm}(C,X)=0$ (resp. ${\rm Ext}^1_{\Lm}(X,C)=0$) for all $X\in \mathcal{C}_{v}$. If $C$ is $\mathcal{C}_{v}$-projective and $\mathcal{C}_{v}$-injective, $C$ is said to be $\mathcal{C}_{v}$-{\it projective-injective}.

\begin{thm}$\cite{BIRS, GLS}$
The category $\mathcal{C}_{v}$ has $r$ indecomposable $\mathcal{C}_{v}$-projective-injective modules, which are the indecomposable direct summand of $ I_{v}$.
\end{thm}

\begin{ex}\label{initialex1}
Let {\rm \bf{i}} be the sequence in $(\ref{redwords2})$, and calculate $V_k=V_{{\rm \bf{i}},k}$ $(1\leq k\leq 2r)$. Let $j_k$ be the $k$ th index of {\rm \bf{i}} from the right, that is, ${\rm \bf{i}}=(j_{2r},\cdots,j_{r+1},j_r,\cdots,j_1)$. For example, $j_1=r-1$ if $r$ is even and $j_1=r$ if $r$ is odd. Note that $j_l=j_{l+r}$ $(1\leq l\leq r)$.

First, to calculate $V_1={\rm soc}_{(j_1)}(\hat{I}_{j_1})$, we consider the chain $0=X_0\subset X_1\subset \hat{I}_{j_1}$
such that $X_1/X_0=X_1={\rm soc}_{j_1}(\hat{I}_{j_1})=S_{j_1}$. By definition, we get $V_1=X_1=S_{j_1}$.

Next, for $2\leq k\leq \lfloor \frac{r+1}{2} \rfloor$, to calculate $V_k={\rm soc}_{(j_k,j_{k-1},\cdots,j_1)}(\hat{I}_{j_k})$, we consider the chain 
\[ 0=X_0\subset X_1\subset X_2\subset\cdots\subset X_k \subset\hat{I}_{j_k} \]
such that $X_1={\rm soc}_{j_k}(\hat{I}_{j_k})=S_{j_k}$, $X_2/X_1={\rm soc}_{j_{k-1}}(\hat{I}_{j_k}/X_1)$, $X_3/X_2={\rm soc}_{j_{k-2}}(\hat{I}_{j_k}/X_2),\cdots$, $X_{k}/X_{k-1}={\rm soc}_{j_1}(\hat{I}_{j_k}/X_{k-1})$. By $(\ref{inj-mod})$, the module $\hat{I}_{j_k}/S_{j_k}$ has simple submodules isomorphic to $S_{j_k-1}$ and $S_{j_k+1}$. Since $\hat{I}_{j_k}/S_{j_k}$ has no simple submodules isomorphic to $S_{j_{k-1}},\ S_{j_{k-2}}, \cdots,\ S_{j_1}$, we have $X_1=X_2=\cdots=X_k$ and then
\begin{equation}\label{iniex1-1}
V_k=X_k=S_{j_k}.
\end{equation}
Next, for $\lfloor \frac{r+1}{2} \rfloor +1\leq k\leq r$, we consider the chain 
\[ 0=X_0\subset X_1\subset X_2\subset\cdots\subset X_k \subset\hat{I}_{j_k} \]
such that $X_1={\rm soc}_{j_k}(\hat{I}_{j_k})=S_{j_k}$, $X_2/X_1={\rm soc}_{j_{k-1}}(\hat{I}_{j_k}/X_1),\cdots$. In the same way as in $(\ref{iniex1-1})$, we get
$X_1=X_2=\cdots=X_{\lfloor \frac{r}{2} \rfloor}=S_{j_k}$. And $X_{\lfloor \frac{r}{2} \rfloor+1}/X_{\lfloor \frac{r}{2} \rfloor}={\rm soc}_{j_{k-\lfloor \frac{r}{2} \rfloor}}(\hat{I}_{j_k}/X_{\lfloor \frac{r}{2} \rfloor})=S_{j_k-1}$. So the module $X_{\lfloor \frac{r}{2} \rfloor+1}$ is described as
\[
\begin{xy}
(65,103) *{j_k-1}="2r-2k+1",
(75,97)*{j_k}="2r-2k+2",
(55,103) *{}="12",
\ar@{->} "2r-2k+1";"2r-2k+2"
\end{xy}
\]
Similarly, we obtain $X_{\lfloor \frac{r}{2} \rfloor+2}/X_{\lfloor \frac{r}{2} \rfloor+1}={\rm soc}_{j_{k-\lfloor \frac{r}{2} \rfloor-1}}(\hat{I}_{j_k}/X_{\lfloor \frac{r}{2} \rfloor+1})=S_{j_k+1}$. In the same way as in $(\ref{iniex1-1})$, we have $V_k=X_{k}=X_{k-1}=\cdots=X_{\lfloor \frac{r}{2} \rfloor+2}$. Thus, the module $V_k$ is described as
\begin{equation}\label{iniex1-2}
\begin{xy}
(65,103)*{j_k-1}="2r-2k+1",(85,103)*{j_k+1}="2r-2k+3",
(75,97)*{j_k}="2r-2k+2",
(55,103) *{}="12",
\ar@{->} "2r-2k+1";"2r-2k+2"
\ar@{->} "2r-2k+3";"2r-2k+2"
\end{xy}
\end{equation}
Next, for $r+1\leq k\leq \lfloor \frac{r+1}{2} \rfloor +r$, we consider the chain 
\[ 0=X_0\subset X_1\subset X_2\subset\cdots\subset X_k \subset\hat{I}_{j_k} \]
such that $X_1={\rm soc}_{j_k}(\hat{I}_{j_k})=S_{j_k}$, $X_2/X_1={\rm soc}_{j_{k-1}}(\hat{I}_{j_k}/X_1),\cdots$. In the same way as in $(\ref{iniex1-1})$, we get
$X_1=X_2=\cdots=X_{\lfloor \frac{r+1}{2} \rfloor-1}=S_{j_k}$. And $X_{\lfloor \frac{r+1}{2} \rfloor}/X_{\lfloor \frac{r+1}{2} \rfloor-1}={\rm soc}_{j_{k-\lfloor \frac{r+1}{2} \rfloor+1}}(\hat{I}_{j_k}/X_{\lfloor \frac{r+1}{2} \rfloor-1})=S_{j_k-1}$, where we set $S_j:=0$ for $j\leq 0$. We also get
$X_{\lfloor \frac{r+1}{2} \rfloor+1}/X_{\lfloor \frac{r+1}{2} \rfloor}={\rm soc}_{j_{k-\lfloor \frac{r+1}{2} \rfloor}}(\hat{I}_{j_k}/X_{\lfloor \frac{r+1}{2} \rfloor})=S_{j_k+1}$, and
\[ X_{\lfloor \frac{r+1}{2} \rfloor+1}=X_{\lfloor \frac{r+1}{2} \rfloor+2}=\cdots =X_{r-1}. \]
We also obtain $X_r/X_{r-1}={\rm soc}_{j_{k-r+1}}(\hat{I}_{j_k}/X_{r-1})=S_{j_k-2}$, $X_{r+1}/X_{r}=S_{j_k}$, $X_{r+2}/X_{r+1}=S_{j_k+2}$ and
$X_{r+2}=X_{r+3}=\cdots=X_{k}$.
Therefore, the module $X_k=V_k$ is described as
\begin{equation}\label{iniex1-3}
\begin{xy}
(50,113)*{j_k-2}="3r-2k-1",(70,113)*{j_k}="3r-2k+1a",(90,113)*{j_k+2}="3r-2k+3",
(60,107)*{j_k-1}="3r-2k",(80,107)*{j_k+1}="3r-2k+2",
(70,101)*{j_k}="3r-2k+1",
(35,107)*{}="12",
\ar@{->} "3r-2k";"3r-2k+1"
\ar@{->} "3r-2k+2";"3r-2k+1"
\ar@{->} "3r-2k-1";"3r-2k"
\ar@{->} "3r-2k+1a";"3r-2k"
\ar@{->} "3r-2k+1a";"3r-2k+2"
\ar@{->} "3r-2k+3";"3r-2k+2"
\end{xy}
\end{equation}
Finally, for $\lfloor \frac{r+1}{2} \rfloor +r+1\leq k\leq 2r$, we can verify that the module $V_{k}$ is described as:
\begin{equation}\label{iniex1-4}
\begin{xy}
(35,123)*{_{j_k-3}}="4r-2k-1",(55,123)*{_{j_k-1}}="4r-2k+1b",(75,123)*{_{j_k+1}}="4r-2k+3b",
(95,123)*{_{j_k+3}}="4r-2k+5",
(45,118)*{_{j_k-2}}="4r-2k",(65,118)*{_{j_k}}="4r-2k+2a",(85,118)*{_{j_k+2}}="4r-2k+4",
(55,113)*{_{j_k-1}}="4r-2k+1",(75,113)*{_{j_k+1}}="4r-2k+3",
(65,108)*{_{j_k}}="4r-2k+2",
(35,108)*{}="12",
\ar@{->} "4r-2k+1";"4r-2k+2"
\ar@{->} "4r-2k+3";"4r-2k+2"
\ar@{->} "4r-2k";"4r-2k+1"
\ar@{->} "4r-2k+2a";"4r-2k+1"
\ar@{->} "4r-2k+2a";"4r-2k+3"
\ar@{->} "4r-2k+4";"4r-2k+3"
\ar@{->} "4r-2k-1";"4r-2k"
\ar@{->} "4r-2k+1b";"4r-2k"
\ar@{->} "4r-2k+1b";"4r-2k+2a"
\ar@{->} "4r-2k+3b";"4r-2k+2a"
\ar@{->} "4r-2k+3b";"4r-2k+4"
\ar@{->} "4r-2k+5";"4r-2k+4"
\end{xy}
\end{equation}
by the same argument as in $(\ref{iniex1-1})$, $(\ref{iniex1-2})$ and $(\ref{iniex1-3})$. In this case, we have $I_{c^2}=I_{\bf{{\rm i}}}=V_{r+1}\oplus\cdots\oplus V_{2r}$.
\end{ex}

\begin{rem}\label{diarem}
When we see the quiver
\[
\begin{xy}
(35,123)*{_{j-3}}="4r-2k-1",(55,123)*{_{j-1}}="4r-2k+1b",(75,123)*{_{j+1}}="4r-2k+3b",
(95,123)*{_{j+3}}="4r-2k+5",
(45,118)*{_{j-2}}="4r-2k",(65,118)*{_{j}}="4r-2k+2a",(85,118)*{_{j+2}}="4r-2k+4",
(55,113)*{_{j-1}}="4r-2k+1",(75,113)*{_{j+1}}="4r-2k+3",
(65,108)*{_{j}}="4r-2k+2",
(35,108)*{}="12",
\ar@{->} "4r-2k+1";"4r-2k+2"
\ar@{->} "4r-2k+3";"4r-2k+2"
\ar@{->} "4r-2k";"4r-2k+1"
\ar@{->} "4r-2k+2a";"4r-2k+1"
\ar@{->} "4r-2k+2a";"4r-2k+3"
\ar@{->} "4r-2k+4";"4r-2k+3"
\ar@{->} "4r-2k-1";"4r-2k"
\ar@{->} "4r-2k+1b";"4r-2k"
\ar@{->} "4r-2k+1b";"4r-2k+2a"
\ar@{->} "4r-2k+3b";"4r-2k+2a"
\ar@{->} "4r-2k+3b";"4r-2k+4"
\ar@{->} "4r-2k+5";"4r-2k+4"
\end{xy}
\]
or its subquiver, if $j=1$, $2$ or $3$, we understand it means
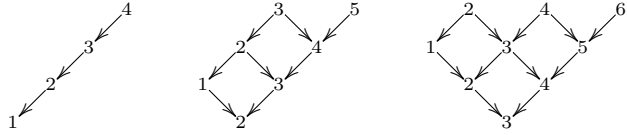
\[
\begin{xy}
(15,95)*{_{4}}="4r-2k+5",
(10,90)*{_{3}}="4r-2k+4",
(5,85)*{_{2}}="4r-2k+3",
(0,80)*{_{1}}="4r-2k+2",
(35,95)*{_{3}}="3a",(45,95)*{_{5}}="5",
(30,90)*{_{2}}="2a",(40,90)*{_{4}}="4",
(25,85)*{_{1}}="1",(35,85)*{_{3}}="3",
(30,80)*{_{2}}="2",
(60,95)*{_{2}}="2bb",(70,95)*{_{4}}="4bb",(80,95)*{_{6}}="6b",
(55,90)*{_{1}}="1b",(65,90)*{_{3}}="3bb",(75,90)*{_{5}}="5bb",
(60,85)*{_{2}}="2b",(70,85)*{_{4}}="4b",
(65,80)*{_{3}}="3b",
\ar@{->} "4r-2k+3";"4r-2k+2"
\ar@{->} "4r-2k+4";"4r-2k+3"
\ar@{->} "4r-2k+5";"4r-2k+4"
\ar@{->} "3a";"4"
\ar@{->} "3a";"2a"
\ar@{->} "5";"4"
\ar@{->} "2a";"1"
\ar@{->} "2a";"3"
\ar@{->} "4";"3"
\ar@{->} "1";"2"
\ar@{->} "3";"2"
\ar@{->} "2bb";"1b"
\ar@{->} "2bb";"3bb"
\ar@{->} "4bb";"3bb"
\ar@{->} "4bb";"5bb"
\ar@{->} "6b";"5bb"
\ar@{->} "1b";"2b"
\ar@{->} "3bb";"2b"
\ar@{->} "3bb";"4b"
\ar@{->} "5bb";"4b"
\ar@{->} "2b";"3b"
\ar@{->} "4b";"3b"
\end{xy}
\]
respectively. Similarly, if $j=r$, $r-1$ or $r-2$, we understand it means
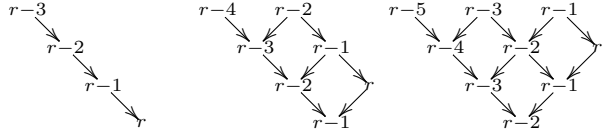
\[
\begin{xy}
(0,95)*{_{r-3}}="4r-2k+5",
(5,90)*{_{r-2}}="4r-2k+4",
(10,85)*{_{r-1}}="4r-2k+3",
(15,80)*{_{r}}="4r-2k+2",
(25,95)*{_{r-4}}="3a",(35,95)*{_{r-2}}="5",
(30,90)*{_{r-3}}="2a",(40,90)*{_{r-1}}="4",
(35,85)*{_{r-2}}="1",(45,85)*{_{r}}="3",
(40,80)*{_{r-1}}="2",
(50,95)*{_{r-5}}="2bb",(60,95)*{_{r-3}}="4bb",(70,95)*{_{r-1}}="6b",
(55,90)*{_{r-4}}="1b",(65,90)*{_{r-2}}="3bb",(75,90)*{_{r}}="5bb",
(60,85)*{_{r-3}}="2b",(70,85)*{_{r-1}}="4b",
(65,80)*{_{r-2}}="3b",
\ar@{->} "4r-2k+3";"4r-2k+2"
\ar@{->} "4r-2k+4";"4r-2k+3"
\ar@{->} "4r-2k+5";"4r-2k+4"
\ar@{->} "5";"2a"
\ar@{->} "3a";"2a"
\ar@{->} "5";"4"
\ar@{->} "2a";"1"
\ar@{->} "4";"1"
\ar@{->} "4";"3"
\ar@{->} "1";"2"
\ar@{->} "3";"2"
\ar@{->} "2bb";"1b"
\ar@{->} "4bb";"1b"
\ar@{->} "4bb";"3bb"
\ar@{->} "6b";"3bb"
\ar@{->} "6b";"5bb"
\ar@{->} "1b";"2b"
\ar@{->} "3bb";"2b"
\ar@{->} "3bb";"4b"
\ar@{->} "5bb";"4b"
\ar@{->} "2b";"3b"
\ar@{->} "4b";"3b"
\end{xy}
\]
\end{rem}

\subsection{Mutation}

For a $\Lambda$-module $T$ in mod$(\Lambda)$, let add$(T)$ denote the subcategory of mod$(\Lambda)$ whose objects are all $\Lambda$-modules which are isomorphic to finite direct sums of direct
summands of $T$.

\begin{defn}$\cite{ASS, LD, GLS}$
\begin{enumerate}
\item A $\Lambda$-module $T$ is {\it rigid} if ${\rm Ext}^1_{\Lambda}(T,T)=0$. 
\item For a rigid module $T$ in $\mathcal{C}_v$, we say $T$ is a $\mathcal{C}_v$-{\it cluster-tilting} module if ${\rm Ext}^1_{\Lambda}(T,X)=0$ with $X\in \mathcal{C}_v$ implies $X\in {\rm add}(T)$.
\item A $\Lambda$-module $T$ is said to be {\it basic}, if it is decomposed to a direct sum
of pairwise non-isomorphic indecomposable modules.
\item Let $T$, $X$ and $Y\in {\rm mod}(\Lambda)$. A morphism $f\in{\rm Hom}_{\Lambda}(X,Y)$ (resp. $f\in{\rm Hom}_{\Lambda}(Y,X)$) is said to be a {\it left} (resp. {\it right}) {\rm add}(T)-{\it approximation} of $X$ if $Y\in{\rm add}(T)$ and for an arbitrary $Y'\in {\rm add}(T)$ and $f'\in{\rm Hom}_{\Lambda}(X,Y')$ (resp. $f'\in{\rm Hom}_{\Lambda}(Y',X)$), there exists $g\in {\rm Hom}_{\Lambda}(Y,Y')$ (resp. $g\in {\rm Hom}_{\Lambda}(Y',Y)$) and $f'=g\circ f$ (resp. $f'=f\circ g$).
\item For $V$, $W\in {\rm mod}(\Lambda)$, a morphism $f\in{\rm Hom}_{\Lambda}(V,W)$ is said to be {\it left} (resp. {\it right}) {\it minimal} if every endomorphism $g\in{\rm End}_{\Lambda}(W)$ (resp. $g\in{\rm End}_{\Lambda}(V)$) such that $g\circ f=f$ (resp. $f\circ g=f$) is an isomorphism.
\end{enumerate}
\end{defn}

\begin{prop}\label{exseqprop}$\cite{LD,GLS,GLS2}$
Let $T=T_1\oplus T_2\oplus \cdots\oplus T_n$ be a basic $\mathcal{C}_{v}$-cluster-tilting object. We suppose that the $\{T_i\}_{i=1,2,\cdots,n}$ are indecomposable summands of $T$ and $T_{n-r+1}, \cdots, T_n$ are the $\mathcal{C}_{v}$-projective-injective modules. Then for $k \in\{1,2,\cdots, n-r\}$, there exists a short exact sequence 
\begin{equation}\label{exseqdef}
 0\rightarrow T_k\overset{f}{\rightarrow} \ovl{T_k}\overset{g}{\rightarrow} T^*_k\rightarrow 0 
\end{equation}
such that
\begin{enumerate}
\item $f$ is a left minimal left {\rm add}$(T/T_k)$-approximation,
\item $g$ is a right minimal right {\rm add}$(T/T_k)$-approximation,
\item $T^*_k$ is indecomposable,
\item $T^*_k\notin {\rm add}(T)$,
\item $T/T_k\oplus T^*_k$ is basic $\mathcal{C}_{v}$-cluster-tilting.
\end{enumerate}
\end{prop}

\begin{defn}$\cite{LD, GLS}$
In the setting of the previous proposition, the mutation $\mu_{T_k}(T)$ of $T$ in direction $T_k$ is defined as
\begin{equation}\label{mudiretk}
 \mu_{T_k}(T):=T/T_k\oplus T^*_k. 
\end{equation}
We call the short exact sequence (\ref{exseqdef}) in Proposition \ref{exseqprop} the {\it exchange sequence} associated to the direct summand $T_k$ of $T$.
\end{defn}

For a basic module $T=T_1\oplus\cdots\oplus T_n$ in $\mathcal{C}_{v}$, let $\Gamma_{T}$ be the quiver of ${\rm End}_{\Lambda}(T)^{{\rm op}}$, that is, ${\rm End}_{\Lambda}(T)^{{\rm op}}\cong\mathbb{C}\Gamma_{T}/(R)$ with an admissible ideal $(R)$ \cite{ASS}. Setting
\[ {\rm Rad}(T_i,T_j)=
\begin{cases}
{\rm Hom}_{\Lambda}(T_i,T_j) & {\rm if}\ i\neq j,\\
\{{\rm nilpotent\ elements\ of}\ {\rm End}_{\Lambda}(T_i)\} & {\rm if}\ i= j,
\end{cases}
\]
we have the following:
\begin{lem}\label{titj}$\cite{ASS, LD}$
The quiver $\Gamma_{T}$ has $n$ vertices indexed by $\{1,2,\cdots,n\}$, and for $1\leq i, j\leq n$, the number of arrows $j\rightarrow i$ is equal to the dimension of the space 
\[ \frac{{\rm Rad}(T_i,T_j)}{\sum^{n}_{k=1}{\rm Rad}(T_k,T_j)\circ {\rm Rad}(T_i,T_k)}. \]
\end{lem}

\begin{defn}
Let $T=T_1\oplus\cdots\oplus T_n$ be a basic module in $\mathcal{C}_{v}$. For $i,j\in[1,n]$ and a non-zero homomorphism $f\in {\rm Hom}_{\Lambda}(T_i,T_j)$, it is said that $f$ is {\it factorizable} in the direct summands of $T$ if it belongs to $\sum^{n}_{k=1}{\rm Rad}(T_k,T_j)\circ {\rm Rad}(T_i,T_k)$.
\end{defn}

Let $B(\Gamma_{T})=(b_{i,j})$ denote $n\times (n-r)$-matrix defined by
\[ b_{i,j}=({\rm number\ of\ arrows}\ j\rightarrow i\ {\rm in}\ \Gamma_{T})-({\rm number\ of\ arrows}\ i\rightarrow j\ {\rm in}\ \Gamma_{T}).
\]

For $\textbf{i}=(j_n,\cdots,j_1)\in Q^n_0$, we define a quiver $\overline{\Gamma}_{\textbf{i}}$ as follows: We use the notation $k^-$ in \ref{cAi}. We also denote $k^+$ the smallest index $l$ such that $k<l$ and $|j_l|=|j_k|$ if it exists. If it does not exist, we set $k^+=n+1$. The
vertices of $\overline{\Gamma}_{\textbf{i}}$ are $1,2,\cdots,n$. For two vertices $k,l\in[1,n]$ with $l<k$, there exists an arrow $k\rightarrow l$ (resp. $l\rightarrow k$) if and only if $l=k^-$ (resp. $k<l^+\leq k^+$ and $a_{i_k,i_l}<0$).

\begin{thm}\label{GLSthm}$\cite{BIRS, GLS, GLS2}$
Let $n=l(v)$ and ${\rm \bf{i}}=(j_n,\cdots,j_1)$ be a reduced word of $v$.
\begin{enumerate}
\item The module $V_{{\rm \bf{i}}}$ defined in \ref{ppalg-cat-sub}
 is a basic $\mathcal{C}_{v}$-cluster-tilting object and $\Gamma_{V_{{\rm \bf{i}}}}=\overline{\Gamma}_{{\rm \bf{i}}}$.
\item Let $T=T_1\oplus T_2\oplus \cdots\oplus T_n$ be a basic $\mathcal{C}_{v}$-cluster-tilting object. For $1\leq k \leq n-r$, we have $B(\Gamma_{\mu_{T_k}(T)})=\mu_k(B(\Gamma_{T}))$.
\item For a basic $\mathcal{C}_{v}$-cluster-tilting object $T=T_1\oplus T_2\oplus \cdots\oplus T_n$ and $1\leq k \leq n-r$, the exchange sequence associated to the direct summand $T_k$ of $T$ is
\[ 0\rightarrow T_k \rightarrow \bigoplus_{i\rightarrow k\ {\rm in}\ \Gamma_T}T_i\rightarrow T_k^*\rightarrow 0. \]
\end{enumerate}
\end{thm}

\begin{ex}

Let ${\rm \bf{i}}$ be the reduced word in $(\ref{redwords2})$. By Theorem \ref{GLSthm} (i), for $1\leq k\leq\lfloor \frac{r+1}{2} \rfloor$, the quiver
$\Gamma_{V_{{\rm \bf{i}}}}$ is described as
\begin{equation}\label{gammavex-1}
\begin{xy}
(100,86)*{\cdots}="emp1",
(85,90) *{r+\lfloor \frac{r}{2} \rfloor + k}="3",
(85,82)*{\lfloor \frac{r}{2} \rfloor +k}="-4",
(60,90) *{r+k}="k",
(60,82)*{k}="-j_k",
(35,90) *{r+\lfloor \frac{r}{2} \rfloor +k+1}="4",
(35,82)*{\lfloor \frac{r}{2} \rfloor +k+1}="-2",
(10,90) *{r+k+1}="2",
(10,82)*{k+1}="-1",
(0,86)*{\cdots}="emp",
\ar@{->} "3";"-4"
\ar@{->} "k";"-j_k"
\ar@{->} "-j_k";"-4"
\ar@{->} "-j_k";"-2"
\ar@{->} "-1";"-2"
\ar@{->} "4";"-2"
\ar@{->} "2";"-1"
\ar@{->} "k";"3"
\ar@{->} "k";"4"
\ar@{->} "2";"4"
\ar@{->} "-4";"k"
\ar@{->} "-2";"k"
\ar@{->} "-2";"2"
\end{xy} 
\end{equation}

\end{ex}

\begin{ex}\label{initialex2}
In the setting of Example \ref{initialex1}, let us consider the mutation of $V_{\rm \bf{i}}$ in direction $V_k$ $(1\leq k\leq r)$. Let us constitute the exchange sequence
\[ 0\rightarrow V_k \rightarrow \ovl{V_k}\rightarrow V^*_k\rightarrow 0 \]
associated to the direct summand $V_k$ of $V_{{\rm \bf{i}}}$.

For $1\leq k\leq \lfloor \frac{r+1}{2} \rfloor$, recall that $V_k=S_{j_k}$. In $\{V_i|\ 1\leq i\leq 2r,\ i\neq k\}$, the module $V_{r+k}$ has the simple socle isomorphic to $S_{j_k}$ and the others do not so since their simple socles are $S_l$ $(l\neq j_k)$ by $(\ref{iniex1-2})$, $(\ref{iniex1-3})$, $(\ref{iniex1-4})$. The module $V_{r+k}$ is described as
\[ \begin{xy}
(50,113)*{j_k-2}="3r-2k-1",(70,113)*{j_k}="3r-2k+1a",(90,113)*{j_k+2}="3r-2k+3",
(60,107)*{j_k-1}="3r-2k",(80,107)*{j_k+1}="3r-2k+2",
(70,101)*{j_k}="3r-2k+1",
(35,107)*{}="12",
\ar@{->} "3r-2k";"3r-2k+1"
\ar@{->} "3r-2k+2";"3r-2k+1"
\ar@{->} "3r-2k-1";"3r-2k"
\ar@{->} "3r-2k+1a";"3r-2k"
\ar@{->} "3r-2k+1a";"3r-2k+2"
\ar@{->} "3r-2k+3";"3r-2k+2"
\end{xy} \]
and bottom $j_k$ means a basis generating the simple socle isomorphic to $S_{j_k}$ $((\ref{inj-mod}),\ (\ref{iniex1-3}))$. Hence, there exists an injective homomorphism $V_k\rightarrow V_{r+k}$, and its image is the simple socle. By the above argument, we have {\rm Hom}$(V_k,V_{r+k})\cong\mathbb{C}$ and {\rm Rad}$(V_k,V_{t})=\{0\}$ for $t\neq r+k$. We get $\ovl{V_k}=V_{r+k}$ by Lemma \ref{titj} and Theorem \ref{GLSthm} (iii), which yields 
that $V^*_k$ $(1\leq k< \lfloor \frac{r+1}{2} \rfloor)$ and $V^*_{\lfloor \frac{r+1}{2} \rfloor}$ are described as
\begin{equation}\label{mutex1}
\begin{xy}
(60,113)*{_{j_k-2}}="3r-2k-1",(70,113)*{_{j_k}}="3r-2k+1a",(80,113)*{_{j_k+2}}="3r-2k+3",
(65,106)*{_{j_k-1}}="3r-2k",(75,106)*{_{j_k+1}}="3r-2k+2",
(87,106)*{,}=",",
(100,113)*{_3}="3",
(95,106)*{_2}="2",
\ar@{->} "3r-2k-1";"3r-2k"
\ar@{->} "3r-2k+1a";"3r-2k"
\ar@{->} "3r-2k+1a";"3r-2k+2"
\ar@{->} "3r-2k+3";"3r-2k+2"
\ar@{->} "3";"2"
\end{xy} 
\end{equation}
respectively.

Next, for $\lfloor \frac{r+1}{2} \rfloor+1\leq k \leq r$, the module $V_k$ is given as $(\ref{iniex1-2})$. 
The module $V_{r+k}$ is described as $(\ref{iniex1-4})$ and it has the submodule isomorphic to $V_k$, which is generated by the basis $j_k-1$, $j_k$ and $j_k+1$ lower one in $(\ref{iniex1-4})$. Let $c_{j_k-1}$, $c_{j_k}$ and $c_{j_k+1}$ denote these three bases. Thus, there exists an injective homomorphism $V_k\rightarrow V_{r+k}$. Since $V_k$ has the simple quotients isomorphic to $S_{j_k-1}$, $S_{j_k+1}$, there exist surjective homomorphisms $V_k\rightarrow V_{k-\lfloor \frac{r}{2} \rfloor}=S_{j_k-1}$ and $V_k\rightarrow V_{k-\lfloor \frac{r}{2} \rfloor-1}=S_{j_k+1}$ $($note that 
$j_{k-\lfloor \frac{r}{2} \rfloor}=j_k-1$ and $j_{k-\lfloor \frac{r}{2} \rfloor-1}=j_k+1)$. The modules $V_{r+k-\lfloor \frac{r}{2} \rfloor}$ and $V_{r+k-\lfloor \frac{r}{2} \rfloor-1}$ have the simple submodules isomorphic to $S_{j_k-1}$ and $S_{j_k+1}$ respectively. However, homomorphisms $V_{k}\rightarrow V_{r+k-\lfloor \frac{r}{2} \rfloor}$ and $V_{k}\rightarrow V_{r+k-\lfloor \frac{r}{2} \rfloor-1}$ are factorizable in the direct summands of $V_{{\rm \bf{i}}}$ since they are equal to the composite maps $V_k\rightarrow V_{k-\lfloor \frac{r}{2} \rfloor}\rightarrow V_{r+k-\lfloor \frac{r}{2} \rfloor}$ and $V_k\rightarrow V_{k-\lfloor \frac{r}{2} \rfloor-1}\rightarrow V_{r+k-\lfloor \frac{r}{2} \rfloor-1}$ respectively. Moreover, we see that {\rm Rad}$(V_k,V_t)=0$ for $t\neq r+k,\ k-\lfloor \frac{r}{2} \rfloor,\ k-\lfloor \frac{r}{2} \rfloor-1$ since $V_t$ does not have submodule isomorphic to $V_k$, $S_{j_k-1}$ and $S_{j_k+1}$. From this, the homomorphisms
$V_k\rightarrow V_{r+k}$, $V_k\rightarrow V_{k-\lfloor \frac{r}{2} \rfloor}$ and $V_k\rightarrow V_{k-\lfloor \frac{r}{2} \rfloor-1}$ are not factorizable in the direct summands of $V_{{\rm \bf{i}}}$. Therefore, the exchange sequence
associated to the direct summand $V_k$ of $V_{{\rm \bf{i}}}$ is
\[ 0\rightarrow V_k\rightarrow V_{r+k}\oplus S_{j_k-1}\oplus S_{j_k+1}\rightarrow V^*_k\rightarrow 0 \]
 by Lemma \ref{titj} and Theorem \ref{GLSthm} (iii). The image of the homomorphism $V_k\rightarrow V_{r+k}\oplus S_{j_k-1}\oplus S_{j_k+1}$ is $3$-dimensional and it can be explicitly written as $
\mathbb{C}(c_{j_k-1}+d_{j_k-1})\oplus\mathbb{C}(c_{j_k})\oplus\mathbb{C}(c_{j_k+1}+e_{j_k+1})$ with some non-zero elements $d_{j_k-1}\in S_{j_k-1}$, $e_{j_k+1}\in S_{j_k+1}$. By the above argument, the module $V^*_k=(V_{r+k}\oplus S_{j_k-1}\oplus S_{j_k+1})/(\mathbb{C}(c_{j_k-1}\oplus d_{j_k-1}\oplus e_{j_k+1})\oplus\mathbb{C}(c_{j_k}\oplus d_{j_k-1}\oplus e_{j_k+1})\oplus\mathbb{C}(c_{j_k+1}\oplus d_{j_k-1}\oplus e_{j_k+1}))$ is described as follows:
\begin{equation}\label{mutex2}
\begin{xy}
(35,123)*{j_k-3}="4r-2k-1",(55,123)*{j_k-1}="4r-2k+1b",(75,123)*{j_k+1}="4r-2k+3b",
(95,123)*{j_k+3}="4r-2k+5",
(45,115)*{j_k-2}="4r-2k",(65,115)*{j_k}="4r-2k+2a",(85,115)*{j_k+2}="4r-2k+4",
(55,107)*{j_k-1}="4r-2k+1",(75,107)*{j_k+1}="4r-2k+3",
(35,107)*{}="12",
\ar@{->} "4r-2k";"4r-2k+1"
\ar@{->} "4r-2k+2a";"4r-2k+1"
\ar@{->} "4r-2k+2a";"4r-2k+3"
\ar@{->} "4r-2k+4";"4r-2k+3"
\ar@{->} "4r-2k-1";"4r-2k"
\ar@{->} "4r-2k+1b";"4r-2k"
\ar@{->} "4r-2k+1b";"4r-2k+2a"
\ar@{->} "4r-2k+3b";"4r-2k+2a"
\ar@{->} "4r-2k+3b";"4r-2k+4"
\ar@{->} "4r-2k+5";"4r-2k+4"
\end{xy}
\end{equation}
\end{ex}

By the same way in this example, we have the following proposition:

\begin{prop}\label{remaindia}
The modules $(\mu_{V_r}\mu_{V_{\lfloor \frac{r+1}{2} \rfloor}}V_{{\rm \bf{i}}})_{r}$ and $(\mu_{V_{k-\lfloor \frac{r}{2} \rfloor-1}}\mu_{V_k}V_{{\rm \bf{i}}})_{k-\lfloor \frac{r}{2} \rfloor-1}$ $(\lfloor \frac{r+1}{2} \rfloor+2\leq k\leq r)$ are described as
\[
\begin{xy}
(35,95)*{_{3}}="3a",(45,95)*{_{5}}="5",
(30,90)*{_{2}}="2a",(40,90)*{_{4}}="4",
(35,85)*{_{3}}="3",
(50,85)*{,}=",",
(55,95)*{_{j_k-3}}="4r-2k-1",(65,95)*{_{j_k-1}}="4r-2k+1b",(75,95)*{_{j_k+1}}="4r-2k+3b",
(85,95)*{_{j_k+3}}="4r-2k+5",
(60,90)*{_{j_k-2}}="4r-2k",(70,90)*{_{j_k}}="4r-2k+2a",(80,90)*{_{j_k+2}}="4r-2k+4",
(65,85)*{_{j_k-1}}="4r-2k+1",
\ar@{->} "3a";"4"
\ar@{->} "3a";"2a"
\ar@{->} "5";"4"
\ar@{->} "2a";"3"
\ar@{->} "4";"3"
\ar@{->} "4r-2k";"4r-2k+1"
\ar@{->} "4r-2k+2a";"4r-2k+1"
\ar@{->} "4r-2k-1";"4r-2k"
\ar@{->} "4r-2k+1b";"4r-2k"
\ar@{->} "4r-2k+1b";"4r-2k+2a"
\ar@{->} "4r-2k+3b";"4r-2k+2a"
\ar@{->} "4r-2k+3b";"4r-2k+4"
\ar@{->} "4r-2k+5";"4r-2k+4"
\end{xy}
\]
respectively. Note that $j_{\lfloor \frac{r+1}{2} \rfloor}=1$.
\end{prop}

\subsection{Cluster algebra structure of $\mathbb{C}[L^{e,v}]$}

For a $\Lambda$-module $X$ and a sequence $\textbf{k}=(k_1,\cdots,k_s)$ $(k_t\in[1,r])$, let $\cF_{\textbf{k},X}$ denote the projective variety of composition series of $X$:
\[ 0=X_0\subset X_1\subset X_2\subset\cdots\subset X_{s}=X,  \]
such that each subfactor $X_t/X_{t-1}$ is isomorphic
to the simple $\Lambda$-module $S_{k_t}$ $(1\leq t\leq s)$. Recall that we set $x_i(t):={\rm exp}(te_i)$ in (\ref{xiyidef}).

\begin{prop}\label{dualsemi}$\cite{LD, GLS}$
For each $\Lambda$-module $X$ in ${\rm mod}(\Lambda)$, there exists a unique function $\varphi_X\in\mathbb{C}[N]$ such that for any sequence ${\rm \bf{i}}=(i_1,\cdots,i_k)$ $(1\leq i_1,\cdots,i_k\leq r)$,
\[ \varphi_X (x_{i_1}(t_1)x_{i_2}(t_2)\cdots x_{i_k}(t_k))=\sum_{\textbf{a}=(a_1,\cdots,a_k)\in (\mathbb{Z}_{\geq0})^k}
\chi_c(\mathcal{F}_{{\rm \bf{i}}^{\textbf{a}},X})\frac{t^{a_1}_1\cdots t^{a_k}_k}{a_1!\cdots a_k!}, \]
where $\chi_c$ is the Euler characteristic, and for $\textbf{a}=(a_1,a_2,\cdots,a_k)$, 
\[{\rm \bf{i}}^{\bf{a}}:=(\underbrace{i_1,\cdots,i_1}_{a_1},\underbrace{i_2,\cdots,i_2}_{a_2},\cdots,\underbrace{i_k,\cdots,i_k}_{a_k}).\]
\end{prop}

Note that we can write $x_{i_1}(t_1)x_{i_2}(t_2)\cdots x_{i_k}(t_k)=x^G_{\rm \bf{i}}(1;t_1,\cdots,t_k)$, where $1$ is the identity element of $H$ and $x^G_{\rm \bf{i}}$ is defined in (\ref{xgdef}).

For a $\Lambda$-module $X$ in ${\rm mod}(\Lambda)$ and ${\rm \bf{i}}=(i_1,\cdots,i_k)$, $\textbf{a}=(a_1,a_2,\cdots,a_k)\in(\mathbb{Z}_{\geq 0})^{k}$, let $\mathcal{F}_{{\rm \bf{i}},{\rm \bf{a}},X}$ be the projective variety of partial composition series of $X$
\[ 0=X_0\subset X_1\subset X_2\subset\cdots\subset X_{k}=X  \]
such that each subfactor $X_t/X_{t-1}$ is isomorphic to $S^{a_t}_{i_t}$ for all $1\leq t\leq k$. Then we have $\chi_c(\mathcal{F}_{{\rm \bf{i}}^{\rm \bf{a}},X})=\chi_c(\mathcal{F}_{{\rm \bf{i}},{\rm \bf{a}},X})
a_1! a_2! \cdots a_k!$ \cite{GLS}. Therefore, in the setting of Proposition \ref{dualsemi},
\begin{equation}\label{iatoia}
 \varphi_X (x_{i_1}(t_1)x_{i_2}(t_2)\cdots x_{i_k}(t_k))=\sum_{\textbf{a}=(a_1,\cdots,a_k)\in (\mathbb{Z}_{\geq0})^k}
\chi_c(\mathcal{F}_{{\rm \bf{i}},{\textbf{a}},X})t^{a_1}_1\cdots t^{a_k}_k. 
\end{equation}

\begin{ex}\label{initialex3}
In the setting of Example \ref{initialex1} and \ref{initialex2}, let us calculate $\varphi_{V_k}$ $(1\leq k\leq r)$ and $\varphi_{(\mu_k{V_{{\rm \bf{i}}}})_k}$ $(1\leq k\leq r)$. We set $\textbf{Y}:=(Y_{1,j_{r}},\cdots,Y_{1,j_{1}},Y_{2,j_r},\cdots,Y_{2,j_{2}},Y_{2,j_1})$.

For ${\rm \bf{i}}$ in $(\ref{redwords2})$, let us consider the variety of flags $\cF_{{\rm \bf{i}}^{{\rm \bf{a}}},V_k}$. Let $j_k$ be the $k$ th index of ${\rm \bf{i}}$ from the right. We write ${\rm \bf{a}}\in (\mathbb{Z}_{\geq0})^{2r}$ as follows: 
\[
{\rm \bf{a}}=(a_{1,j_{r}},\cdots,a_{1,j_{2}},a_{1,j_{1}},
a_{2,j_r},\cdots, a_{2,j_2},a_{2,j_1})  .
\]
By Example \ref{initialex1}, for $1\leq k\leq \lfloor \frac{r+1}{2} \rfloor$, since $V_k=S_{j_k}$, if $\cF_{{\rm \bf{i}}^{{\rm \bf{a}}},V_k}\neq \phi$ then ${\rm \bf{i}}^{{\rm \bf{a}}}=(j_k)$, which implies $a_{1,j_k}=1$ and other $a_{1,j},a_{2,j}$ are equal to $0$, or $a_{2,j_k}=1$ and other $a_{1,j},a_{2,j}$ are equal to $0$. In this case, $\cF_{{\rm \bf{i}}^{{\rm \bf{a}}},V_k}$ is a point $(=(0\subset S_{j_k}=V_k))$. Thus, Proposition \ref{dualsemi} means that
\[\varphi_{V_k} (x^G_{\rm \bf{i}}(1;\textbf{Y}))=Y_{1,j_k}+Y_{2,j_k}. \]

Next, for $\lfloor \frac{r+1}{2} \rfloor+1\leq k\leq r$, the module $V_k$ is described as $(\ref{iniex1-2})$. If $\cF_{{\rm \bf{i}}^{{\rm \bf{a}}},V_k}\neq \phi$ then ${\rm \bf{i}}^{{\rm \bf{a}}}=(j_k,j_k-1,j_k+1)$ or ${\rm \bf{i}}^{{\rm \bf{a}}}=(j_k,j_k+1,j_k-1)$, which implies 
\[a_{1,j_k}=a_{1,j_k-1}=a_{1,j_k+1}=1,\quad {\rm or}\quad \ a_{1,j_k}=a_{1,j_k-1}=a_{2,j_k+1}=1,\]
\[ {\rm or}\ \ a_{1,j_k}=a_{2,j_k-1}=a_{2,j_k+1}=1,\quad {\rm or}\quad a_{2,j_k}=a_{2,j_k-1}=a_{2,j_k+1}=1,\]
\[ {\rm or}\ \ a_{1,j_k}=a_{1,j_k+1}=a_{2,j_k-1}=1,\]
and the all others are equal to $0$. Thus, by Proposition \ref{dualsemi},
\begin{eqnarray}\varphi_{V_k} (x^G_{\rm \bf{i}}(1;\textbf{Y})) 
&=&Y_{1,j_k}Y_{1,j_k-1}Y_{1,j_k+1}+Y_{1,j_k}Y_{1,j_k-1}Y_{2,j_k+1}+
Y_{1,j_k}Y_{2,j_k-1}Y_{2,j_k+1}\nonumber \\ 
& &+Y_{2,j_k}Y_{2,j_k-1}Y_{2,j_k+1}+Y_{1,j_k}Y_{2,j_k-1}Y_{1,j_k+1}. \label{iniex3-1}
\end{eqnarray}

Similarly, it follows from $(\ref{mutex1})$ that for $1\leq k< \lfloor \frac{r+1}{2} \rfloor$, 
\begin{equation}\label{iniex3-2}
\varphi_{(\mu_k V)_k}(x^G_{\rm \bf{i}}(1;\textbf{Y}))
=\sum_{(*)}Y_{a,j_k-1}Y_{b,j_k+1}Y_{c,j_k-2}Y_{d,j_k}Y_{e,j_k+2},
\end{equation}
where $(*)$ is condition for $a,b,c,d$ and $e$ : $1\leq a\leq c\leq 2$, $1\leq b\leq d\leq 2$, $1\leq a\leq d\leq 2$ and $1\leq b\leq e\leq 2$. And
\begin{equation}\label{iniex3-2ano}
\varphi_{(\mu_{\lfloor \frac{r+1}{2} \rfloor} V)_{\lfloor \frac{r+1}{2} \rfloor}}(x^G_{\rm \bf{i}}(1;\textbf{Y}))
=\sum_{1\leq a\leq b\leq 2}Y_{a,2}Y_{b,3}.
\end{equation}
For $\lfloor \frac{r+1}{2} \rfloor+1\leq k\leq r$, it follows from $(\ref{mutex2})$ that
\begin{eqnarray}
& &\varphi_{(\mu_k V)_k}(x^G_{\rm \bf{i}}(1;\textbf{Y}))\nonumber \\
&=&
Y_{1,j_k-1}Y_{1,j_k+1}Y_{2,j_k-2}Y_{2,j_k}Y_{2,j_k+2}
Y_{2,j_k-3}Y_{2,j_k-1}Y_{2,j_k+1}Y_{2,j_k+3}. \label{iniex3-3}
\end{eqnarray}
\end{ex}

For two basic $\mathcal{C}_v$-cluster-tilting modules $R$, $R'$, we denote $R\sim R'$ if $R$ is obtained from $R'$ by a sequence of mutations (\ref{mudiretk}).

For $v\in W$, let $L(\mathcal{C}_v):=L(\mathcal{C}_v,V_{\textbf{i}})$ be the subalgebra of $\mathbb{C}[N]$ generated by $\{\varphi_{R_1},\varphi_{R_2},\cdots,\varphi_{R_n}| R_1\oplus R_2\oplus\cdots\oplus R_n\in {\rm Ob}(\mathcal{C}_v)\ \sim V_{\textbf{i}}\}$. Let $\tilde{L}(\mathcal{C}_v)$ be the cluster algebra obtained from $L(\mathcal{C}_v)$ by formally inverting the elements $\varphi_P$ for all $\mathcal{C}_v$-projective-injective module $P$. That is, $\tilde{L}(\mathcal{C}_v)$ is the localization of the ring $L(\mathcal{C}_v)$ with respect to $\varphi_P$.

\begin{thm}\label{GLSthm2}$\cite{GLS}$
For $v\in W$ and its reduced word ${\rm \bf{i}}=(j_n,\cdots,j_1)$, the coordinate ring $\mathbb{C}[L^{e,v}]$ have a cluster algebra structure, and the pair $((\varphi_{V_{{\rm \bf{i}},n}},\cdots,\varphi_{V_{{\rm \bf{i}},1}}),B(\Gamma_{V_{{\rm \bf{i}}}}))$ is its initial cluster. Moreover, we have
\[ \mathbb{C}[L^{e,v}]= \tilde{L}(\mathcal{C}_v).\] 
Furthermore, using the notation as in $(\ref{inc})$, we have $\varphi_{V_{{\rm \bf{i}},k}}=\Delta_{\Lambda_{j_k},v_{>n-k}\Lambda_{j_k}}|_{L^{e,v}}$. 
\end{thm}

\section{Monomial realizations and Demazure crystals}

In Sect.\ref{gmc}, we shall describe cluster variables in a cluster algebra of finite type in terms of the {\it monomial realizations} of Demazure crystals. Let us recall the notion of crystal base and its monomial realization in this section. Let $\ge$ be a complex simple Lie algebra and $I=\{1,2,\cdots,r\}$ the index set.

\subsection{Monomial realizations of crystals}\label{monoreacry}

In this subsection, we review the monomial realizations of crystals~\cite{K, K2, Nj}. First, let us recall the crystals.

\begin{defn}\label{defcry}$\cite{K0}$
A~{\it crystal} associated with a Cartan matrix~$A$ is a~set~$B$ together with the maps $\text{wt}: B \rightarrow P$,
$\tilde{e_{i}}$, $\tilde{f_{i}}: B \cup \{0\} \rightarrow B \cup \{0\}$ and $\varepsilon_i$,
$\varphi_i: B \rightarrow {\mathbb Z} \cup \{-\infty \}$, $i \in I$, satisfying some properties.
\end{defn}
We call $\{\tilde{e}_i,\ \tilde{f}_i\}$ the {\it Kashiwara operators}. Let $U_q(\mathfrak g)$ be the quantum enveloping algebra \cite{K0} associated with the Cartan matrix~$A$, that is,  $U_q(\mathfrak g)$ has generators $\{e_i,\ f_i,\ h_i |\ i\in I\}$ over $\mathbb{C}(q)$ satisfying some relations, where $q$ is an indeterminate. Let $V(\lm)$ ($\lm\in P^+=\oplus_{i\in I}\mathbb{Z}_{\geq0}\Lambda_i$) be the finite dimensional irreducible representation of $U_q(\mathfrak g)$ which has the highest weight vector $v_{\lm}$, and $B(\lm)$ be the crystal base of $V(\lm)$. The crystal base $B(\lm)$ has a crystal structure.

Let us introduce monomial realizations which realize each element of $B(\lambda)$ as a~certain Laurent monomial. First, fix a cyclic sequence of the indices $\cdots(i_1,i_2,\cdots,i_r)(i_1,i_2,\cdots,i_r)\cdots$ such that $\{i_1,i_2,\cdots,i_r\}=I$. And we can associate this sequence with a~set of integers $p=(p_{j,i})_{j,i \in I,\; j \neq i}$ such that
\begin{gather*}
p_{i_a,i_b}=
\begin{cases}
1 & \text{if} \quad  a<b,
\\
0 & \text{if} \quad  a>b.
\end{cases}
\end{gather*}

Second, for doubly-indexed variables $\{Y_{s,i} \,|\, i \in I$, $s\in \mathbb{Z}\}$, we def\/ine the set of monomials
\begin{gather*}
{\mathcal Y}:=\left\{Y=\prod\limits_{s \in \mathbb{Z},\ i \in I}
Y_{s,i}^{\zeta_{s,i}}\, \Bigg| \,\zeta_{s,i} \in \mathbb{Z},\
\zeta_{s,i} =0~\text{except for f\/initely many}~(s,i) \right\}.
\end{gather*}
Finally, we def\/ine maps $\text{wt}: {\mathcal Y} \rightarrow P$, $\varepsilon_i$, $\varphi_i: {\mathcal Y} \rightarrow
\mathbb{Z}$, $i \in I$. 
For $Y=\prod\limits_{s \in \mathbb{Z},\; i \in I} Y_{s,i}^{\zeta_{s,i}}\in {\mathcal Y}$,
\begin{gather}\label{wtph}
\text{wt}(Y):= \sum\limits_{i,s}\zeta_{s,i}\Lambda_i,\!
\quad
\varphi_i(Y):=\max\left\{\! \sum\limits_{k\leq s}\zeta_{k,i}  \,|\, s\in \mathbb{Z} \!\right\},\!
\quad
\varepsilon_i(Y):=\varphi_i(Y)-\text{wt}(Y)(h_i).
\end{gather}

We set
\begin{gather}
\label{asidef}
A_{s,i}:=Y_{s,i}Y_{s+1,i}\prod\limits_{j\neq i}Y_{s+p_{j,i},j}^{a_{j,i}}
\end{gather}
and def\/ine the Kashiwara operators as follows
\begin{gather*}
\tilde{f}_iY=
\begin{cases}
A_{n_{f_i},i}^{-1}Y & \text{if} \quad  \varphi_i(Y)>0,
\\
0 & \text{if} \quad  \varphi_i(Y)=0,
\end{cases}
\qquad
\tilde{e}_iY=
\begin{cases}
A_{n_{e_i},i}Y & \text{if} \quad  \varepsilon_i(Y)>0,
\\
0 & \text{if} \quad  \varepsilon_i(Y)=0,
\end{cases}
\end{gather*}
where
\begin{gather*}
n_{f_i}:=\min \left\{n \,\Bigg|\, \varphi_i(Y)= \sum\limits_{k\leq n}\zeta_{k,i}\right\},
\qquad
n_{e_i}:=\max \left\{n \,\Bigg|\, \varphi_i(Y)= \sum\limits_{k\leq n}\zeta_{k,i}\right\}.
\end{gather*}

Then the following theorem holds:

\begin{thm}$\cite{K, Nj}$\label{monorealmain}\quad
\begin{enumerate}\itemsep=0pt
\item[(i)] For the set $p=(p_{j,i})$ as above, $({\mathcal Y}, {\rm wt}, \varphi_i, \varepsilon_i,\tilde{f}_i,
\tilde{e}_i)_{i\in I}$ is a~crystal.
When we emphasize~$p$, we write ${\mathcal Y}$ as ${\mathcal Y}(p)$.
\item[(ii)] If a~monomial $Y \in {\mathcal Y}(p)$ satisfies $\varepsilon_i(Y)=0$ for all $i \in I$,
then the connected component containing~$Y$ is isomorphic to $B({\rm wt}(Y))$.
\end{enumerate}
\end{thm}

\subsection{Demazure crystals}\label{Demcrysub}

For $w\in W$ and $\lambda\in P^+$, a {\it Demazure crystal}
$B(\lambda)_w\subset B(\lambda)$ is inductively defined as follows. 

\begin{defn}
Let $u_{\lambda}$ be the highest weight vector of $B(\lambda)$. For
 the identity element $e$ of $W$, we set
 $B(\lambda)_e:=\{u_{\lambda}\}$. 
For $w\in W$, if $s_iw<w$, 
\[ B(\lambda)_w:=\{\tilde{f}_i^kb\ |\ k\geq0,\ b\in
 B(\lambda)_{s_iw},\ \tilde{e}_ib=0\}\setminus \{0\}. 
\]
\end{defn}

\begin{thm}\label{kashidem}${\cite{K3}}$
For $w\in W$, let $w=s_{i_1}\cdots s_{i_n}$ be an arbitrary reduced
 expression. Let $u_{\lambda}$ be the highest
weight vector of $B(\lambda)$. Then 
\[B(\lambda)_w=\{\tilde{f}_{i_1}^{a(1)}\cdots
\tilde{f}_{i_n}^{a(n)}u_{\lambda}|a(1),\cdots,a(n)\in\mathbb{Z}_{\geq0}
\}
\setminus\{0\}.
\]
\end{thm}

\begin{ex}\label{monorealex1}
Let us consider the case of type ${\rm A}_{r}$ and cyclic sequence is 
\[
\begin{cases}
(2,4,\cdots,r,1,3,5,\cdots,r-1) & {\rm if}\ r\ {\rm is\ even}, \\
(2,4,\cdots,r-1,1,3,5,\cdots,r) & {\rm if}\ r\ {\rm is\ odd}.
\end{cases}
\]
In the notation of $(\ref{asidef})$, we can write
\begin{equation}\label{asidefsp}
 A_{1,i}=
\begin{cases}
\frac{Y_{1,i}Y_{2,i}}{Y_{1,i-1}Y_{1,i+1}} & {\rm if}\ i\ {\rm is\ even}, \\
\frac{Y_{1,i}Y_{2,i}}{Y_{2,i-1}Y_{2,i+1}} & {\rm if}\ i\ {\rm is\ odd}. 
\end{cases}
\end{equation}
In general, if each factor of a monomial $Y\in\mathcal{Y}$ has non-negative degree, then $\varepsilon_i(Y)=0$ for all $i\in I$.
Therefore, we have $\varepsilon_i(Y_{1,j})=0$. Thus, we can consider the monomial realization of crystal base $B(\Lambda_j)$ such that the highest weight vector in $B(\Lambda_j)$ is realized by $Y_{1,j}$. The following is a part of it:
\[
\begin{xy}
(-40,85)*{Y_{1,j}}="h",
(0,85)*{Y_{1,j}A^{-1}_{1,j}}="1j",
(0,70)*{Y_{1,j}A^{-1}_{1,j}A^{-1}_{1,j-1}}="1j-1",
(40,85)*{Y_{1,j}A^{-1}_{1,j}A^{-1}_{1,j+1}}="1j+1",
(40,70)*{Y_{1,j}A^{-1}_{1,j}A^{-1}_{1,j+1}A^{-1}_{1,j-1}}="1j+1-1",
\ar@{->} "h";"1j"^{\tilde{f}_j}
\ar@{->} "1j";"1j+1"^{\tilde{f}_{j+1}}
\ar@{->} "1j";"1j-1"_{\tilde{f}_{j-1}}
\ar@{->} "1j+1";"1j+1-1"^{\tilde{f}_{j-1}}
\ar@{->} "1j-1";"1j+1-1"_{\tilde{f}_{j+1}}
\end{xy}
\]

\end{ex}

\section{Cluster variables and crystals}\label{gmc}

For $a,b\in\mathbb{Z}_{\geq0}$ with $a\leq b$, we set $[a,b]:=\{a,a+1,\cdots,b\}$. Let $G$ be a complex algebraic group of type A. In this section, we describe the cluster variables on some double Bruhat cell by the total sum of monomial realizations of Demazure crystals. In the rest of the paper, we only treat the Coxeter element $c\in W$ such that a reduced word ${\rm \bf{i}}$ of $c^2$ can be written as (\ref{redwords2}). Let $j_k$ be the $k$-th index of ${\rm \bf{i}}$ from the right, and consider the monomial realization associated with the sequence $(j_r,\cdots,j_2,j_1)$ (Sect.\ref{monoreacry}). The setting below is the same as in Example \ref{monorealex1}.

Let $\mathbb{V}:=((\varphi_{\mathbb{V}})_{2r},\cdots,(\varphi_{\mathbb{V}})_{r+1},(\varphi_{\mathbb{V}})_{r},\cdots,(\varphi_{\mathbb{V}})_{1},(\varphi_{\mathbb{V}})_{-r},\cdots,(\varphi_{\mathbb{V}})_{-1})$, where $(\varphi_{\mathbb{V}})_k\in \mathbb{C}[G^{e,c^2}]$ are defined as follows:
\[
(\varphi_{\mathbb{V}})_k=
\begin{cases}
\Delta_{\Lambda_{j_k},c^{2}_{>2r-k}\Lambda_{j_k}} & {\rm if}\ 1\leq k\leq 2r, \\
\Delta_{\Lambda_{|k|},\Lambda_{|k|}} & {\rm if}\ -r\leq k\leq -1.
\end{cases}
\]
By Theorem \ref{clmainthm}, Example \ref{clmainthmex} and Theorem \ref{finthm2}, we can regard $\mathbb{C}[G^{e,c^2}]$ as a cluster algebra of finite type and $\mathbb{V}$ as its initial cluster. Moreover, $(\varphi_{\mathbb{V}})_{2r},\cdots,(\varphi_{\mathbb{V}})_{r+1}$ and $(\varphi_{\mathbb{V}})_{-r},\cdots,(\varphi_{\mathbb{V}})_{-1}$ are frozen.
From Theorem \ref{GLSthm2}, for $k\in[1,2r]$, 
\begin{equation}\label{genminorgl}
 (\varphi_{\mathbb{V}})_{k}|_{L^{e,c^2}}=\varphi_{V_k}.
\end{equation}

Using these notation, we can rewrite (\ref{gammaex-1}) as
\begin{equation}\label{iniquiver}
\begin{xy}
(115,90)*{\cdots}="emp1",
(100,100)*{(\varphi_{\mathbb{V}})_{r+\lfloor \frac{r}{2} \rfloor + k}}="7",
(100,90) *{(\varphi_{\mathbb{V}})_{\lfloor \frac{r}{2} \rfloor + k}}="3",
(100,80)*{(\varphi_{\mathbb{V}})_{-j_k-1}}="-4",
(70,100)*{(\varphi_{\mathbb{V}})_{r+k}}="r+k",
(70,90) *{(\varphi_{\mathbb{V}})_k}="k",
(70,80)*{(\varphi_{\mathbb{V}})_{-j_k}}="-j_k",
(40,100)*{(\varphi_{\mathbb{V}})_{r+\lfloor \frac{r}{2} \rfloor +k+1}}="8",
(40,90) *{(\varphi_{\mathbb{V}})_{\lfloor \frac{r}{2} \rfloor +k+1}}="4",
(40,80)*{(\varphi_{\mathbb{V}})_{-j_k+1}}="-2",
(10,100)*{(\varphi_{\mathbb{V}})_{r+k+1}}="6",
(10,90) *{(\varphi_{\mathbb{V}})_{k+1}}="2",
(10,80)*{(\varphi_{\mathbb{V}})_{-j_k+2}}="-1",
(0,90)*{\cdots}="emp",
\ar@{->} "7";"3"
\ar@{->} "r+k";"k"
\ar@{->} "8";"4"
\ar@{->} "6";"2"
\ar@{->} "3";"-4"
\ar@{->} "k";"-j_k"
\ar@{->} "4";"-2"
\ar@{->} "2";"-1"
\ar@{->} "k";"3"
\ar@{->} "k";"4"
\ar@{->} "2";"4"
\ar@{->} "3";"r+k"
\ar@{->} "4";"r+k"
\ar@{->} "4";"6"
\ar@{->} "-4";"k"
\ar@{->} "-2";"k"
\ar@{->} "-2";"2"
\end{xy}
\end{equation}
Comparing with (\ref{gammavex-1}), we see that the matrix $B(\Gamma_{V_{{\rm \bf{i}}}})$ is a submatrix of $-\tilde{B}({\rm \bf{i})}$, which is obtained by deleting rows labelled by 
$(\varphi_{\mathbb{V}})_{-r},\cdots,(\varphi_{\mathbb{V}})_{-1}$ (Note that - sign of $-\tilde{B}({\rm \bf{i})}$ is needed to match the setting of \cite{A-F-Z} and \cite{GLS}). Note also that there are some differences between the quiver $\Gamma_{V_{\textbf{i}}}$ in (\ref{gammavex-1}) and the quiver obtained from $\Gamma_{\textbf{i}}$ by deleting the bottom row, that is, the arrows between the frozen cluster variables such as $r+\lfloor \frac{r}{2}\rfloor+k+1$, $r+k$, $r+\lfloor \frac{r}{2}\rfloor+k$ in (\ref{gammavex-1}). 

 In the rest of the paper, when we write a cluster in $\mathbb{C}[G^{e,c^2}]$, we drop frozen variables. For example, $\mathbb{V}=((\varphi_{\mathbb{V}})_{r},\cdots,(\varphi_{\mathbb{V}})_{1})$.
We will order the cluster variables $(\varphi_{\mathbb{V}})_{1},\cdots,(\varphi_{\mathbb{V}})_{r}$ from the right in $\mathbb{V}$ as above, and let $\mu_k$ denote the mutation of the $k$-th cluster variable from the right. For a cluster $\mathbb{T}$ in $\mathbb{C}[G^{e,c^2}]$, let $(\varphi_{\mathbb{T}})_{k}$ denote the $k$-th (non-frozen) cluster variable from the right:
\[ \mathbb{T}:=((\varphi_{\mathbb{T}})_{r},\cdots,(\varphi_{\mathbb{T}})_{1}).\] 
Each cluster variable is a regular function on $G^{e,c^2}$. On the other hand, by Proposition \ref{gprime}
, it can be seen as a function on
$H\times
(\mathbb{C}^{\times})^{2r}$ 
. Then, let us consider the following change of variables:
\begin{defn}\label{gendef}
Along (\ref{redwords2}), we set the variables $\textbf{Y}\in (\mathbb{C}^{\times})^{2r}$ as
\begin{equation}\label{yset}
\textbf{Y}:=
\begin{cases}
(Y_{1,2},Y_{1,4},\cdots,Y_{1,r},Y_{1,1},Y_{1,3},\cdots, Y_{1,r-1},Y_{2,2},\cdots,Y_{2,r}, Y_{2,1},\cdots, Y_{2,r-1}) & r\ {\rm is\ even,}\\
(Y_{1,2},Y_{1,4},\cdots,Y_{1,r-1},Y_{1,1},Y_{1,3},\cdots, Y_{1,r},Y_{2,2},\cdots,Y_{2,r-1}, Y_{2,1},\cdots, Y_{2,r}) & r\ {\rm is\ odd.}
\end{cases}
\end{equation}
Then for $a\in H$ and cluster $\mathbb{T}$ in $\mathbb{C}[G^{e,c^2}]$, we define
\[ 
(\varphi^G_{\mathbb{T}})_k(a;\textbf{Y}):=(\varphi_{\mathbb{T}})_k\circ
 \ovl{x}^G_{{\rm \bf{i}}}(a;\textbf{Y}),\quad (1\leq k\leq r),
\] 
where $\ovl{x}^G_{{\rm \bf{i}}}$ is as in \ref{fuctorisec}.
\end{defn}

\begin{ex}\label{initialex4}
In the above setting, 
\[ (\varphi^G_{\mathbb{V}})_k(a;\textbf{Y})=(\varphi_{\mathbb{V}})_k\circ x^G_{{\rm \bf{i}}}\circ \phi(a;\textbf{Y}), \]
where $\phi:H\times (\mathbb{C}^{\times})^{2r}\rightarrow H\times (\mathbb{C}^{\times})^{2r}$, 
\[
\phi(a;\textbf{Y})=(\Phi_H(a;\textbf{Y});\Phi_{1,j_{r}}(\textbf{Y}),\cdots,\Phi_{1,j_{1}}(\textbf{Y}),
\Phi_{2,j_r}(\textbf{Y}),\cdots,\Phi_{2,j_2}(\textbf{Y}),\Phi_{2,j_1}(\textbf{Y}))
\]
 is the map in the proof of Proposition \ref{gprime}. Since $(\varphi_{\mathbb{V}})_k$ is the generalized minor $\Delta_{\Lambda_{j_k},c^2_{>2r-k}\Lambda_{j_k}}$ (Theorem \ref{clmainthm}), we have
\begin{equation}\label{gtol}
 (\varphi_{\mathbb{V}})_k\circ x^G_{{\rm \bf{i}}}(a;\textbf{Y})=a^{\Lambda_{j_k}}
(\varphi_{\mathbb{V}})_k\circ x^G_{{\rm \bf{i}}}(1;\textbf{Y})
, \end{equation}
where $\textbf{Y}:=(Y_{1,j_{2r}},\cdots,Y_{1,j_{r+1}},Y_{2,j_r},\cdots,Y_{2,j_1})$ and $1$ is the identity element of $H$. By $(\ref{genminorgl})$, we obtain $(\varphi_{\mathbb{V}})_k\circ x^G_{{\rm \bf{i}}}(1;\textbf{Y})=\varphi_{V_k}\circ x^G_{{\rm \bf{i}}}(1;\textbf{Y})$. In Example \ref{initialex3}, we have calculated $\varphi_{V_k}\circ x^G_{{\rm \bf{i}}}(1;\textbf{Y})$.

If $1\leq k\leq \lfloor \frac{r+1}{2} \rfloor$, $j_k$ is odd. Then, since we know that $\varphi_{V_k}\circ x^G_{{\rm \bf{i}}}(1;\textbf{Y})=Y_{1,j_k}+Y_{2,j_k}$, it follows from $(\ref{mbasea})$, $(\ref{mbase0})$ and $(\ref{mbase01})$ that 
\begin{eqnarray*}
& &\hspace{-30pt}(\varphi^G_{\mathbb{V}})_k(a;\textbf{Y})=(\Phi_H (a;\textbf{Y}))^{\Lambda_{j_k}} (\Phi_{1,j_k} (a;\textbf{Y})+\Phi_{2,j_k} (a;\textbf{Y}))\\
&=&a^{\Lambda_{j_k}}(Y_{1,j_k}Y_{2,j_k})\left(\frac{Y_{2,j_k-1}Y_{2,j_k+1}}{Y_{1,j_k}Y^2_{2,j_k}}+\frac{1}{Y_{2,j_k}}\right)
=a^{\Lambda_{j_k}}Y_{1,j_k}(1+A^{-1}_{1,j_k}),
\end{eqnarray*}
where $A_{1,j}$ is given in $(\ref{asidefsp})$. By Theorem \ref{kashidem} and Example \ref{monorealex1}, the set of monomials $\{Y_{1,j_k},Y_{1,j_k}A^{-1}_{1,j_k}\}$ coincides with the monomial realization of Demazure crystal $B(\Lambda_{j_k})_{s_{j_k}}$, where the monomial corresponding to the highest weight vector is $Y_{1,j_k}$.

For $\lfloor \frac{r+1}{2} \rfloor+1\leq k\leq r$, $j_k$ is even. In this case, we calculated $\varphi_{V_k}\circ x^G_{{\rm \bf{i}}}(1;\textbf{Y})$ in $(\ref{iniex3-1})$. Thus, using $(\ref{mbasea})$, $(\ref{mbase0})$ and $(\ref{mbase01})$,
\begin{eqnarray*}
& &\hspace{-20pt}(\varphi^G_{\mathbb{V}})_k(a;\textbf{Y})=a^{\Lambda_{j_k}}(Y_{1,j_k}Y_{2,j_k})\\
&\times &
\left(\frac{Y_{2,j_k-2}Y_{2,j_k+2}}{Y_{1,j_k}Y_{2,j_k-1}Y_{2,j_k+1}}+\frac{Y_{1,j_k+1}Y_{2,j_k-2}}{Y_{1,j_k}Y_{2,j_k}Y_{2,j_k-1}}
+\frac{Y_{1,j_k-1}Y_{1,j_k+1}}{Y_{1,j_k}Y^2_{2,j_k}}+\frac{1}{Y_{2,j_k}}+\frac{Y_{1,j_k-1}Y_{2,j_k+2}}{Y_{1,j_k}Y_{2,j_k}Y_{2,j_k+1}}\right)\\
&=&a^{\Lambda_{j_k}}\left(Y_{1,j_k}+\frac{Y_{1,j_k-1}Y_{1,j_k+1}}{Y_{2,j_k}}
+\frac{Y_{1,j_k+1}Y_{2,j_k-2}}{Y_{2,j_k-1}}+\frac{Y_{1,j_k-1}Y_{2,j_k+2}}{Y_{2,j_k+1}}+
\frac{Y_{2,j_k-2}Y_{2,j_k}Y_{2,j_k+2}}{Y_{2,j_k-1}Y_{2,j_k+1}}
\right)\\
&=&a^{\Lambda_{j_k}}Y_{1,j_k}(1+A^{-1}_{1,j_k}+A^{-1}_{1,j_k}A^{-1}_{1,j_k-1}+A^{-1}_{1,j_k}A^{-1}_{1,j_k+1}
+A^{-1}_{1,j_k}A^{-1}_{1,j_k-1}A^{-1}_{1,j_k+1}).
\end{eqnarray*}
By Theorem \ref{kashidem} and Example \ref{monorealex1}, the set of monomials $\{Y_{1,j_k},\frac{Y_{1,j_k-1}Y_{1,j_k+1}}{Y_{2,j_k}},\frac{Y_{1,j_k+1}Y_{2,j_k-2}}{Y_{2,j_k-1}}$,
$\frac{Y_{1,j_k-1}Y_{2,j_k+2}}{Y_{2,j_k+1}},\frac{Y_{2,j_k-2}Y_{2,j_k}Y_{2,j_k+2}}{Y_{2,j_k-1}Y_{2,j_k+1}}\}$ coincides with the monomial realization of Demazure crystal $B(\Lambda_{j_k})_{s_{j_k+1}s_{j_k-1}s_{j_k}}$, where the monomial corresponding to the highest weight vector is $Y_{1,j_k}$.

Next, let us consider the mutation in direction $k$ of $\mathbb{V}$ by calculating $(\varphi^G_{(\mu_k\mathbb{V})})_k(a;\textbf{Y})$ $(1\leq k\leq r)$. If $1\leq k\leq \lfloor \frac{r+1}{2} \rfloor$, by $(\ref{gammaex-1})$,
\begin{eqnarray}
& &\hspace{-20pt} (\varphi_{(\mu_k\mathbb{V})})_k\circ x^G_{{\rm \bf{i}}}(a;\textbf{Y})\nonumber\\
&=&\left(
\frac{(\varphi_{\mathbb{V}})_{r+k}(\varphi_{\mathbb{V}})_{-j_k+1}(\varphi_{\mathbb{V}})_{-j_k-1}+
(\varphi_{\mathbb{V}})_{-j_k}(\varphi_{\mathbb{V}})_{ \lfloor \frac{r}{2} \rfloor+k} (\varphi_{\mathbb{V}})_{ \lfloor \frac{r}{2} \rfloor+k+1}
}{(\varphi_{\mathbb{V}})_k}\right) \circ x^G_{{\rm \bf{i}}}(a;\textbf{Y})\nonumber\\
&=& \frac{a^{\Lambda_{j_k-1}}a^{\Lambda_{j_k}}a^{\Lambda_{j_k+1}}}{a^{\Lambda_{j_k}}}\cdot
\left(
\frac{\varphi_{V_{r+k}}+
\varphi_{V_{ \lfloor \frac{r}{2} \rfloor+k}} \varphi_{V_{ \lfloor \frac{r}{2} \rfloor+k+1}}
}{\varphi_{V_k}}\right) \circ x^G_{{\rm \bf{i}}}(1;\textbf{Y})\nonumber\\
&=& a^{\Lambda_{j_k-1}}a^{\Lambda_{j_k+1}}\cdot
(\varphi_{(\mu_k V)_k})\circ x^G_{{\rm \bf{i}}}(1;\textbf{Y})\label{mukcoeff}, 
\end{eqnarray}
where we use $(\ref{genbasic})$. From $(\ref{iniex3-2})$, for $1\leq k< \lfloor \frac{r+1}{2} \rfloor$, we get
\begin{eqnarray*}
& &\hspace{-30pt} (\varphi^G_{(\mu_k\mathbb{V})})_k (a;\textbf{Y})
=(\varphi_{(\mu_k\mathbb{V})})_k\circ x^G_{{\rm \bf{i}}}\circ \phi(a;\textbf{Y})\\
&=&{\Phi_H}^{\Lambda_{j_k-1}}{\Phi_H}^{\Lambda_{j_k+1}}\cdot
\sum_{(*)}\Phi_{b_1,j_k-1}\Phi_{b_2,j_k+1}\Phi_{b_3,j_k-2}\Phi_{b_4,j_k}\Phi_{b_5,j_k+2},
\end{eqnarray*}
where $(*)$ is the condition for $b_1,b_2,b_3,b_4$ and $b_5$ : $1\leq b_1\leq b_3\leq 2$, $1\leq b_2\leq b_4\leq 2$, $1\leq b_1\leq b_4\leq 2$ and $1\leq b_2\leq b_5\leq 2$. We can easily verify that if $b_4=2$, then $1\leq b_1\leq b_3\leq 2$ and $1\leq b_2\leq b_5\leq 2$. If $b_4=1$, then $b_1=b_2=1$ and $1\leq b_3,\ b_5\leq 2$. By $(\ref{mbase0})$ and $(\ref{mbase01})$, we see that 
\begin{equation}\label{phi12}
\Phi_{1,j}=A^{-1}_{1,j}\Phi_{2,j}
\end{equation}
for $j\in I$. Thus,
\begin{multline*}
\sum_{(*),\ b_4=2}\Phi_{b_1,j_k-1}\Phi_{b_2,j_k+1}\Phi_{b_3,j_k-2}\Phi_{b_4,j_k}\Phi_{b_5,j_k+2}= \\
\Phi_{2,j_k-1}\Phi_{2,j_k+1}\Phi_{2,j_k-2}\Phi_{2,j_k}\Phi_{2,j_k+2}
(1+A^{-1}_{1,j_k-1}+A^{-1}_{1,j_k-1}A^{-1}_{1,j_k-2})(1+A^{-1}_{1,j_k+1}+A^{-1}_{1,j_k+1}A^{-1}_{1,j_k+2})
\\
=\frac{Y_{2,j_k}}{Y_{2,j_k-1}Y_{2,j_k+1}}(1+A^{-1}_{1,j_k-1}+A^{-1}_{1,j_k-1}A^{-1}_{1,j_k-2})(1+A^{-1}_{1,j_k+1}+A^{-1}_{1,j_k+1}A^{-1}_{1,j_k+2})\quad,
\end{multline*}

\begin{eqnarray*}
& &\sum_{(*),\ b_4=1}\Phi_{b_1,j_k-1}\Phi_{b_2,j_k+1}\Phi_{b_3,j_k-2}\Phi_{b_4,j_k}\Phi_{b_5,j_k+2} \\
&=&\Phi_{1,j_k-1}\Phi_{1,j_k+1}\Phi_{2,j_k-2}\Phi_{1,j_k}\Phi_{2,j_k+2}
(1+A^{-1}_{1,j_k-2})(1+A^{-1}_{1,j_k+2})
\\
&=&\frac{Y_{1,j_k-2}Y_{1,j_k}Y_{1,j_k+2}}{Y_{1,j_k-1}Y_{2,j_k-1}Y_{1j_k+1}Y_{2,j_k+1}}(1+A^{-1}_{1,j_k-2})(1+A^{-1}_{1,j_k+2}).
\end{eqnarray*}

By the above argument, we get
\begin{multline*}
(\varphi^G_{(\mu_k\mathbb{V})})_k (a;\textbf{Y})=
{a}^{\Lambda_{j_k-1}}{a}^{\Lambda_{j_k+1}} (Y_{1,j_k-2}Y_{1,j_k}Y_{1,j_k+2}(1+A^{-1}_{1,j_k-2})(1+A^{-1}_{1,j_k+2})\\
+Y_{1,j_k-1}Y_{2,j_k}Y_{1,j_k+1}(1+A^{-1}_{1,j_k-1}+A^{-1}_{1,j_k-1}A^{-1}_{1,j_k-2})(1+A^{-1}_{1,j_k+1}+A^{-1}_{1,j_k+1}A^{-1}_{1,j_k+2})).
\end{multline*}
By the definition of Kashiwara operators in \ref{monoreacry}, we see that $Y_{1,j_k-2}Y_{1,j_k}Y_{1,j_k+2}(1+A^{-1}_{1,j_k-2})(1+A^{-1}_{1,j_k+2})$ are the total sum of the monomial realization of the Demazure crystal $B(\Lambda_{j_k-2}+\Lambda_{j_k}+\Lambda_{j_k+2})_{s_{j_k-2}s_{j_k+2}}$, and $Y_{1,j_k-1}Y_{2,j_k}Y_{1,j_k+1}(1+A^{-1}_{1,j_k-1}+A^{-1}_{1,j_k-1}A^{-1}_{1,j_k-2})(1+A^{-1}_{1,j_k+1}+A^{-1}_{1,j_k+1}A^{-1}_{1,j_k+2})$ is the total sum of the monomial realization of the Demazure crystal $B(\Lambda_{j_k-1}+\Lambda_{j_k}+\Lambda_{j_k+1})_{s_{j_k-2}s_{j_k-1}s_{j_k+2}s_{j_k+1}}$.
Note that $\Lambda_{j_k-2}+\Lambda_{j_k}+\Lambda_{j_k+2}=(\Lambda_{j_k-1}+\Lambda_{j_k}+\Lambda_{j_k+1})-\al_{j_k-1}-\al_{j_k}-\al_{j_k+1}$.

Arguing similarly, we obtain 
\[
(\varphi^G_{(\mu_{\lfloor \frac{r+1}{2} \rfloor}\mathbb{V})})_{\lfloor \frac{r+1}{2} \rfloor} (a;\textbf{Y})=a^{\Lambda_2}Y_{1,2}Y_{2,1}(1+A^{-1}_{1,2}+A^{-1}_{1,2}A^{-1}_{1,3}).
\]
The polynomial $Y_{1,2}Y_{2,1}(1+A^{-1}_{1,2}+A^{-1}_{1,2}A^{-1}_{1,3})$ is the total sum of the monomial realization of the Demazure crystal $B(\Lambda_{1}+\Lambda_2)_{s_3s_2}$.

Similarly, for $\lfloor \frac{r+1}{2} \rfloor < k \leq r$, it follows from 
$(\ref{gammaex-2})$ and $(\ref{iniex3-3})$ that
\begin{eqnarray*}
& &\hspace{-20pt}(\varphi^G_{(\mu_k\mathbb{V})})_k (a;\textbf{Y})=
{\Phi_H}^{\Lambda_{j_k-1}}{\Phi_H}^{\Lambda_{j_k+1}}
\times \Phi_{1,j_k-1}\Phi_{1,j_k+1}\Phi_{2,j_k-2}\Phi_{2,j_k}\Phi_{2,j_k+2}\\
& &\times\Phi_{2,j_k-3}\Phi_{2,j_k-1}\Phi_{2,j_k+1}\Phi_{2,j_k+3}
={a}^{\Lambda_{j_k-1}}{a}^{\Lambda_{j_k+1}} Y_{2,j_k}.
\end{eqnarray*}
The monomial $Y_{2,j_k}$ is the monomial realization of the Demazure crystal $B(\Lambda_{j_k})_e$. 
\end{ex}
Using Proposition \ref{remaindia}, we obtain the following Proposition \ref{remainexp} by the same way as the above example:

\begin{prop}\label{remainexp}
We have the cluster variable
\[(\varphi^G_{(\mu_r\mu_{\lfloor \frac{r+1}{2} \rfloor}\mathbb{V})})_r (a;\textbf{Y})=a^{\Lambda_3}Y_{2,1},\]
and the monomial realization of the Demazure crystal $B(\Lambda_1)_e$. For $\lfloor \frac{r+1}{2} \rfloor+2\leq k\leq r$, we obtain
\begin{eqnarray*}
& &\hspace{-30pt}(\varphi^G_{(\mu_{k-\lfloor \frac{r}{2} \rfloor-1}\mu_{k}\mathbb{V})})_{k-\lfloor \frac{r}{2} \rfloor-1} (a;\textbf{Y})\\
&=&a^{\Lambda_{j_k-1}+\Lambda_{j_k+2}}Y_{2,j_k+1}Y_{1,j_k+2}(1+A^{-1}_{1,j_k+2}+A^{-1}_{1,j_k+2}A^{-1}_{1,j_k+3}),
\end{eqnarray*}
and the set $\{Y_{2,j_k+1}Y_{1,j_k+2},Y_{2,j_k+1}Y_{1,j_k+2}A^{-1}_{1,j_k+2},Y_{2,j_k+1}Y_{1,j_k+2}A^{-1}_{1,j_k+2}A^{-1}_{1,j_k+3}\}$ coincides with the monomial realization of the Demazure crystal $B(\Lambda_{j_k+1}+\Lambda_{j_k+2})_{s_{j_k+3}s_{j_k+2}}$, where $Y_{2,j_k+1}Y_{1,j_k+2}$ is the monomial corresponding to the highest weight vector in $B(\Lambda_{j_k+1}+\Lambda_{j_k+2})$.
\end{prop}

In the following Proposition \ref{spprop1}, \ref{spprop2} and \ref{spprop3}, we shall give explicit expressions of all the other cluster variables in $\mathbb{C}[G^{e,c^2}]$. We use the notation as in (\ref{mbasea}), (\ref{mbase0}), (\ref{mbase01}) and (\ref{asidefsp}), and set $\phi(\textbf{Y}):=(\Phi_{1,j_{r}}(\textbf{Y}),\cdots,\Phi_{1,j_{1}}(\textbf{Y}),\Phi_{2,j_r}(\textbf{Y}),\cdots,\Phi_{2,j_1}(\textbf{Y}))$. We abbreviate $\Phi_H(a;\textbf{Y})$ to $\Phi_H$. For the integers $b$, $c$ $(b<c)$ and $x$, we set
\[ A[b,c;x]:=\left(\prod^{c-1}_{s=b} A_{1,x-2s-2}A_{1,x-2s-3}\right)^{-1} A_{1,x-2c-2}^{-1}. \]
For $p\in\mathbb{Z}_{> 0}$ and $\textbf{b}=(b_i)^p_{i=1}\in (\mathbb{Z}_{\geq 0})^p$,
$\textbf{c}=(c_i)^p_{i=1}\in (\mathbb{Z}_{\geq 0})^p$ such that $b_i<c_i$ $(1\leq i\leq p)$, we also set
\[
\al[{\rm \bf{b}},{\rm \bf{c}};x]:=\sum^{x-2b_1-2}_{t=x-2c_1-2}\al_t+\cdots+\sum^{x-2b_p-2}_{t=x-2c_p-2}\al_t,
\]
where when $s\leq0$, we understand $A_{1,s}=1$, $\al_s=0$.
For $l\in\mathbb{Z}_{\geq 0}$, we define
\[ R^p_{l}:=\{(\textbf{b},\textbf{c})\in (\mathbb{Z}_{\geq 0})^p\times (\mathbb{Z}_{\geq 0})^p|\ 
\textbf{b}=(b_i)^p_{i=1},\ \textbf{c}=(c_i)^p_{i=1},\ 0\leq b_1<c_1<\cdots<b_p<c_p\leq l\}. \]
For $(\textbf{b},\textbf{c})\in R^p_{l}$, we define $[\textbf{b},\textbf{c}]:=[b_1,c_1]\cup\cdots\cup[b_p,c_p]$.
\begin{prop}\label{spprop1} For $k\in[\lfloor \frac{r+1}{2} \rfloor+1,r-1]$ and $l\in[0,r-k-1]$,
let $\mu[l]$ be the following iteration of mutations
\[\mu[l]:=
(\mu_{k-\lfloor \frac{r}{2} \rfloor +l}\mu_{k+l+1}\mu_{k+l})\cdots(\mu_{k-\lfloor \frac{r}{2} \rfloor +1}\mu_{k+2}\mu_{k+1})(\mu_{k-\lfloor \frac{r}{2} \rfloor }\mu_{k+1}\mu_{k}),
\]
where $\lfloor \ \rfloor $ is the Gaussian symbol. 
\begin{enumerate}
\item[$(a)$]We have the cluster variable
\begin{eqnarray*}
& &\hspace{-20pt}(\varphi^G_{(\mu[l]\mathbb{V})})_{k-\lfloor \frac{r}{2} \rfloor +l}(a;\textbf{Y})\\
&=&\Phi^{(\sum^{l+2}_{s=0}\Lambda_{j_k-2s+1})+(\sum^{l-1}_{s=0}\Lambda_{j_k-2s-2})}_H \varphi_{(\mu[l]V)_{k-\lfloor \frac{r}{2} \rfloor +l}}\circ x^G_{{\rm \bf{i}}}(1;\phi(\textbf{Y})) \\
&=&a^{(\sum^{l+2}_{s=0}\Lambda_{j_k-2s+1})+(\sum^{l-1}_{s=0}\Lambda_{j_k-2s-2})}H_1\left(
\prod_{q\in[0,l-1]}(1+A_{1,j_k-2q-2}^{-1})\right.\\
&+&\left.\sum_{p>0}\sum_{({\rm \bf{b}},{\rm \bf{c}})\in R^p_{l-1}} 
\prod^p_{i=1} A[b_i,c_i;j_k]
\prod_{q\in[0,l-1]\setminus[{\rm \bf{b}},{\rm \bf{c}}]}(1+A_{1,j_k-2q-2}^{-1})\right).
\end{eqnarray*}
\item[$(b)$]We also obtain the cluster variable
\begin{eqnarray*}
& &\hspace{-20pt}(\varphi^G_{(\mu_{k+l+1}\mu[l]\mathbb{V})})_{k+l+1}(a;\textbf{Y})\\
&=& \Phi^{\sum^{l}_{s=0}(\Lambda_{j_k-2s+1}+\Lambda_{j_k-2s-2})}_H \varphi_{(\mu_{k+l+1}\mu[l]V)_{k+l+1}}\circ x^G_{{\rm \bf{i}}}(1;\phi(\textbf{Y}))\\
&=& a^{\sum^{l}_{s=0}(\Lambda_{j_k-2s+1}+\Lambda_{j_k-2s-2})} H_2\\
&\times&\left( 
(1+A^{-1}_{1,j_k-2l-2}+A^{-1}_{1,j_k-2l-2}A^{-1}_{1,j_k-2l-3})
\prod_{q\in[0,l-1]}(1+A_{1,j_k-2q-2}^{-1})
\right.\\
&+&\sum_{p>0}
\sum_{({\rm \bf{b}},{\rm \bf{c}})\in R^p_l}(1-\delta_{c_p,l}+A^{-1+\delta_{c_p,l}}_{1,j_k-2l-2}(1+A^{-1}_{1,j_k-2l-3}))\\
&\times& \left. \prod^p_{i=1} A[b_i,c_i;j_k]
\prod_{q\in[0,l-1]\setminus[{\rm \bf{b}},{\rm \bf{c}}]}(1+A_{1,j_k-2q-2}^{-1})\right)
,\end{eqnarray*}
where
\[H_1:=\left(\prod^{l-1}_{t=0} Y_{1,j_k-2t-2}\right)\left(\prod^{l}_{t=0}Y_{2,j_k-2t-1}\right),
\ \ H_2:=\left(\prod^{l}_{t=0} Y_{1,j_k-2t-2}Y_{2,j_k-2t-1}\right).
\]
\end{enumerate}
\end{prop}

\begin{ex}\label{r10ex}
If $r=10$, $k=6$ and $l=2$, then $\mu[2]=\mu_3\mu_9\mu_8\mu_{2}\mu_8\mu_7\mu_{1}\mu_7\mu_6$, $j_6=10$ and $H_1=Y_{2,5}Y_{1,6}Y_{2,7}Y_{1,8}Y_{2,9}$ in the notation of Proposition \ref{spprop1}. Note that 
\[
R^p_1=
\begin{cases}
\{(0,1)\} & {\rm if}\ p=1,\\
\phi & {\rm otherwise}.
\end{cases}
\]
It follows from Proposition \ref{spprop1} (a) that
\begin{eqnarray}
(\varphi^G_{(\mu[2]\mathbb{V})})_{3}(a;\textbf{Y})&=&
a^{\Lm_3+\Lm_5+\Lm_6+\Lm_7+\Lm_8+\Lm_9}Y_{2,5}Y_{1,6}Y_{2,7}Y_{1,8}Y_{2,9}(1+A^{-1}_{1,8})(1+A^{-1}_{1,6})\nonumber\\
& &+ a^{\Lm_3+\Lm_5+\Lm_6+\Lm_7+\Lm_8+\Lm_9}Y_{2,5}Y_{1,6}Y_{2,7}Y_{1,8}Y_{2,9}A[0,1;10]\nonumber\\
&=&a^{\Lm_3+\Lm_5+\Lm_6+\Lm_7+\Lm_8+\Lm_9}Y_{2,5}Y_{1,6}Y_{2,7}Y_{1,8}Y_{2,9}(1+A^{-1}_{1,8}+A^{-1}_{1,6}+A^{-1}_{1,6}A^{-1}_{1,8})\nonumber\\
& &+a^{\Lm_3+\Lm_5+\Lm_6+\Lm_7+\Lm_8+\Lm_9}Y_{1,5}Y_{2,5}Y_{1,7}Y_{1,9}Y_{2,9}. \label{r10ex-1}
\end{eqnarray}
In the same setting, let us calculate $(\varphi^G_{(\mu_9\mu[2]\mathbb{V})})_{9}(a;\textbf{Y})$. Note that
\[
R^p_2=
\begin{cases}
\{(0,1), (0,2), (1,2)\} & {\rm if}\ p=1,\\
\phi & {\rm otherwise},
\end{cases}
\]
and $H_2=Y_{1,4}Y_{2,5}Y_{1,6}Y_{2,7}Y_{1,8}Y_{2,9}$. Thus, by Proposition \ref{spprop1} (b), 
\begin{eqnarray}
& &\hspace{-20pt}a^{-(\Lm_4+\Lm_6+\Lm_7+\Lm_8+\Lm_9)}(\varphi^G_{(\mu_9\mu[2]\mathbb{V})})_{9}(a;\textbf{Y})\nonumber\\
&=&
Y_{1,4}Y_{2,5}Y_{1,6}Y_{2,7}Y_{1,8}Y_{2,9}(1+A^{-1}_{1,8})(1+A^{-1}_{1,6})(1+A^{-1}_{1,4}+A^{-1}_{1,4}A^{-1}_{1,3})\nonumber\\
& &+ Y_{1,4}Y_{2,5}Y_{1,6}Y_{2,7}Y_{1,8}Y_{2,9}A[0,1;10](1+A^{-1}_{1,4}+A^{-1}_{1,4}A^{-1}_{1,3})\nonumber\\
& &+ Y_{1,4}Y_{2,5}Y_{1,6}Y_{2,7}Y_{1,8}Y_{2,9}A[0,2;10](1+A^{-1}_{1,3})\nonumber\\
& &+ Y_{1,4}Y_{2,5}Y_{1,6}Y_{2,7}Y_{1,8}Y_{2,9}A[1,2;10](1+A^{-1}_{1,3})\nonumber\\
&=&
Y_{1,4}Y_{2,5}Y_{1,6}Y_{2,7}Y_{1,8}Y_{2,9}(1+A^{-1}_{1,6}+A^{-1}_{1,8}+A^{-1}_{1,6}A^{-1}_{1,8})(1+A^{-1}_{1,4}+A^{-1}_{1,4}A^{-1}_{1,3})\nonumber\\
& &+ Y_{1,4}Y_{1,5}Y_{2,5}Y_{1,7}Y_{1,9}Y_{2,9}(1+A^{-1}_{1,4}+A^{-1}_{1,4}A^{-1}_{1,3})\nonumber\\
& &+ Y_{1,3}Y_{1,5}Y_{2,6}Y_{1,7}Y_{1,9}Y_{2,9}(1+A^{-1}_{1,3})\label{r10ex-2}\\
& &+ Y_{1,3}Y_{1,5}Y_{1,7}Y_{2,7}Y_{1,8}Y_{2,9}(1+A^{-1}_{1,3}).\nonumber
\end{eqnarray}

\end{ex}

\begin{prop}\label{spprop2} For $k\in[1,\lfloor \frac{r+1}{2} \rfloor-2]$ and $l\in[0,\lfloor \frac{r+1}{2} \rfloor -k-2]$,
let $\mu'[l]$ be the following iteration of mutations 
\[\mu'[l]:=
(\mu_{k+l+1}\mu_{\lfloor \frac{r}{2} \rfloor +k+l+2}\mu_{\lfloor \frac{r}{2} \rfloor +k+l+1})\cdots(\mu_{k+2}\mu_{\lfloor \frac{r}{2} \rfloor +k+3}\mu_{\lfloor \frac{r}{2} \rfloor +k+2})(\mu_{k+1}\mu_{\lfloor \frac{r}{2} \rfloor +k+2}\mu_{\lfloor \frac{r}{2} \rfloor +k+1})\mu_{k}.
\]
\begin{enumerate}
\item[$(a)$]If $j_k<r$, we have the cluster variable
\begin{eqnarray*}
& &\hspace{-20pt}(\varphi^G_{(\mu'[l]\mathbb{V})})_{k+l+1}(a;\textbf{Y})
=a^{\sum^{l+1}_{s=0}\Lambda_{j_k-2s-2}+\Lambda_{j_k-2s+1}}H_3\\
&\times&\left( (1+A^{-1}_{1,j_k+1}+A^{-1}_{1,j_k+1}A^{-1}_{1,j_k+2})
\prod_{q\in[1,l+1]}(1+A_{1,j_k-2q+1}^{-1})
\right.\\
&+&
\sum_{p>0}\sum_{({\rm \bf{b}},{\rm \bf{c}})\in R^p_{l+1}}
(1-\delta_{b_1,0}+A^{-1+\delta_{b_1,0}}_{1,j_k+1}(1+A^{-1}_{1,j_k+2})) \\
&\times& \left. \prod^p_{i=1} A[b_i,c_i;j_k+3]
\prod_{q\in[1,l+1]\setminus([{\rm \bf{b}},{\rm \bf{c}}])}(1+A_{1,j_k-2q+1}^{-1})\right),
\end{eqnarray*}
and if $j_k=r$, we have
\begin{eqnarray*}
& &\hspace{-20pt}(\varphi^G_{(\mu'[l]\mathbb{V})})_{k+l+1}(a;\textbf{Y})=a^{\sum^{l+1}_{s=0}\Lambda_{r-2s-2}+\Lambda_{r-2s+1}}H_3 \left(\prod_{q\in[1,l+1]}(1+A_{1,r-2q+1}^{-1}) \right.\\
&+&\left.\sum_{p>0}\sum_{({\rm \bf{b}},{\rm \bf{c}})\in R^p_{l+1},b_1>0} \prod^p_{i=1} A[b_i,c_i;r+3]
\prod_{q\in[1,l+1]\setminus([{\rm \bf{b}},{\rm \bf{c}}])}(1+A_{1,r-2q+1}^{-1})\right).
\end{eqnarray*}
\item[$(b)$]If $j_k<r$, we also obtain the cluster variable
\begin{eqnarray*}
& &\hspace{-20pt}(\varphi^G_{(\mu_{\lfloor \frac{r}{2} \rfloor +k+l+2}\mu'[l]\mathbb{V})})_{\lfloor \frac{r}{2} \rfloor +k+l+2}(a;\textbf{Y})=a^{\sum^{l-1}_{s=0}\Lambda_{j_k-2s-2}+\sum^{l+2}_{s=0}\Lambda_{j_k-2s+1}}H_4 \\
&\times&\left( \prod_{q\in[1,l+1]}(1+A_{1,j_k-2q+1}^{-1})
(1+A^{-1}_{1,j_k-2l-3}+A^{-1}_{1,j_k-2l-3}A^{-1}_{1,j_k-2l-4})
\right. \\
&\times&
(1+A^{-1}_{1,j_k+1}+A^{-1}_{1,j_k+1}A^{-1}_{1,j_k+2})+
\sum_{p>0}\sum_{({\rm \bf{b}},{\rm \bf{c}})\in R^p_{l+2}}(1-\delta_{b_1,0}+A^{-1+\delta_{b_1,0}}_{1,j_k+1}(1+A^{-1}_{1,j_k+2}))\\
&\times&
(1-\delta_{c_p,l+2}+A^{-1+\delta_{c_p,l+2}}_{1,j_k-2l-3}(1+A^{-1}_{1,j_k-2l-4}))\\
&\times&\left. \prod^p_{i=1} A[b_i,c_i;j_k+3]
\prod_{q\in[1,l+1]\setminus([{\rm \bf{b}},{\rm \bf{c}}])}(1+A_{1,j_k-2q+1}^{-1})\right),
\end{eqnarray*}
and if $j_k=r$, we have
\begin{eqnarray*}
& &\hspace{-20pt}(\varphi^G_{(\mu_{\lfloor \frac{r}{2} \rfloor +k+l+2}\mu'[l]\mathbb{V})})_{\lfloor \frac{r}{2} \rfloor +k+l+2}(a;\textbf{Y})=a^{\sum^{l-1}_{s=0}\Lambda_{r-2s-2}+\sum^{l+2}_{s=0}\Lambda_{r-2s+1}} H_4\\
&\times&\left(
(1+A^{-1}_{1,r-2l-3}+A^{-1}_{1,r-2l-3}A^{-1}_{1,r-2l-4})
\prod_{q\in[1,l+1]}(1+A_{1,r-2q+1}^{-1})
\right.\\
&+&\sum_{p>0}\sum_{({\rm \bf{b}},{\rm \bf{c}})\in R^p_{l+2},\ b_1>0}
\prod^p_{i=1} A[b_i,c_i;r+3]\\
&\times& \left.(1-\delta_{c_p,l+2}+A^{-1+\delta_{c_p,l+2}}_{1,r-2l-3}(1+A^{-1}_{1,r-2l-4}))
\prod_{q\in[1,l+1]\setminus([{\rm \bf{b}},{\rm \bf{c}}])}(1+A_{1,r-2q+1}^{-1})\right),
\end{eqnarray*}
where $H_3:=\left(\prod^{l+1}_{t=0}Y_{1,j_k-2t+1}Y_{2,j_k-2t}\right)$, $H_4:=\left(\prod^{l+2}_{t=0}Y_{1,j_k-2t+1}\right)\left(\prod^{l+1}_{t=0} Y_{2,j_k-2t}\right)
=H_3\times Y_{1,j_k-2l-3}$.. If $j_k=r$, then we understand $Y_{1,j_k+1}=1$ and $\Lambda_{j_k+1}=0$. 
\end{enumerate}
\end{prop}

\begin{prop}\label{spprop3} For $l\in[0,\lfloor \frac{r}{2} \rfloor -2]$,
let $\mu''[l]$ be the following iteration of mutations 
\[\mu''[l]:=(\mu_{\lfloor \frac{r+1}{2} \rfloor-l-1}\mu_{r-l-1}\mu_{r-l})\cdots(\mu_{\lfloor \frac{r+1}{2} \rfloor-2}\mu_{r-2}\mu_{r-1})(\mu_{\lfloor \frac{r+1}{2} \rfloor -1}\mu_{r-1}\mu_{r})\mu_{\lfloor \frac{r+1}{2} \rfloor}.
\]
\begin{enumerate}
\item[$(a)$]We have the cluster variable
\begin{eqnarray*}
& &\hspace{-20pt}(\varphi^G_{(\mu''[l]\mathbb{V})})_{\lfloor \frac{r+1}{2} \rfloor -l-1}(a;\textbf{Y})=a^{\sum^{l+1}_{s=0}\Lambda_{2s+3}+\sum^{l}_{s=0}\Lambda_{2s+2}}\cdot H_5
\left(\prod_{q\in[1,l+1]}(1+A_{1,2q}^{-1})
\right. \\
&+&\left.
\sum_{p>0}\sum_{({\rm \bf{b}},{\rm \bf{c}})\in R^p_{l+1}, b_1>0}
\prod^p_{i=1} A[-c_i,-b_i;2]
\prod_{q\in[1,l+1]\setminus[{\rm \bf{b}},{\rm \bf{c}}]}(1+A_{1,2q}^{-1})\right).
\end{eqnarray*}
\item[$(b)$]We also obtain the cluster variable
\begin{eqnarray*}
& &\hspace{-40pt}(\varphi^G_{(\mu_{r-l-1}\mu''[l]\mathbb{V})})_{r-l-1}(a;\textbf{Y})=
a^{\sum^{l-1}_{s=0}\Lambda_{2s+3}+\sum^{l+1}_{s=0}\Lambda_{2s+2}} H_6 \\
&\times &\left(
(1+A^{-1}_{1,2l+4}+A^{-1}_{1,2l+4}A^{-1}_{2l+5}))
\prod_{q\in[1,l+1]}(1+A_{1,2q}^{-1})\right.\\
&+&\sum_{p>0}\sum_{({\rm \bf{b}},{\rm \bf{c}})\in R^p_{l+2},b_1>0}
(1-\delta_{c_p,l+2}+A^{-1+\delta_{c_p,l+2}}_{1,2l+4}(1+A^{-1}_{2l+5}))\\
&\times& \left.\prod^p_{i=1} A[-c_i,-b_i;2]
\prod_{q\in[1,l+1]\setminus[{\rm \bf{b}},{\rm \bf{c}}]}(1+A_{1,2q}^{-1})\right)
,\end{eqnarray*}
where 
$H_5:=\left(\prod^{l}_{t=0} Y_{1,2t+2}\right)\left(\prod^{l+1}_{t=0}Y_{2,2t+1}\right)$, $H_6:=\left(\prod^{l+1}_{t=0}Y_{1,2t+2}Y_{2,2t+1}\right)
=H_5\times Y_{1,2t+4}$.
\end{enumerate}
Furthermore, if $r$ is odd, then we get the cluster variable
\begin{eqnarray*}
& &\hspace{-20pt}(\varphi^G_{(\mu_1\mu_{\frac{r+3}{2}}\mu[\frac{r-1}{2}-2]\mathbb{V})})_{1}(a;\textbf{Y})\\
&=&a^{\sum^{\frac{r-5}{2}}_{s=0}\Lambda_{2s+3}+\sum^{\frac{r-3}{2}}_{s=0}\Lambda_{2s+2}}
\left(\prod^{\frac{r-3}{2}}_{t=0}Y_{1,2t+2}\right)
\left(\prod^{\frac{r-1}{2}}_{t=0}Y_{2,2t+1}\right)
\left(\prod_{q\in[1,\frac{r-1}{2}]}(1+A_{1,2q}^{-1})
\right.\\
&+&\left.
\sum_{p> 0}\sum_{({\rm \bf{b}},{\rm \bf{c}})\in R^p_{\frac{r-1}{2}},\ b_1>0}\prod^p_{i=1} A[-c_i,-b_i;2]
\prod_{q\in[1,\frac{r-1}{2}]\setminus[{\rm \bf{b}},{\rm \bf{c}}]}(1+A_{1,2q}^{-1})\right).
\end{eqnarray*}
\end{prop}

The following theorem is the main result, which means a relation between all the cluster variables in $\mathbb{C}[G^{e,c^2}]$ and Demazure crystals. We use the notation as in Proposition \ref{spprop1}, \ref{spprop2} and \ref{spprop3}.

\begin{thm}\label{thm1}
\begin{enumerate}
\item[$(1)$]Let $k\in[\lfloor \frac{r+1}{2} \rfloor+1,r-1]$ and $l\in[0,r-k-1]$.

$(a)$ The cluster variable
$(\varphi^G_{(\mu[l]\mathbb{V})})_{k-\lfloor \frac{r}{2} \rfloor +1}(a;\textbf{Y})$ is the total sum of monomial realizations of the Demazure crystals
\[ B\left(\sum^{j_k-1}_{s=j_k-2l-1}\Lambda_s\right)_{w_1}
\oplus \bigoplus_{\begin{array}{c}\vspace*{-3pt} _{p>0} \\ \vspace*{-5pt}
_{({\rm \bf{b}},{\rm \bf{c}})\in R^p_{l-1}} \end{array}}\hspace{-10pt}
B\left(\left(\sum^{j_k-1}_{s=j_k-2l-1}\Lambda_s\right)-\al[{\rm \bf{b}},{\rm \bf{c}};j_k]\right)_{w_1({\rm \bf{b}},{\rm \bf{c}})}, 
\]
where $w_1:=\prod_{q\in[0,l-1]}s_{j_k-2q-2}$, $w_1({\rm \bf{b}},{\rm \bf{c}}):=\prod_{q\in[0,l-1]\setminus[{\rm \bf{b}},{\rm \bf{c}}]}s_{j_k-2q-2}$, and the highest weight vectors in $B(\sum^{j_k-1}_{s=j_k-2l-1}\Lambda_s)$ and  $B((\sum^{j_k-1}_{s=j_k-2l-1}\Lambda_s)-\al[{\rm \bf{b}},{\rm \bf{c}};j_k])$ are realized by the monomials $H_1$ and $H_1\cdot
A[b_1,c_1;j_k]\cdots A[b_p,c_p;j_k]$, respectively.

$(b)$ The cluster variable
$(\varphi^G_{(\mu_{k+l+1}\mu[l]\mathbb{V})})_{k+l+1}(a;\textbf{Y})$ is the total sum of monomial realizations of the Demazure crystals
\[ 
B\left(\sum^{j_k-1}_{s=j_k-2l-2}\Lambda_s\right)_{w_2}\oplus\bigoplus_{\begin{array}{c}\vspace*{-3pt} _{p>0} \\ \vspace*{-5pt}
_{({\rm \bf{b}},{\rm \bf{c}})\in R^p_{l}} \end{array}}\hspace{-10pt}
B\left(\left(\sum^{j_k-1}_{s=j_k-2l-2}\Lambda_s\right)-\al[{\rm \bf{b}},{\rm \bf{c}};j_k]\right)_{w_2({\rm \bf{b}},{\rm \bf{c}})}, 
\]
where $w_2:=s_{j_k-2l-3}s_{j_k-2l-2}\prod_{q\in[0,l-1]}s_{j_k-2q-2}$, $w_2({\rm \bf{b}},{\rm \bf{c}}):=s_{j_k-2l-3}s^{1-\delta_{c_p,l}}_{j_k-2l-2}\prod_{q\in[0,l-1]\setminus[{\rm \bf{b}},{\rm \bf{c}}]}s_{j_k-2q-2}$, and the highest weight vectors in $B(\sum^{j_k-1}_{s=j_k-2l-2}\Lambda_s)$ and $B((\sum^{j_k-1}_{s=j_k-2l-2}\Lambda_s)-\al[{\rm \bf{b}},{\rm \bf{c}};j_k])$ are realized by the monomials $H_2$ and $H_2\cdot
A[b_1,c_1;j_k]\cdots A[b_p,c_p;j_k]$, respectively. 

\item[$(2)$]Let $k\in[1,\lfloor \frac{r+1}{2} \rfloor-2]$ and $l\in[0,\lfloor \frac{r+1}{2} \rfloor -k-2]$.

$(a)$ If $j_k<r$, the cluster variable
$(\varphi^G_{(\mu'[l]\mathbb{V})})_{k+l+1}(a;\textbf{Y})$ is the total sum of monomial realizations of the Demazure crystals
\[ B\left(\sum^{j_k+1}_{s=j_k-2l-2}\Lambda_s\right)_{w_3}
\oplus
\bigoplus_{\begin{array}{c}\vspace*{-3pt} _{p>0} \\ \vspace*{-5pt}
_{({\rm \bf{b}},{\rm \bf{c}})\in R^p_{l+1}} \end{array}}\hspace{-10pt}
B\left(\left(\sum^{j_k+1}_{s=j_k-2l-2}\Lambda_s\right)-\al[{\rm \bf{b}},{\rm \bf{c}};j_k+3]\right)_{w_3({\rm \bf{b}},{\rm \bf{c}})}, 
\]
if $j_k=r$, $(\varphi^G_{(\mu'[l]\mathbb{V})})_{k+l+1}(a;\textbf{Y})$ is the total sum of monomial realizations of the Demazure crystals
\[ 
B\left(\sum^{r}_{s=r-2l-2}\Lambda_s\right)_{w_3}
\oplus
\bigoplus_{\begin{array}{c}\vspace*{-3pt} _{p>0} \\ \vspace*{-5pt}
_{({\rm \bf{b}},{\rm \bf{c}})\in R^p_{l+1},b_1>0} \end{array}}\hspace{-10pt}
B\left(\left(\sum^{r}_{s=r-2l-2}\Lambda_s\right)-\al[{\rm \bf{b}},{\rm \bf{c}};r+3]\right)_{w_3({\rm \bf{b}},{\rm \bf{c}})}, 
\]
where $w_3:=s_{j_k+2}s_{j_k+1}\prod_{q\in[1,l+1]}s_{j_k-2q+1}$, $w_3({\rm \bf{b}},{\rm \bf{c}}):=s_{j_k+2}s^{1-\delta_{b_1,0}}_{j_k+1}\prod_{q\in[1,l+1]\setminus[{\rm \bf{b}},{\rm \bf{c}}]}s_{j_k-2q+1}$, and the highest weight vectors in $B(\sum^{j_k+1}_{s=j_k-2l-2}\Lambda_s)$ and $B((\sum^{j_k+1}_{s=j_k-2l-2}\Lambda_s)-\al[{\rm \bf{b}},{\rm \bf{c}};j_k+3])$ are realized by the monomials $H_3$ and $H_3\cdot
A[b_1,c_1;j_k+3]\cdots A[b_p,c_p;j_k+3]$, respectively.

$(b)$ If $j_k<r$, the cluster variable
$(\varphi^G_{(\mu_{\lfloor \frac{r}{2} \rfloor +k+l+2}\mu'[l]\mathbb{V})})_{\lfloor \frac{r}{2} \rfloor +k+l+2}(a;\textbf{Y})$ is the total sum of monomial realizations of the Demazure crystals
\[ 
B\left(\sum^{j_k+1}_{s=j_k-2l-3}\Lambda_s\right)_{w_4}
\oplus
\bigoplus_{\begin{array}{c}\vspace*{-3pt} _{p>0} \\ \vspace*{-5pt}
_{({\rm \bf{b}},{\rm \bf{c}})\in R^p_{l+2}} \end{array}}\hspace{-10pt}
B\left(\left(\sum^{j_k+1}_{s=j_k-2l-3}\Lambda_s\right)-\al[{\rm \bf{b}},{\rm \bf{c}};j_k+3]\right)_{w_4({\rm \bf{b}},{\rm \bf{c}})}, 
\]
if $j_k=r$, $(\varphi^G_{(\mu_{\lfloor \frac{r}{2} \rfloor +k+l+2}\mu'[l]\mathbb{V})})_{\lfloor \frac{r}{2} \rfloor +k+l+2}(a;\textbf{Y})$ is the total sum of monomial realizations of the Demazure crystals
\[ 
B\left(\sum^{r}_{s=r-2l-3}\Lambda_s\right)_{w_4}
\oplus
\bigoplus_{\begin{array}{c}\vspace*{-3pt} _{p>0} \\ \vspace*{-5pt}
_{({\rm \bf{b}},{\rm \bf{c}})\in R^p_{l+2},b_1>0} \end{array}}\hspace{-10pt}
B\left(\left(\sum^{r}_{s=r-2l-3}\Lambda_s\right)-\al[{\rm \bf{b}},{\rm \bf{c}};r+3]\right)_{w_4({\rm \bf{b}},{\rm \bf{c}})}, 
\]
where $w_4:=s_{j_k+2}s_{j_k+1}s_{j_k-2l-4}s_{j_k-2l-3}\prod_{q\in[1,l+1]}s_{j_k-2q+1}$, $w_4({\rm \bf{b}},{\rm \bf{c}}):=s_{j_k+2}s^{1-\delta_{b_1,0}}_{j_k+1}s_{j_k-2l-4}s^{1-\delta_{c_p,l+2}}_{j_k-2l-3}\prod_{q\in[1,l+1]\setminus[{\rm \bf{b}},{\rm \bf{c}}]}s_{j_k-2q+1}$, and the highest weight vectors in $B(\sum^{j_k+1}_{s=j_k-2l-3}\Lambda_s)$ and $B((\sum^{j_k+1}_{s=j_k-2l-3}\Lambda_s)-\al[{\rm \bf{b}},{\rm \bf{c}};j_k+3])$ are realized by the monomials $H_4$ and $H_4\cdot
A[b_1,c_1;j_k+3]\cdots A[b_p,c_p;j_k+3]$, respectively. If $j_k=r$, then we understand $Y_{1,j_k+1}=1$, $\Lambda_{j_k+1}=0$ and $s_{j_k+1}=s_{j_k+2}=e$.

\item[$(3)$] Let $l\in[0,\lfloor \frac{r}{2} \rfloor -2]$.

$(a)$ The cluster variable
$(\varphi^G_{(\mu''[l]\mathbb{V})})_{\lfloor \frac{r+1}{2} \rfloor -l-1}(a;\textbf{Y})$ is the total sum of monomial realizations of the Demazure crystals
\[ 
B\left(\sum^{2l+3}_{s=1}\Lambda_s\right)_{w_5}
 \oplus \hspace{-12pt}
\bigoplus_{\begin{array}{c}\vspace*{-3pt} _{p>0} \\ \vspace*{-5pt}
_{({\rm \bf{b}},{\rm \bf{c}})\in R^p_{l+1},b_1>0} \end{array}}\hspace{-20pt}
B\left(\left(\sum^{2l+3}_{s=1}\Lambda_s\right)-\al[-{\rm \bf{c}},-{\rm \bf{b}};2]\right)_{w_5({\rm \bf{b}},{\rm \bf{c}})}, 
\]
where $w_5:=\prod_{q\in[1,l+1]}s_{2q}$, $w_5({\rm \bf{b}},{\rm \bf{c}}):=\prod_{q\in[1,l+1]\setminus[{\rm \bf{b}},{\rm \bf{c}}]}s_{2q}$ and the highest weight vectors in $B(\sum^{2l+3}_{s=1}\Lambda_s)$ and $B((\sum^{2l+3}_{s=1}\Lambda_s)-\al[-{\rm \bf{c}},-{\rm \bf{b}};2])$ are realized by the monomials $H_5$ and $H_5\cdot
A[-c_1,-b_1;2]\cdots A[-c_p,-b_p;2]$, respectively.

$(b)$ The cluster variable
$(\varphi^G_{(\mu_{r-l-1}\mu''[l]\mathbb{V})})_{r-l-1}(a;\textbf{Y})$ is the total sum of monomial realizations of the Demazure crystals
\[ 
B\left(\sum^{2l+4}_{s=1}\Lambda_s\right)_{w_6}
\oplus \hspace{-12pt}
\bigoplus_{\begin{array}{c}\vspace*{-3pt} _{p>0} \\ \vspace*{-5pt}
_{({\rm \bf{b}},{\rm \bf{c}})\in R^p_{l+2},b_1>0} \end{array}}\hspace{-20pt}
B\left(\left(\sum^{2l+4}_{s=1}\Lambda_s\right)-\al[-{\rm \bf{c}},-{\rm \bf{b}};2]\right)_{w_6({\rm \bf{b}},{\rm \bf{c}})}, 
\]
where $w_6:=s_{2l+5}s_{2l+4}\prod_{q\in[1,l+1]}s_{2q}$, $w_6({\rm \bf{b}},{\rm \bf{c}}):=s_{2l+5}s^{1-\delta_{c_p,l+2}}_{2l+4}\prod_{q\in[1,l+1]\setminus[{\rm \bf{b}},{\rm \bf{c}}]}s_{2q}$ and the highest weight vectors in $B(\sum^{2l+4}_{s=1}\Lambda_s)$ and $B((\sum^{2l+4}_{s=1}\Lambda_s)-\al[-{\rm \bf{c}},-{\rm \bf{b}};2])$ are realized by the monomials $H_6$ and $H_6\cdot
A[-c_1,-b_1;2]\cdots A[-c_p,-b_p;2]$, respectively.

Furthermore, if $r$ is odd, then the cluster variable
$(\varphi^G_{(\mu_1\mu_{\frac{r+3}{2}}\mu[\frac{r-1}{2}-2]\mathbb{V})})_{1}(a;\textbf{Y})$ is the total sum of monomial realizations of the Demazure crystals
\[
B\left(\sum^{r}_{s=1}\Lambda_s\right)_{\prod_{q\in[1,\frac{r-1}{2}]}s_{2q}}
\oplus \hspace{-12pt}
\bigoplus_{\begin{array}{c}\vspace*{-3pt} _{p>0} \\ \vspace*{-5pt}
_{({\rm \bf{b}},{\rm \bf{c}})\in R^p_{\frac{r-1}{2}},b_1>0} \end{array}}\hspace{-20pt}
B\left(\left(\sum^{r}_{s=1}\Lambda_s\right)-\al[-{\rm \bf{c}},-{\rm \bf{b}};2]\right)_{\prod_{q\in[1,\frac{r-1}{2}]\setminus[{\rm \bf{b}},{\rm \bf{c}}]}s_{2q}}, 
\]
where the highest weight vectors in $B(\sum^{r}_{s=1}\Lambda_s)$ and $B((\sum^{r}_{s=1}\Lambda_s)-\al[-{\rm \bf{c}},-{\rm \bf{b}};2])$ are realized by the monomials 
$\prod^{\frac{r-3}{2}}_{t=0}Y_{1,2t+2}
\prod^{\frac{r-1}{2}}_{t=0}Y_{2,2t+1}$ and
$\prod^{\frac{r-3}{2}}_{t=0}Y_{1,2t+2}
\prod^{\frac{r-1}{2}}_{t=0}Y_{2,2t+1}\cdot
A[-c_1,-b_1;2]\cdots A[-c_p,-b_p;2]$, respectively. 
\end{enumerate}
\end{thm}

We obtain the following theorem from Example \ref{initialex4}, Proposition \ref{remainexp} and Theorem \ref{thm1}. Let $\Xi$ be the set of the non-frozen cluster variables in $\mathbb{C}[G^{e,c^2}]$.
\begin{thm}\label{maincor}
Each initial cluster variable $\varphi_{\mathbb{V}_k}$ in $\mathbb{C}[G^{e,c^2}]$ is the total sum of monomials in the Demazure crystal $B(\Lambda_{j_k})_{c^2_{>2r-k}}$, where we use the notation as in (\ref{inc}).
For each non-initial cluster variable $\varphi$ in $\mathbb{C}[G^{e,c^2}]$, there uniquely exist $p\geq0$, $w,w[i]\in W$, $\lambda:=\sum^{b}_{j=a}\Lambda_{j}$ $(1\leq a\leq b\leq r)$ and $\lambda_i\in P^+$ such that $\lambda-\lambda_i\in \oplus_{s\in I}\mathbb{Z}_{\geq0}\al_s$ and $\varphi$ is the total sum of monomials in Demazure crystals in the form
\[ B(\lambda)_w\oplus\bigoplus^p_{i=1}B(\lambda_{i})_{w[i]}. \]
Then, let $\varphi_{\lambda}$ denote this non-initial cluster variable $\varphi$.

In particular, the set $\{\varphi_{\lambda}\}_{\lambda\in \Phi_{\geq -1}}$ exhausts all the cluster variables in $C[G^{e,c^2}]$. More precisely, the map $\Phi_{\geq -1}\rightarrow \Xi$,
\[ -\al_{j_k}\mapsto \varphi_{\mathbb{V}_k},\ \ \sum^b_{j=a}\al_{j}\mapsto \varphi_{\sum^b_{j=a}\Lambda_j} \]
is a bijection between the set $\Phi_{\geq -1}$ of almost positive roots and $\Xi$.
\end{thm}

\begin{rem}
The correspondence between $\Phi_{\geq -1}$ and the set $\Xi$ of cluster variables in Theorem \ref{maincor} is different from the one of \cite{FZ3}. 
\end{rem}

\begin{ex}\label{r10ex2}
We consider the same setting as in Example \ref{r10ex}. Let $\mu$ (resp. $\mu'$) denote the monomial realization of crystal $B(\Lm_5+\Lm_6+\Lm_7+\Lm_8+\Lm_9)$ (resp. $B(2\Lm_5+\Lm_7+2\Lm_9)$) such that the highest weight vector is realized as $Y_{2,5}Y_{1,6}Y_{2,7}Y_{1,8}Y_{2,9}$ (resp. $Y_{1,5}Y_{2,5}Y_{1,7}Y_{1,9}Y_{2,9}$). It follows from Theorem \ref{monorealmain}, \ref{kashidem} and (\ref{r10ex-1}) that
\begin{eqnarray*}
& &\hspace{-20pt}(\varphi^G_{(\mu[2]\mathbb{V})})_{3}(a;\textbf{Y})\\
&=&a^{\Lm_3+\Lm_5+\Lm_6+\Lm_7+\Lm_8+\Lm_9}\left(\sum_{b\in B(\Lm_5+\Lm_6+\Lm_7+\Lm_8+\Lm_9)_{s_6s_8}} \mu(b) 
+\sum_{b\in B(2\Lm_5+\Lm_7+2\Lm_9)_{e}} \mu'(b) \right).
\end{eqnarray*}
Similarly, using (\ref{r10ex-2}),
\begin{eqnarray*}
& &a^{-(\Lm_4+\Lm_6+\Lm_7+\Lm_8+\Lm_9)}(\varphi^G_{(\mu_9\mu[2]\mathbb{V})})_{9}(a;\textbf{Y})\nonumber\\
&=&\sum_{b\in B(\Lm_4+\Lm_5+\Lm_6+\Lm_7+\Lm_8+\Lm_9)_{s_3s_4s_6s_8}}\mu(b)
+\sum_{b\in B(\Lm_4+2\Lm_5+\Lm_7+2\Lm_9)_{s_3s_4}}\mu'(b)\nonumber\\
& &+\sum_{b\in B(\Lm_3+\Lm_5+\Lm_6+\Lm_7+2\Lm_9)_{s_3}}\mu''(b)
+\sum_{b\in B(\Lm_3+\Lm_5+2\Lm_7+\Lm_8+\Lm_9)_{s_3}}\mu'''(b),\nonumber
\end{eqnarray*}
where $\mu$, $\mu'$, $\mu''$ and $\mu'''$ are the monomial realizations such that the highest weight vectors are realized by $Y_{1,4}Y_{2,5}Y_{1,6}Y_{2,7}Y_{1,8}Y_{2,9}$, $Y_{1,4}Y_{1,5}Y_{2,5}Y_{1,7}Y_{1,9}Y_{2,9}$, $Y_{1,3}Y_{1,5}Y_{2,6}Y_{1,7}Y_{1,9}Y_{2,9}$ and $Y_{1,3}Y_{1,5}Y_{1,7}Y_{2,7}Y_{1,8}Y_{2,9}$, respectively.
\end{ex}

\section{The proof of main theorem}

In this section, we prove Proposition \ref{spprop1}, \ref{spprop2}, \ref{spprop3} and Theorem \ref{thm1}. For $k_1,\cdots,k_s\in[1,r]$, let $\mu_{k_1}\cdots\mu_{k_s}\Gamma_{{\rm \bf{i}}}$ be the quiver of the seed $(\mu_{k_1}\cdots\mu_{k_s}(\mathbb{V}),\mu_{k_1}\cdots\mu_{k_s}(\tilde{B}_{{\rm \bf{i}}}))$ (Sect.\ref{CluSect}).

\begin{lem}\label{mut0}
In the setting of Proposition \ref{spprop1}, we have
\[
(\varphi^G_{(\mu_{k+l+2}\mu_{k+l+1}\mu[l]\mathbb{V})})_{k+l+2}(a;\textbf{Y})
=(\varphi^G_{(\mu_{k+l+2}\mathbb{V})})_{k+l+2}(a;\textbf{Y}).
\]
\end{lem}
\nd
{\sl [Proof.]}

In the quiver $\Gamma_{{\rm \bf{i}}}$, by (\ref{iniquiver}), the vertices and arrows around $(\varphi_{\mathbb{V}})_{k+l+2}$ are
\[
\begin{xy}
(115,90)*{\cdots}="emp1",
(100,100)*{(\varphi_{\mathbb{V}})_{r+k+l+1}}="7",
(100,90) *{(\varphi_{\mathbb{V}})_{k+l+1}}="3",
(100,80)*{(\varphi_{\mathbb{V}})_{-j_k-2l-2}}="-4",
(70,100)*{(\varphi_{\mathbb{V}})_{r+k-\lfloor \frac{r}{2} \rfloor+l+1}}="r+k",
(70,90) *{(\varphi_{\mathbb{V}})_{k-\lfloor \frac{r}{2} \rfloor+l+1}}="k",
(70,80)*{(\varphi_{\mathbb{V}})_{-j_k-2l-3}}="-j_k",
(40,100)*{(\varphi_{\mathbb{V}})_{r+k+l+2}}="8",
(40,90) *{(\varphi_{\mathbb{V}})_{k+l+2}}="4",
(40,80)*{(\varphi_{\mathbb{V}})_{-j_k-2l-4}}="-2",
(10,100)*{(\varphi_{\mathbb{V}})_{r+k-\lfloor \frac{r}{2} \rfloor+l+2}}="6",
(10,90) *{(\varphi_{\mathbb{V}})_{k-\lfloor \frac{r}{2} \rfloor+l+2}}="2",
(10,80)*{(\varphi_{\mathbb{V}})_{-j_k-2l-5}}="-1",
(-5,90)*{\cdots}="emp",
\ar@{->} "7";"3"
\ar@{->} "r+k";"k"
\ar@{->} "8";"4"
\ar@{->} "6";"2"
\ar@{->} "3";"-4"
\ar@{->} "k";"-j_k"
\ar@{->} "4";"-2"
\ar@{->} "2";"-1"
\ar@{->} "k";"3"
\ar@{->} "k";"4"
\ar@{->} "2";"4"
\ar@{->} "3";"r+k"
\ar@{->} "4";"r+k"
\ar@{->} "4";"6"
\ar@{->} "-4";"k"
\ar@{->} "-2";"k"
\ar@{->} "-2";"2"
\end{xy}
\]
The initial cluster variables changed by $\mu_{k+l+1}\mu[l]$ are $(\varphi_{\mathbb{V}})_{k}, (\varphi_{\mathbb{V}})_{k+1},\cdots,(\varphi_{\mathbb{V}})_{k+l+1}$ and $(\varphi_{\mathbb{V}})_{k-\lfloor \frac{r}{2} \rfloor}$,$(\varphi_{\mathbb{V}})_{k-\lfloor \frac{r}{2} \rfloor+1},\cdots,(\varphi_{\mathbb{V}})_{k-\lfloor \frac{r}{2} \rfloor+l}$, which are not connected with $(\varphi_{\mathbb{V}})_{k+l+2}$ in the above quiver. Therefore, Lemma \ref{mutgamlem} says that the arrows incident to $(\varphi_{\mathbb{V}})_{k+l+2}$ in $\Gamma_{{\rm \bf{i}}}$ coincide with the ones in $\mu_{k+l+1}\mu[l]\Gamma_{{\rm \bf{i}}}$. Therefore, we get
\begin{eqnarray*}
& &\hspace{-20pt}(\varphi^G_{(\mu_{k+l+2}\mu_{k+l+1}\mu[l]\mathbb{V})})_{k+l+2}\\
&=&\frac{1}{(\varphi^G_{\mathbb{V}})_{k+l+2}}\left(
(\varphi^G_{\mathbb{V}})_{r+k-\lfloor \frac{r}{2} \rfloor+l+1}
(\varphi^G_{\mathbb{V}})_{r+k-\lfloor \frac{r}{2} \rfloor+l+2}
(\varphi^G_{\mathbb{V}})_{-j_k-2l-4}\right.\\
& &\left.+
(\varphi^G_{\mathbb{V}})_{k-\lfloor \frac{r}{2} \rfloor+l+1}
(\varphi^G_{\mathbb{V}})_{k-\lfloor \frac{r}{2} \rfloor+l+2}
(\varphi^G_{\mathbb{V}})_{r+k+l+2}\right)
=(\varphi^G_{(\mu_{k+l+2}\mathbb{V})})_{k+l+2}.
\end{eqnarray*}
\qed

Next, we will order the indecomposable direct summands $V_1,\cdots,V_{2r}$ of $V_{{\rm \bf{i}}}$ from the right:
\[ V_{{\rm \bf{i}}}=V_{2r}\oplus \cdots \oplus V_1. \]
For a basic $\mathcal{C}_{c^2}$-cluster-tilting $\Lambda$-module $T=T_{2r}\oplus \cdots \oplus T_1$, we write
\[ \mu_k(T):=\mu_{T_k}(T)=T_{2r}\oplus \cdots\oplus T_{k+1}\oplus T^*_k\oplus T_{k-1}\oplus\cdots \oplus T_1, \] 
for $k\in[1,r]$. Let $(\mu_k(T))_{l}$ denote the $l$-th indecomposable direct summand of $\mu_k(T)$ from the right. 

In the following Lemma \ref{mutlem1}-\ref{mutlem3}, the notation as in Remark \ref{diarem} is applied.

\begin{lem}\label{mutlem1}
We use the notation as in Proposition \ref{spprop1} and let $j_k$ be the $k$-th index of {\rm \bf{i}} in $(\ref{redwords2})$ from the right. 
\begin{enumerate}
\item[$(a)$] The module $(\mu[l](V_{\rm \bf{i}}))_{k-\lfloor \frac{r}{2} \rfloor +l}$ is described as follows:
\begin{equation}\label{mutlemclaim1}
\begin{xy}
(63,123)*{_{j_k-3}}="a6",(71,123)*{_{j_k-1}}="a7",(79,123)*{_{j_k+1}}="a8",
(87,123)*{_{j_k+3}}="a9",
(67,113)*{_{j_k-2}}="a3",(75,113)*{_{j_k}}="a4",(83,113)*{_{j_k+2}}="a5",
(71,103)*{_{j_k-1}}="a1",(79,103)*{ _{j_k+1}}="a2",
(32,123)*{_{j_k-5}}="b6",(40,123)*{_{j_k-3}}="b7",(48,123)*{_{j_k-1}}="b8",
(56,123)*{_{j_k+1}}="b9",
(36,113)*{_{j_k-4}}="b3",(44,113)*{_{j_k-2}}="b4",(52,113)*{_{j_k}}="b5",
(40,103)*{_{j_k-3}}="b1",
(-17,123)*{_{j_{k+l}-5}}="c6",(-7,123)*{_{j_{k+l}-3}}="c7",(3,123)*{_{j_{k+l}-1}}="c8",
(13,123)*{_{j_{k+l}+1}}="c9",
(-12,113)*{_{j_{k+l}-4}}="c3",(-2,113)*{_{j_{k+l}-2}}="c4",(8,113)*{_{j_{k+l}}}="c5",
(-7,103)*{_{j_{k+l}-3}}="c1",(22,103)*{_{j_{k+l}-1}}="z1",(20,113)*{\cdots}="dot",
(26,113)*{\ }="hs",
\ar@{->} "hs";"b1"
\ar@{->} "dot";"b1"
\ar@{->} "a3";"a1"
\ar@{->} "a4";"a1"
\ar@{->} "a4";"a2"
\ar@{->} "a5";"a2"
\ar@{->} "a6";"a3"
\ar@{->} "a7";"a3"
\ar@{->} "a7";"a4"
\ar@{->} "a8";"a4"
\ar@{->} "a8";"a5"
\ar@{->} "a9";"a5"
\ar@{->} "b3";"b1"
\ar@{->} "b4";"b1"
\ar@{->} "b4";"a1"
\ar@{->} "b5";"a1"
\ar@{->} "b6";"b3"
\ar@{->} "b7";"b3"
\ar@{->} "b7";"b4"
\ar@{->} "b8";"b4"
\ar@{->} "b8";"b5"
\ar@{->} "b9";"b5"
\ar@{->} "c3";"c1"
\ar@{->} "c4";"c1"
\ar@{->} "c4";"z1"
\ar@{->} "c5";"z1"
\ar@{->} "c6";"c3"
\ar@{->} "c7";"c3"
\ar@{->} "c7";"c4"
\ar@{->} "c8";"c4"
\ar@{->} "c8";"c5"
\ar@{->} "c9";"c5"
\end{xy}
\end{equation}

Note that $j_{k+l}=j_k-2l$ from $(\ref{redwords2})$.

\item[$(b)$] The module $(\mu_{k+l+1}\mu[l](V_{\rm \bf{i}}))_{k+l+1}$ is described as follows:
\begin{equation}\label{mutlemclaim2}
\begin{xy}
(63,123)*{_{j_k-3}}="a6",(71,123)*{_{j_k-1}}="a7",(79,123)*{_{j_k+1}}="a8",
(87,123)*{_{j_k+3}}="a9",
(67,113)*{_{j_k-2}}="a3",(75,113)*{_{j_k}}="a4",(83,113)*{_{j_k+2}}="a5",
(71,103)*{_{j_k-1}}="a1",(79,103)*{ _{j_k+1}}="a2",
(32,123)*{_{j_k-5}}="b6",(40,123)*{_{j_k-3}}="b7",(48,123)*{_{j_k-1}}="b8",
(56,123)*{_{j_k+1}}="b9",
(36,113)*{_{j_k-4}}="b3",(44,113)*{_{j_k-2}}="b4",(52,113)*{_{j_k}}="b5",
(40,103)*{_{j_k-3}}="b1",
(-17,123)*{_{j_{k+l}-3}}="c6",(-7,123)*{_{j_{k+l}-1}}="c7",(3,123)*{_{j_{k+l}+1}}="c8",
(13,123)*{_{j_{k+l}+3}}="c9",
(-12,113)*{_{j_{k+l}-2}}="c3",(-2,113)*{_{j_{k+l}}}="c4",(8,113)*{_{j_{k+l}+2}}="c5",
(22,103)*{_{j_{k+l}+1}}="z1",(20,113)*{\cdots}="dot",
(26,113)*{\ }="hs",
\ar@{->} "hs";"b1"
\ar@{->} "dot";"b1"
\ar@{->} "a3";"a1"
\ar@{->} "a4";"a1"
\ar@{->} "a4";"a2"
\ar@{->} "a5";"a2"
\ar@{->} "a6";"a3"
\ar@{->} "a7";"a3"
\ar@{->} "a7";"a4"
\ar@{->} "a8";"a4"
\ar@{->} "a8";"a5"
\ar@{->} "a9";"a5"
\ar@{->} "b3";"b1"
\ar@{->} "b4";"b1"
\ar@{->} "b4";"a1"
\ar@{->} "b5";"a1"
\ar@{->} "b6";"b3"
\ar@{->} "b7";"b3"
\ar@{->} "b7";"b4"
\ar@{->} "b8";"b4"
\ar@{->} "b8";"b5"
\ar@{->} "b9";"b5"
\ar@{->} "c4";"z1"
\ar@{->} "c5";"z1"
\ar@{->} "c6";"c3"
\ar@{->} "c7";"c3"
\ar@{->} "c7";"c4"
\ar@{->} "c8";"c4"
\ar@{->} "c8";"c5"
\ar@{->} "c9";"c5"
\end{xy}
\end{equation}
\end{enumerate}
\end{lem}
\nd
{\sl [Proof.]}

Using the induction on $l$, we shall prove (a) and (b) simultaneously. 

First, let us prove (a) and (b) for $l=0$. As have seen in Example \ref{initialex1} (\ref{iniex1-1}), $(\mu_{k+1}\mu_k(V_{\rm \bf{i}}))_{k-\lfloor \frac{r}{2} \rfloor}=(V_{\rm \bf{i}})_{k-\lfloor \frac{r}{2} \rfloor}=S_{j_k-1}$. We have already obtained $(\mu_k(V_{\rm \bf{i}}))_k=V^*_k$ in Example \ref{initialex2} (\ref{mutex2}). Similarly, $(\mu_{k+1}\mu_k(V_{\rm \bf{i}}))_{k+1}$ is

\begin{equation}\label{mutlem1-pr1}
\begin{xy}
(35,123)*{j_k-5}="4r-2k-1",(55,123)*{j_k-3}="4r-2k+1b",(75,123)*{j_k-1}="4r-2k+3b",
(95,123)*{j_k+1}="4r-2k+5",
(45,118)*{j_k-4}="4r-2k",(65,118)*{j_k-2}="4r-2k+2a",(85,118)*{j_k}="4r-2k+4",
(55,113)*{j_k-3}="4r-2k+1",(75,113)*{j_k-1}="4r-2k+3",
(35,113)*{}="12",
\ar@{->} "4r-2k";"4r-2k+1"
\ar@{->} "4r-2k+2a";"4r-2k+1"
\ar@{->} "4r-2k+2a";"4r-2k+3"
\ar@{->} "4r-2k+4";"4r-2k+3"
\ar@{->} "4r-2k-1";"4r-2k"
\ar@{->} "4r-2k+1b";"4r-2k"
\ar@{->} "4r-2k+1b";"4r-2k+2a"
\ar@{->} "4r-2k+3b";"4r-2k+2a"
\ar@{->} "4r-2k+3b";"4r-2k+4"
\ar@{->} "4r-2k+5";"4r-2k+4"
\end{xy}
\end{equation}

Hence, the modules $(\mu_k(V_{\rm \bf{i}}))_k$ and $(\mu_{k+1}\mu_k(V_{\rm \bf{i}}))_{k+1}$ have the simple submodule isomorphic to $S_{j_k-1}$. So there exist injective homomorphisms 
\[(\mu_{k+1}\mu_k(V_{\rm \bf{i}}))_{k-\lfloor \frac{r}{2} \rfloor}=S_{j_k-1}\rightarrow (\mu_k(V_{\rm \bf{i}}))_k,\quad (\mu_{k+1}\mu_k(V_{\rm \bf{i}}))_{k-\lfloor \frac{r}{2} \rfloor} \rightarrow (\mu_{k+1}\mu_k(V_{\rm \bf{i}}))_{k+1}.\]
 Let $e_{j_k-1}$ denote a basis vector in $(\mu_{k+1}\mu_k(V_{\rm \bf{i}}))_{k-\lfloor \frac{r}{2} \rfloor}=S_{j_k-1}$, and let $e'_{j_k-1}\in (\mu_k(V_{\rm \bf{i}}))_k$ and $e''_{j_k-1}\in (\mu_{k+1}\mu_k(V_{\rm \bf{i}}))_{k+1}$ be the images of $e_{j_k-1}$ respectively. Note that since $j_{r+k-\lfloor \frac{r}{2} \rfloor}=j_{k-\lfloor \frac{r}{2} \rfloor}=j_k-1$, the module
$V_{r+k-\lfloor \frac{r}{2} \rfloor}$ is described as

\begin{equation}\label{mutlem1-pr1b}
\begin{xy}
(65,123)*{j_k-3}="4r-2k+1b",(85,123)*{j_k-1}="4r-2k+3b",
(105,123)*{j_k+1}="4r-2k+5",
(75,118)*{j_k-2}="4r-2k+2a",(95,118)*{j_k}="4r-2k+4",
(85,113)*{j_k-1}="4r-2k+3",
(35,113)*{}="12",
\ar@{->} "4r-2k+2a";"4r-2k+3"
\ar@{->} "4r-2k+4";"4r-2k+3"
\ar@{->} "4r-2k+1b";"4r-2k+2a"
\ar@{->} "4r-2k+3b";"4r-2k+2a"
\ar@{->} "4r-2k+3b";"4r-2k+4"
\ar@{->} "4r-2k+5";"4r-2k+4"
\end{xy}
\end{equation}
and has the simple socle isomorphic to $S_{j_k-1}$ (\ref{iniex1-3}). So there exists an injective homomorphism $S_{j_k-1}\rightarrow V_{r+k-\lfloor \frac{r}{2} \rfloor}$. However, this map is factorizable in the direct summands of $\mu_{k+1}\mu_k(V_{\rm \bf{i}})$ since it 
is the same as the composite map $S_{j_k-1}\rightarrow (\mu_{k+1}\mu_k(V_{\rm \bf{i}}))_k=(\mu_k(V_{\rm \bf{i}}))_k \rightarrow V_{r+k-\lfloor \frac{r}{2} \rfloor}$. Moreover, we can verify that Hom$(S_{j_k-1},V_t)=\{0\}$ for $t\neq r+k-\lfloor \frac{r}{2} \rfloor$, $k-\lfloor \frac{r}{2} \rfloor$, $k$, $k+1$. 
From Lemma \ref{titj} and Theorem \ref{GLSthm} (iii), the exchange sequence associated to the direct summand $S_{j_k-1}$ of 
$\mu_{k+1}\mu_k(V_{\rm \bf{i}})$ is as follows:
\[ 0 \rightarrow S_{j_k-1}\rightarrow (\mu_k(V_{\rm \bf{i}}))_k\oplus (\mu_{k+1}\mu_k(V_{\rm \bf{i}}))_{k+1}\rightarrow (\mu_{k-\lfloor \frac{r}{2} \rfloor}\mu_{k+1}\mu_k(V_{\rm \bf{i}}))_{k-\lfloor \frac{r}{2} \rfloor}\rightarrow 0, \]
where the image of the injective homomorphism $S_{j_k-1}\rightarrow (\mu_k(V_{\rm \bf{i}}))_k\oplus (\mu_{k+1}\mu_k(V_{\rm \bf{i}}))_{k+1}$ is $\mathbb{C}(e'_{j_k-1}+ e''_{j_k-1})$. 
Therefore, the module 
\begin{equation}\label{mutlem1-pr-seq1}
(\mu_{k-\lfloor \frac{r}{2} \rfloor}\mu_{k+1}\mu_k(V_{\rm \bf{i}}))_{k-\lfloor \frac{r}{2} \rfloor}=((\mu_k(V_{\rm \bf{i}}))_k\oplus (\mu_{k+1}\mu_k(V_{\rm \bf{i}}))_{k+1})/\mathbb{C}(e'_{j_k-1}+ e''_{j_k-1})
\end{equation}
is described as follows:
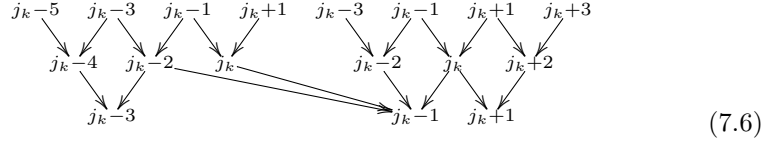
\begin{equation}\label{mutlem1-pr2}
\begin{xy}
(75,123)*{_{j_k-3}}="a6",(85,123)*{_{j_k-1}}="a7",(95,123)*{_{j_k+1}}="a8",
(105,123)*{_{j_k+3}}="a9",
(80,116)*{_{j_k-2}}="a3",(90,116)*{_{j_k}}="a4",(100,116)*{_{j_k+2}}="a5",
(85,109)*{_{j_k-1}}="a1",(95,109)*{_{j_k+1}}="a2",
(35,123)*{_{j_k-5}}="b6",(45,123)*{_{j_k-3}}="b7",(55,123)*{_{j_k-1}}="b8",
(65,123)*{_{j_k+1}}="b9",
(40,116)*{_{j_k-4}}="b3",(50,116)*{_{j_k-2}}="b4",(60,116)*{_{j_k}}="b5",
(45,109)*{_{j_k-3}}="b1",
\ar@{->} "a3";"a1"
\ar@{->} "a4";"a1"
\ar@{->} "a4";"a2"
\ar@{->} "a5";"a2"
\ar@{->} "a6";"a3"
\ar@{->} "a7";"a3"
\ar@{->} "a7";"a4"
\ar@{->} "a8";"a4"
\ar@{->} "a8";"a5"
\ar@{->} "a9";"a5"
\ar@{->} "b3";"b1"
\ar@{->} "b4";"b1"
\ar@{->} "b4";"a1"
\ar@{->} "b5";"a1"
\ar@{->} "b6";"b3"
\ar@{->} "b7";"b3"
\ar@{->} "b7";"b4"
\ar@{->} "b8";"b4"
\ar@{->} "b8";"b5"
\ar@{->} "b9";"b5"
\end{xy}
\end{equation}
Since $j_{k+1}=j_k-2$, we have the claim (a) for $l=0$.

Next, let us prove the claim (b) for $l=0$. We have seen that $(\mu_{k-\lfloor \frac{r}{2} \rfloor}\mu_{k+1}\mu_k(V_{\rm \bf{i}}))_{k+1}=(\mu_{k+1}\mu_k(V_{\rm \bf{i}}))_{k+1}$ is described as (\ref{mutlem1-pr1}). It follows from (\ref{mutlem1-pr2}) that $(\mu_{k-\lfloor \frac{r}{2} \rfloor}\mu_{k+1}\mu_k(V_{\rm \bf{i}}))_{k-\lfloor \frac{r}{2} \rfloor}$ has the submodule isomorphic to $(\mu_{k+1}\mu_k(V_{\rm \bf{i}}))_{k+1}$. Hence, we can find an injective homomorphism 
\begin{equation}\label{compmapfact}
(\mu_{k+1}\mu_k(V_{\rm \bf{i}}))_{k+1}\rightarrow (\mu_{k-\lfloor \frac{r}{2} \rfloor}\mu_{k+1}\mu_k(V_{\rm \bf{i}}))_{k-\lfloor \frac{r}{2} \rfloor}.
\end{equation}
 By the descriptions $(\ref{mutlem1-pr1})$ and $(\ref{mutlem1-pr1b})$, we see that the module $(\mu_{k+1}\mu_k(V_{\rm \bf{i}}))_{k+1}$ has the quotient isomorphic to $V_{r+k-\lfloor \frac{r}{2} \rfloor}$. Then, we have a surjective homomorphism $(\mu_{k+1}\mu_k(V_{\rm \bf{i}}))_{k+1}\rightarrow V_{r+k-\lfloor \frac{r}{2} \rfloor}$, which is, indeed, factorizable in the direct summands of $(\mu_{k-\lfloor \frac{r}{2} \rfloor}\mu_{k+1}\mu_k(V_{\rm \bf{i}}))$ since it can be written as a composite map as follows: We label
each basis of $V_{r+k-\lfloor \frac{r}{2} \rfloor}$ $(\ref{mutlem1-pr1b})$ as
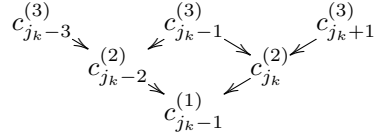
\[
\begin{xy}
(65,123)*{c^{(3)}_{j_k-3}}="4r-2k+1b",(85,123)*{c^{(3)}_{j_k-1}}="4r-2k+3b",
(105,123)*{c^{(3)}_{j_k+1}}="4r-2k+5",
(75,117)*{c^{(2)}_{j_k-2}}="4r-2k+2a",(95,117)*{c^{(2)}_{j_k}}="4r-2k+4",
(85,111)*{c^{(1)}_{j_k-1}}="4r-2k+3",
(35,111)*{}="12",
\ar@{->} "4r-2k+2a";"4r-2k+3"
\ar@{->} "4r-2k+4";"4r-2k+3"
\ar@{->} "4r-2k+1b";"4r-2k+2a"
\ar@{->} "4r-2k+3b";"4r-2k+2a"
\ar@{->} "4r-2k+3b";"4r-2k+4"
\ar@{->} "4r-2k+5";"4r-2k+4"
\end{xy}
\]
and each basis of $(\mu_{k-\lfloor \frac{r}{2} \rfloor}\mu_{k+1}\mu_k(V_{\rm \bf{i}}))_{k-\lfloor \frac{r}{2} \rfloor}$ (\ref{mutlem1-pr2}) as 
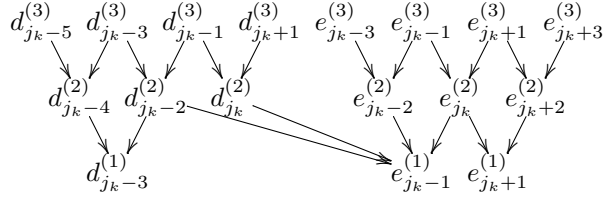
\[
\begin{xy}
(75,123)*{e^{(3)}_{j_k-3}}="a6",(85,123)*{e^{(3)}_{j_k-1}}="a7",(95,123)*{e^{(3)}_{j_k+1}}="a8",
(105,123)*{e^{(3)}_{j_k+3}}="a9",
(80,113)*{e^{(2)}_{j_k-2}}="a3",(90,113)*{e^{(2)}_{j_k}}="a4",(100,113)*{e^{(2)}_{j_k+2}}="a5",
(85,103)*{e^{(1)}_{j_k-1}}="a1",(95,103)*{e^{(1)}_{j_k+1}}="a2",
(35,123)*{d^{(3)}_{j_k-5}}="b6",(45,123)*{d^{(3)}_{j_k-3}}="b7",(55,123)*{d^{(3)}_{j_k-1}}="b8",
(65,123)*{d^{(3)}_{j_k+1}}="b9",
(40,113)*{d^{(2)}_{j_k-4}}="b3",(50,113)*{d^{(2)}_{j_k-2}}="b4",(60,113)*{d^{(2)}_{j_k}}="b5",
(45,103)*{d^{(1)}_{j_k-3}}="b1",
\ar@{->} "a3";"a1"
\ar@{->} "a4";"a1"
\ar@{->} "a4";"a2"
\ar@{->} "a5";"a2"
\ar@{->} "a6";"a3"
\ar@{->} "a7";"a3"
\ar@{->} "a7";"a4"
\ar@{->} "a8";"a4"
\ar@{->} "a8";"a5"
\ar@{->} "a9";"a5"
\ar@{->} "b3";"b1"
\ar@{->} "b4";"b1"
\ar@{->} "b4";"a1"
\ar@{->} "b5";"a1"
\ar@{->} "b6";"b3"
\ar@{->} "b7";"b3"
\ar@{->} "b7";"b4"
\ar@{->} "b8";"b4"
\ar@{->} "b8";"b5"
\ar@{->} "b9";"b5"
\end{xy}
\]
Then we can define the surjective homomorphism $(\mu_{k-\lfloor \frac{r}{2} \rfloor}\mu_{k+1}\mu_k(V_{\rm \bf{i}}))_{k-\lfloor \frac{r}{2} \rfloor}\rightarrow V_{r+k-\lfloor \frac{r}{2} \rfloor}$ as $e^{(1)}_{j_k-1}\mapsto c^{(1)}_{j_k-1}$, $d^{(2)}_{j},\ e^{(2)}_{j}\mapsto c^{(2)}_{j}$ ($j=j_k$, $j_k-2$),
$d^{(3)}_{j},\ e^{(3)}_{j}\mapsto c^{(3)}_{j}$ ($j=j_k+1$, $j_k-1$, $j_k-3$) and all others are mapped to $0$. Then the homomorphism $(\mu_{k+1}\mu_k(V_{\rm \bf{i}}))_{k+1}\rightarrow V_{r+k-\lfloor \frac{r}{2} \rfloor}$ coincides with the composite map 
\[(\mu_{k+1}\mu_k(V_{\rm \bf{i}}))_{k+1}\rightarrow (\mu_{k-\lfloor \frac{r}{2} \rfloor}\mu_{k+1}\mu_k(V_{\rm \bf{i}}))_{k-\lfloor \frac{r}{2} \rfloor}\rightarrow V_{r+k-\lfloor \frac{r}{2} \rfloor},\]
where the first map is the one in $(\ref{compmapfact})$.

The other non-zero homomorphisms from $(\mu_{k-\lfloor \frac{r}{2} \rfloor}\mu_{k+1}\mu_k(V_{\rm \bf{i}}))_{k+1}=(\mu_{k+1}\mu_k(V_{\rm \bf{i}}))_{k+1}$ to the direct summands of $(\mu_{k-\lfloor \frac{r}{2} \rfloor}\mu_{k+1}\mu_k(V_{\rm \bf{i}}))$ are factored through $(\mu_{k-\lfloor \frac{r}{2} \rfloor}\mu_{k+1}\mu_k(V_{\rm \bf{i}}))_{k-\lfloor \frac{r}{2} \rfloor}$. Thus, the homomorphism $(\ref{compmapfact})$ is not factorizable in the direct summands of $(\mu_{k-\lfloor \frac{r}{2} \rfloor}\mu_{k+1}\mu_k(V_{\rm \bf{i}}))$, and the exchange sequence associated to the direct summand $(\mu_{k-\lfloor \frac{r}{2} \rfloor}\mu_{k+1}\mu_k(V_{\rm \bf{i}}))_{k+1}=(\mu_{k+1}\mu_k(V_{\rm \bf{i}}))_{k+1}$ of $\mu_{k-\lfloor \frac{r}{2} \rfloor}\mu_{k+1}\mu_k(V_{\rm \bf{i}})$ is as follows:
\[ 0 \rightarrow (\mu_{k+1}\mu_k(V_{\rm \bf{i}}))_{k+1} \rightarrow (\mu_{k-\lfloor \frac{r}{2} \rfloor}\mu_{k+1}\mu_k(V_{\rm \bf{i}}))_{k-\lfloor \frac{r}{2} \rfloor}\rightarrow (\mu_{k+1}\mu_{k-\lfloor \frac{r}{2} \rfloor}\mu_{k+1}\mu_k(V_{\rm \bf{i}}))_{k+1}
\rightarrow 0. \]
By the above argument, we see that the module $(\mu_{k+1}\mu_{k-\lfloor \frac{r}{2} \rfloor}\mu_{k+1}\mu_k(V_{\rm \bf{i}}))_{k+1}$ is described as
\begin{equation}\label{mutlem1-pr3}
\begin{xy}
(75,123)*{_{j_k-3}}="a6",(85,123)*{_{j_k-1}}="a7",(95,123)*{_{j_k+1}}="a8",
(105,123)*{_{j_k+3}}="a9",
(80,116)*{_{j_k-2}}="a3",(90,116)*{_{j_k}}="a4",(100,116)*{_{j_k+2}}="a5",
(95,109)*{_{j_k+1}}="a2",
\ar@{->} "a4";"a2"
\ar@{->} "a5";"a2"
\ar@{->} "a6";"a3"
\ar@{->} "a7";"a3"
\ar@{->} "a7";"a4"
\ar@{->} "a8";"a4"
\ar@{->} "a8";"a5"
\ar@{->} "a9";"a5"
\end{xy}
\end{equation}
which means the claim (b) for $l=0$.

Next, we assume that the claims (a) and (b) are proven for $0,1,\cdots,l$. Let us consider the claim(a) for $l+1$, and then construct the exchange sequence associated to the direct summand $(\mu_{k+l+2}\mu_{k+l+1}\mu[l](V_{\rm \bf{i}}))_{k-\lfloor \frac{r}{2} \rfloor+l+1}$ of $\mu_{k+l+2}\mu_{k+l+1}\mu[l](V_{\rm \bf{i}})$ as in (\ref{exseq-af1}) below. Note that since 
the mutation $\mu_{k-\lfloor \frac{r}{2} \rfloor+l+1}$ does not appear in
$\mu_{k+l+2}\mu_{k+l+1}\mu[l]$, we have $(\mu_{k+l+2}\mu_{k+l+1}\mu[l](V_{\rm \bf{i}}))_{k-\lfloor \frac{r}{2} \rfloor+l+1}=(V_{\rm \bf{i}})_{k-\lfloor \frac{r}{2} \rfloor+l+1}=S_{j_k-2l-3}$ (see (\ref{iniex1-1})). By the induction hypothesis, the module $(\mu_{k+l+2}\mu_{k+l+1}\mu[l](V_{\rm \bf{i}}))_{k-\lfloor \frac{r}{2} \rfloor+l}=(\mu[l](V_{\rm \bf{i}}))_{k-\lfloor \frac{r}{2} \rfloor+l}$ is described as (\ref{mutlemclaim1}), and it has the simple submodule isomorphic to $S_{j_k-2l-3}$. It follows from Theorem \ref{GLSthm} and a similar argument to the proof of Lemma \ref{mut0} that the module $(\mu_{k+l+2}\mu_{k+l+1}\mu[l](V_{\rm \bf{i}}))_{k+l+2}$ is the same as $(\mu_{k+l+2}(V_{\rm \bf{i}}))_{k+l+2}$, and is described as follows:
\begin{equation}\label{mutlem1-pr4}
\begin{xy}
(35,123)*{_{j_k-2l-7}}="4r-2k-1",(55,123)*{_{j_k-2l-5}}="4r-2k+1b",(75,123)*{_{j_k-2l-3}}="4r-2k+3b",
(95,123)*{_{j_k-2l-1}}="4r-2k+5",
(45,115)*{_{j_k-2l-6}}="4r-2k",(65,115)*{_{j_k-2l-4}}="4r-2k+2a",(85,115)*{_{j_k-2l-2}}="4r-2k+4",
(55,107)*{_{j_k-2l-5}}="4r-2k+1",(75,107)*{_{j_k-2l-3}}="4r-2k+3",
(35,107)*{}="12",
\ar@{->} "4r-2k";"4r-2k+1"
\ar@{->} "4r-2k+2a";"4r-2k+1"
\ar@{->} "4r-2k+2a";"4r-2k+3"
\ar@{->} "4r-2k+4";"4r-2k+3"
\ar@{->} "4r-2k-1";"4r-2k"
\ar@{->} "4r-2k+1b";"4r-2k"
\ar@{->} "4r-2k+1b";"4r-2k+2a"
\ar@{->} "4r-2k+3b";"4r-2k+2a"
\ar@{->} "4r-2k+3b";"4r-2k+4"
\ar@{->} "4r-2k+5";"4r-2k+4"
\end{xy}
\end{equation}
Hence, the module $(\mu_{k+l+2}(V_{\rm \bf{i}}))_{k+l+2}$ has the simple submodule isomorphic to $S_{j_k-2l-3}$. 
It follows from (\ref{iniex1-3}) and $j_{r+k-\lfloor \frac{r}{2} \rfloor+l+1}=j_{k-\lfloor \frac{r}{2} \rfloor+l+1}=j_k-2l-3$ that the module $(\mu_{k+l+2}\mu_{k+l+1}\mu[l](V_{\rm \bf{i}}))_{r+k-\lfloor \frac{r}{2} \rfloor+l+1}=(V_{\rm \bf{i}})_{r+k-\lfloor \frac{r}{2} \rfloor+l+1}$ is described as
\[
\begin{xy}
(65,123)*{j_k-2l-5}="4r-2k+1b",(85,123)*{j_k-2l-3}="4r-2k+3b",
(105,123)*{j_k-2l-1}="4r-2k+5",
(75,116)*{j_k-2l-4}="4r-2k+2a",(95,116)*{j_k-2l-2}="4r-2k+4",
(85,109)*{j_k-2l-3}="4r-2k+3",
(35,109)*{}="12",
\ar@{->} "4r-2k+2a";"4r-2k+3"
\ar@{->} "4r-2k+4";"4r-2k+3"
\ar@{->} "4r-2k+1b";"4r-2k+2a"
\ar@{->} "4r-2k+3b";"4r-2k+2a"
\ar@{->} "4r-2k+3b";"4r-2k+4"
\ar@{->} "4r-2k+5";"4r-2k+4"
\end{xy}
\]
and there exists an injective homomorphism $S_{j_k-2l-3}\rightarrow (V_{\rm \bf{i}})_{r+k-\lfloor \frac{r}{2} \rfloor+l+1}$. But, this map is factorizable since it can be written as composite map $S_{j_k-2l-3}\rightarrow (\mu_{k+l+2}(V_{\rm \bf{i}}))_{k+l+2} \rightarrow (V_{\rm \bf{i}})_{r+k-\lfloor \frac{r}{2} \rfloor+l+1}$. By the induction hypothesis, the other direct summands of $\mu_{k+l+2}\mu_{k+l+1}\mu[l](V_{\rm \bf{i}})$ do not have the simple submodule isomorphic to $S_{j_k-2l-3}$. Thus, the exchange sequence associated to the direct summand $(\mu_{k+l+2}\mu_{k+l+1}\mu[l](V_{\rm \bf{i}}))_{k-\lfloor \frac{r}{2} \rfloor+l+1}=S_{j_k-2l-3}$ of $\mu_{k+l+2}\mu_{k+l+1}\mu[l](V_{\rm \bf{i}})$
is
\begin{equation}\label{exseq-af1}
 0\rightarrow S_{j_k-2l-3}\rightarrow (\mu_{k+l+2}(V_{\rm \bf{i}}))_{k+l+2}\oplus (\mu[l](V_{\rm \bf{i}}))_{k-\lfloor \frac{r}{2} \rfloor+l}\rightarrow \end{equation}
\[ (\mu_{k-\lfloor \frac{r}{2} \rfloor+l+1}\mu_{k+l+2}\mu_{k+l+1}\mu[l](V_{\rm \bf{i}}))_{k-\lfloor \frac{r}{2} \rfloor+l+1}\rightarrow 0.
\]

From $j_{k+l}=j_k-2l$, the module $(\mu[l+1](V_{\rm \bf{i}}))_{k-\lfloor \frac{r}{2} \rfloor+l+1}=(\mu_{k-\lfloor \frac{r}{2} \rfloor+l+1}\mu_{k+l+2}\mu_{k+l+1}\mu[l](V_{\rm \bf{i}}))_{k-\lfloor \frac{r}{2} \rfloor+l+1}$ is described as

\begin{equation}\label{mutlem1-pr5}
\begin{xy}
(63,123)*{_{j_k-3}}="a6",(71,123)*{_{j_k-1}}="a7",(79,123)*{_{j_k+1}}="a8",
(87,123)*{_{j_k+3}}="a9",
(67,113)*{_{j_k-2}}="a3",(75,113)*{_{j_k}}="a4",(83,113)*{_{j_k+2}}="a5",
(71,103)*{_{j_k-1}}="a1",(79,103)*{ _{j_k+1}}="a2",
(32,123)*{_{j_k-5}}="b6",(40,123)*{_{j_k-3}}="b7",(48,123)*{_{j_k-1}}="b8",
(56,123)*{_{j_k+1}}="b9",
(36,113)*{_{j_k-4}}="b3",(44,113)*{_{j_k-2}}="b4",(52,113)*{_{j_k}}="b5",
(40,103)*{_{j_k-3}}="b1",
(-17,123)*{_{j_{k+l}-7}}="c6",(-7,123)*{_{j_{k+l}-5}}="c7",(3,123)*{_{j_{k+l}-3}}="c8",
(13,123)*{_{j_{k+l}-1}}="c9",
(-12,113)*{_{j_{k+l}-6}}="c3",(-2,113)*{_{j_{k+l}-4}}="c4",(8,113)*{_{j_{k+l}-2}}="c5",
(-7,103)*{_{j_{k+l}-5}}="c1",(22,103)*{_{j_{k+l}-3}}="z1",(20,113)*{\cdots}="dot",
(26,113)*{\ }="hs",
\ar@{->} "hs";"b1"
\ar@{->} "dot";"b1"
\ar@{->} "a3";"a1"
\ar@{->} "a4";"a1"
\ar@{->} "a4";"a2"
\ar@{->} "a5";"a2"
\ar@{->} "a6";"a3"
\ar@{->} "a7";"a3"
\ar@{->} "a7";"a4"
\ar@{->} "a8";"a4"
\ar@{->} "a8";"a5"
\ar@{->} "a9";"a5"
\ar@{->} "b3";"b1"
\ar@{->} "b4";"b1"
\ar@{->} "b4";"a1"
\ar@{->} "b5";"a1"
\ar@{->} "b6";"b3"
\ar@{->} "b7";"b3"
\ar@{->} "b7";"b4"
\ar@{->} "b8";"b4"
\ar@{->} "b8";"b5"
\ar@{->} "b9";"b5"
\ar@{->} "c3";"c1"
\ar@{->} "c4";"c1"
\ar@{->} "c4";"z1"
\ar@{->} "c5";"z1"
\ar@{->} "c6";"c3"
\ar@{->} "c7";"c3"
\ar@{->} "c7";"c4"
\ar@{->} "c8";"c4"
\ar@{->} "c8";"c5"
\ar@{->} "c9";"c5"
\end{xy}
\end{equation}
Since $j_{k+l+1}=j_{k+l}-2$, we get (a) for $l+1$.

Next, we consider the claim (b) for $l+1$. The module $(\mu[l+1](V_{\rm \bf{i}}))_{k+l+2}=(\mu_{k+l+2}\mu_{k+l+1}\mu[l](V_{\rm \bf{i}}))_{k+l+2}$ is described as (\ref{mutlem1-pr4}). By the description (\ref{mutlem1-pr5}) of the module $(\mu[l+1](V_{\rm \bf{i}}))_{k-\lfloor \frac{r}{2} \rfloor+l+1}$, it has the submodule isomorphic to $(\mu_{k+l+2}\mu_{k+l+1}\mu[l](V_{\rm \bf{i}}))_{k+l+2}$. Using the same argument in the proof of claim (b) for $l=0$, the other non-zero homomorphisms from $(\mu_{k+l+2}\mu_{k+l+1}\mu[l](V_{\rm \bf{i}}))_{k+l+2}$ to the direct summands of $(\mu[l+1](V_{\rm \bf{i}}))$ are factored through $(\mu[l+1](V_{\rm \bf{i}}))_{k-\lfloor \frac{r}{2} \rfloor+l+1}$. Thus, the exchange sequence associated to the direct summand $(\mu[l+1](V_{\rm \bf{i}}))_{k+l+2}$ of $(\mu[l+1](V_{\rm \bf{i}}))$ is
\begin{equation}\label{exseq-2}
0\rightarrow (\mu[l+1](V_{\rm \bf{i}}))_{k+l+2}\rightarrow (\mu[l+1](V_{\rm \bf{i}}))_{k-\lfloor \frac{r}{2} \rfloor+l+1}
\rightarrow (\mu_{k+l+2}\mu[l+1](V_{\rm \bf{i}}))_{k+l+2} \rightarrow 0,
\end{equation}
which yields the following description of the module $(\mu_{k+l+2}\mu[l+1](V_{\rm \bf{i}}))_{k+l+2}$:
\begin{equation}\label{mutlem1-pr6}
\begin{xy}
(63,123)*{_{j_k-3}}="a6",(71,123)*{_{j_k-1}}="a7",(79,123)*{_{j_k+1}}="a8",
(87,123)*{_{j_k+3}}="a9",
(67,113)*{_{j_k-2}}="a3",(75,113)*{_{j_k}}="a4",(83,113)*{_{j_k+2}}="a5",
(71,103)*{_{j_k-1}}="a1",(79,103)*{ _{j_k+1}}="a2",
(32,123)*{_{j_k-5}}="b6",(40,123)*{_{j_k-3}}="b7",(48,123)*{_{j_k-1}}="b8",
(56,123)*{_{j_k+1}}="b9",
(36,113)*{_{j_k-4}}="b3",(44,113)*{_{j_k-2}}="b4",(52,113)*{_{j_k}}="b5",
(40,103)*{_{j_k-3}}="b1",
(-17,123)*{_{j_{k+l}-5}}="c6",(-7,123)*{_{j_{k+l}-3}}="c7",(3,123)*{_{j_{k+l}-1}}="c8",
(13,123)*{_{j_{k+l}+1}}="c9",
(-12,113)*{_{j_{k+l}-4}}="c3",(-2,113)*{_{j_{k+l}-2}}="c4",(8,113)*{_{j_{k+l}}}="c5",
(22,103)*{_{j_{k+l}-1}}="z1",(20,113)*{\cdots}="dot",
(26,113)*{\ }="hs",
\ar@{->} "hs";"b1"
\ar@{->} "dot";"b1"
\ar@{->} "a3";"a1"
\ar@{->} "a4";"a1"
\ar@{->} "a4";"a2"
\ar@{->} "a5";"a2"
\ar@{->} "a6";"a3"
\ar@{->} "a7";"a3"
\ar@{->} "a7";"a4"
\ar@{->} "a8";"a4"
\ar@{->} "a8";"a5"
\ar@{->} "a9";"a5"
\ar@{->} "b3";"b1"
\ar@{->} "b4";"b1"
\ar@{->} "b4";"a1"
\ar@{->} "b5";"a1"
\ar@{->} "b6";"b3"
\ar@{->} "b7";"b3"
\ar@{->} "b7";"b4"
\ar@{->} "b8";"b4"
\ar@{->} "b8";"b5"
\ar@{->} "b9";"b5"
\ar@{->} "c4";"z1"
\ar@{->} "c5";"z1"
\ar@{->} "c6";"c3"
\ar@{->} "c7";"c3"
\ar@{->} "c7";"c4"
\ar@{->} "c8";"c4"
\ar@{->} "c8";"c5"
\ar@{->} "c9";"c5"
\end{xy}
\end{equation}
Because of $j_{k+l+1}=j_{k+l}-2$, we get (b) for $l+1$. \qed

We can similarly verify the following two lemmas.

\begin{lem}\label{mutlem2}
We use the notation as in Proposition \ref{spprop2} and let $j_k$ be the $k$-th index of {\rm \bf{i}} $(\ref{redwords2})$ from the right.
\begin{enumerate}
\item[$(a)$] The module $(\mu'[l](V_{\rm \bf{i}}))_{k+l+1}$ is described as follows:
\begin{equation}\label{mutlem2claim1}
\begin{xy}
(63,123)*{_{j_k-4}}="a6",(71,123)*{_{j_k-2}}="a7",(79,123)*{_{j_k}}="a8",
(87,123)*{_{j_k+2}}="a9",
(67,113)*{_{j_k-3}}="a3",(75,113)*{_{j_k-1}}="a4",(83,113)*{_{j_k+1}}="a5",
(71,103)*{_{j_k-2}}="a1",
(32,123)*{_{j_k-6}}="b6",(40,123)*{_{j_k-4}}="b7",(48,123)*{_{j_k-2}}="b8",
(56,123)*{_{j_k}}="b9",
(36,113)*{_{j_k-5}}="b3",(44,113)*{_{j_k-3}}="b4",(52,113)*{_{j_k-1}}="b5",
(40,103)*{_{j_k-4}}="b1",
(-17,123)*{_{j_{k+l}-6}}="c6",(-7,123)*{_{j_{k+l}-4}}="c7",(3,123)*{_{j_{k+l}-2}}="c8",
(13,123)*{_{j_{k+l}}}="c9",
(-12,113)*{_{j_{k+l}-5}}="c3",(-2,113)*{_{j_{k+l}-3}}="c4",(8,113)*{_{j_{k+l}-1}}="c5",
(-7,103)*{_{j_{k+l}-4}}="c1",(22,103)*{_{j_{k+l}-2}}="z1",(20,113)*{\cdots}="dot",
(26,113)*{\ }="hs",
\ar@{->} "hs";"b1"
\ar@{->} "dot";"b1"
\ar@{->} "a3";"a1"
\ar@{->} "a4";"a1"
\ar@{->} "a6";"a3"
\ar@{->} "a7";"a3"
\ar@{->} "a7";"a4"
\ar@{->} "a8";"a4"
\ar@{->} "a8";"a5"
\ar@{->} "a9";"a5"
\ar@{->} "b3";"b1"
\ar@{->} "b4";"b1"
\ar@{->} "b4";"a1"
\ar@{->} "b5";"a1"
\ar@{->} "b6";"b3"
\ar@{->} "b7";"b3"
\ar@{->} "b7";"b4"
\ar@{->} "b8";"b4"
\ar@{->} "b8";"b5"
\ar@{->} "b9";"b5"
\ar@{->} "c3";"c1"
\ar@{->} "c4";"c1"
\ar@{->} "c4";"z1"
\ar@{->} "c5";"z1"
\ar@{->} "c6";"c3"
\ar@{->} "c7";"c3"
\ar@{->} "c7";"c4"
\ar@{->} "c8";"c4"
\ar@{->} "c8";"c5"
\ar@{->} "c9";"c5"
\end{xy}
\end{equation}
Note that $j_{k+l}=j_k-2l$ from $(\ref{redwords2})$.

\item[$(b)$] The module $(\mu_{\lfloor \frac{r}{2} \rfloor+k+l+2}\mu'[l](V_{\rm \bf{i}}))_{\lfloor \frac{r}{2} \rfloor+k+l+2}$ is described as follows:
\begin{equation}\label{mutlem2claim2}
\begin{xy}
(63,123)*{_{j_k-4}}="a6",(71,123)*{_{j_k-2}}="a7",(79,123)*{_{j_k}}="a8",
(87,123)*{_{j_k+2}}="a9",
(67,113)*{_{j_k-3}}="a3",(75,113)*{_{j_k-1}}="a4",(83,113)*{_{j_k+1}}="a5",
(71,103)*{_{j_k-2}}="a1",
(32,123)*{_{j_k-6}}="b6",(40,123)*{_{j_k-4}}="b7",(48,123)*{_{j_k-2}}="b8",
(56,123)*{_{j_k}}="b9",
(36,113)*{_{j_k-5}}="b3",(44,113)*{_{j_k-3}}="b4",(52,113)*{_{j_k-1}}="b5",
(40,103)*{_{j_k-4}}="b1",
(-17,123)*{_{j_{k+l}-4}}="c6",(-7,123)*{_{j_{k+l}-2}}="c7",(3,123)*{_{j_{k+l}}}="c8",
(13,123)*{_{j_{k+l}+2}}="c9",
(-12,113)*{_{j_{k+l}-3}}="c3",(-2,113)*{_{j_{k+l}-1}}="c4",(8,113)*{_{j_{k+l}+1}}="c5",
(22,103)*{_{j_{k+l}}}="z1",(20,113)*{\cdots}="dot",
(26,113)*{\ }="hs",
\ar@{->} "hs";"b1"
\ar@{->} "dot";"b1"
\ar@{->} "a3";"a1"
\ar@{->} "a4";"a1"
\ar@{->} "a6";"a3"
\ar@{->} "a7";"a3"
\ar@{->} "a7";"a4"
\ar@{->} "a8";"a4"
\ar@{->} "a8";"a5"
\ar@{->} "a9";"a5"
\ar@{->} "b3";"b1"
\ar@{->} "b4";"b1"
\ar@{->} "b4";"a1"
\ar@{->} "b5";"a1"
\ar@{->} "b6";"b3"
\ar@{->} "b7";"b3"
\ar@{->} "b7";"b4"
\ar@{->} "b8";"b4"
\ar@{->} "b8";"b5"
\ar@{->} "b9";"b5"
\ar@{->} "c4";"z1"
\ar@{->} "c5";"z1"
\ar@{->} "c6";"c3"
\ar@{->} "c7";"c3"
\ar@{->} "c7";"c4"
\ar@{->} "c8";"c4"
\ar@{->} "c8";"c5"
\ar@{->} "c9";"c5"
\end{xy}
\end{equation}
\end{enumerate}
\end{lem}

\begin{lem}\label{mutlem3}
We use the notation as in Proposition \ref{spprop3}.
\begin{enumerate}
\item[$(a)$] The module $(\mu''[l](V_{\rm \bf{i}}))_{\lfloor \frac{r+1}{2} \rfloor-l-1}$ is described as follows:
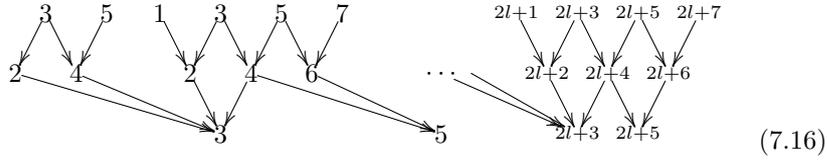
\begin{equation}\label{mutlem3claim1}
\begin{xy}
(77,123)*{_{2l+1}}="x6",(85,123)*{_{2l+3}}="x7",(93,123)*{_{2l+5}}="x8",
(101,123)*{_{2l+7}}="x9",
(81,115)*{_{2l+2}}="x3",(89,115)*{_{2l+4}}="x4",(97,115)*{_{2l+6}}="x5",
(85,107)*{_{2l+3}}="x1",(93,107)*{_{2l+5}}="x2",
(30,123)*{1}="a6",(38,123)*{3}="a7",(46,123)*{5}="a8",
(54,123)*{7}="a9",
(34,115)*{2}="a3",(42,115)*{4}="a4",(50,115)*{6}="a5",
(38,107)*{3}="a1",(67,107)*{5}="a2",(67,115)*{\cdots}="dot",
(71,115)*{\ }="hs",
(15,123)*{3}="b8",
(23,123)*{5}="b9",
(11,115)*{2}="b4",(19,115)*{4}="b5",
\ar@{->} "hs";"x1"
\ar@{->} "dot";"x1"
\ar@{->} "a3";"a1"
\ar@{->} "a4";"a1"
\ar@{->} "a4";"a2"
\ar@{->} "a5";"a2"
\ar@{->} "a6";"a3"
\ar@{->} "a7";"a3"
\ar@{->} "a7";"a4"
\ar@{->} "a8";"a4"
\ar@{->} "a8";"a5"
\ar@{->} "a9";"a5"
\ar@{->} "b4";"a1"
\ar@{->} "b5";"a1"
\ar@{->} "b8";"b4"
\ar@{->} "b8";"b5"
\ar@{->} "b9";"b5"
\ar@{->} "x3";"x1"
\ar@{->} "x4";"x1"
\ar@{->} "x4";"x2"
\ar@{->} "x5";"x2"
\ar@{->} "x6";"x3"
\ar@{->} "x7";"x3"
\ar@{->} "x7";"x4"
\ar@{->} "x8";"x4"
\ar@{->} "x8";"x5"
\ar@{->} "x9";"x5"
\end{xy}
\end{equation}

\item[$(b)$] The module $(\mu_{r-l-1}\mu''[l](V_{\rm \bf{i}}))_{r-l-1}$ is described as follows:
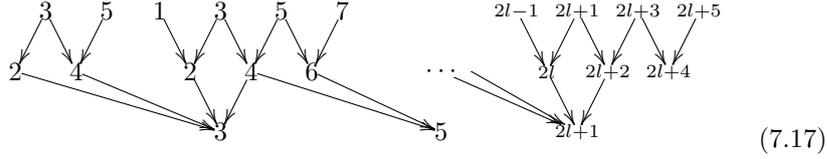
\begin{equation}\label{mutlem3claim2}
\begin{xy}
(77,123)*{_{2l-1}}="x6",(85,123)*{_{2l+1}}="x7",(93,123)*{_{2l+3}}="x8",
(101,123)*{_{2l+5}}="x9",
(81,115)*{_{2l}}="x3",(89,115)*{_{2l+2}}="x4",(97,115)*{_{2l+4}}="x5",
(85,107)*{_{2l+1}}="x1",
(30,123)*{1}="a6",(38,123)*{3}="a7",(46,123)*{5}="a8",
(54,123)*{7}="a9",
(34,115)*{2}="a3",(42,115)*{4}="a4",(50,115)*{6}="a5",
(38,107)*{3}="a1",(67,107)*{5}="a2",(67,115)*{\cdots}="dot",
(71,115)*{\ }="hs",
(15,123)*{3}="b8",
(23,123)*{5}="b9",
(11,115)*{2}="b4",(19,115)*{4}="b5",
\ar@{->} "hs";"x1"
\ar@{->} "dot";"x1"
\ar@{->} "a3";"a1"
\ar@{->} "a4";"a1"
\ar@{->} "a4";"a2"
\ar@{->} "a5";"a2"
\ar@{->} "a6";"a3"
\ar@{->} "a7";"a3"
\ar@{->} "a7";"a4"
\ar@{->} "a8";"a4"
\ar@{->} "a8";"a5"
\ar@{->} "a9";"a5"
\ar@{->} "b4";"a1"
\ar@{->} "b5";"a1"
\ar@{->} "b8";"b4"
\ar@{->} "b8";"b5"
\ar@{->} "b9";"b5"
\ar@{->} "x3";"x1"
\ar@{->} "x4";"x1"
\ar@{->} "x6";"x3"
\ar@{->} "x7";"x3"
\ar@{->} "x7";"x4"
\ar@{->} "x8";"x4"
\ar@{->} "x8";"x5"
\ar@{->} "x9";"x5"
\end{xy}
\end{equation}
\end{enumerate}
Furthermore, if $r$ is odd, then the module $(\mu_1\mu_{\frac{r+3}{2}}\mu[\frac{r-1}{2}-2](V_{\rm \bf{i}}))_1$ is described as
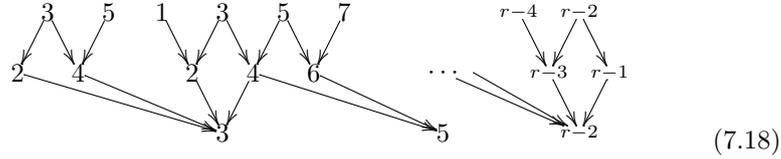
\begin{equation}\label{mutlem3claim1-r}
\begin{xy}
(77,123)*{_{r-4}}="x6",(85,123)*{_{r-2}}="x7",
(81,115)*{_{r-3}}="x3",(89,115)*{_{r-1}}="x4",
(85,107)*{_{r-2}}="x1",
(30,123)*{1}="a6",(38,123)*{3}="a7",(46,123)*{5}="a8",
(54,123)*{7}="a9",
(34,115)*{2}="a3",(42,115)*{4}="a4",(50,115)*{6}="a5",
(38,107)*{3}="a1",(67,107)*{5}="a2",(67,115)*{\cdots}="dot",
(71,115)*{\ }="hs",
(15,123)*{3}="b8",
(23,123)*{5}="b9",
(11,115)*{2}="b4",(19,115)*{4}="b5",
\ar@{->} "hs";"x1"
\ar@{->} "dot";"x1"
\ar@{->} "a3";"a1"
\ar@{->} "a4";"a1"
\ar@{->} "a4";"a2"
\ar@{->} "a5";"a2"
\ar@{->} "a6";"a3"
\ar@{->} "a7";"a3"
\ar@{->} "a7";"a4"
\ar@{->} "a8";"a4"
\ar@{->} "a8";"a5"
\ar@{->} "a9";"a5"
\ar@{->} "b4";"a1"
\ar@{->} "b5";"a1"
\ar@{->} "b8";"b4"
\ar@{->} "b8";"b5"
\ar@{->} "b9";"b5"
\ar@{->} "x3";"x1"
\ar@{->} "x4";"x1"
\ar@{->} "x6";"x3"
\ar@{->} "x7";"x3"
\ar@{->} "x7";"x4"
\end{xy}
\end{equation}
\end{lem}

{\sl[Proof of Proposition \ref{spprop1}.]}

For any cluster $\mathbb{T}$ and $s\in[r+1,2r]\cup [-r,-1]$, we set $(\varphi_{\mathbb{T}})_s=(\varphi_{\mathbb{V}})_s$. 
Using the induction on $l$, let us prove Proposition \ref{spprop1} (a) and (b) simultaneously. For $l=0$, let us calculate $(\varphi^G_{(\mu[0]\mathbb{V})})_{k-\lfloor \frac{r}{2} \rfloor }$. In Sect.\ref{gmc}, we see that the vertices and arrows around the vertex $(\varphi_{\mathbb{V}})_k$ in the quiver $\Gamma_{{\rm \bf{i}}}$ are described as follows:
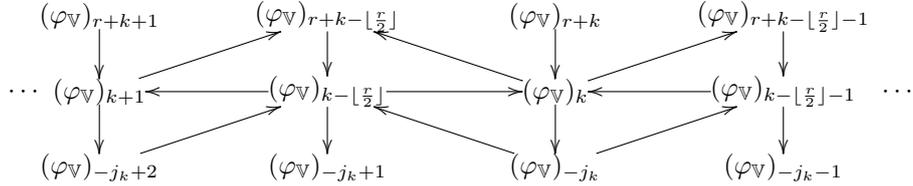
\[
\begin{xy}
(115,90)*{\cdots}="emp1",
(100,100)*{(\varphi_{\mathbb{V}})_{r+k-\lfloor \frac{r}{2} \rfloor -1}}="7",
(100,90) *{(\varphi_{\mathbb{V}})_{k-\lfloor \frac{r}{2} \rfloor -1}}="3",
(100,80)*{(\varphi_{\mathbb{V}})_{-j_k-1}}="-4",
(70,100)*{(\varphi_{\mathbb{V}})_{r+k}}="r+k",
(70,90) *{(\varphi_{\mathbb{V}})_k}="k",
(70,80)*{(\varphi_{\mathbb{V}})_{-j_k}}="-j_k",
(40,100)*{(\varphi_{\mathbb{V}})_{r+k-\lfloor \frac{r}{2} \rfloor}}="8",
(40,90) *{(\varphi_{\mathbb{V}})_{k-\lfloor \frac{r}{2} \rfloor }}="4",
(40,80)*{(\varphi_{\mathbb{V}})_{-j_k+1}}="-2",
(10,100)*{(\varphi_{\mathbb{V}})_{r+k+1}}="6",
(10,90) *{(\varphi_{\mathbb{V}})_{k+1}}="2",
(10,80)*{(\varphi_{\mathbb{V}})_{-j_k+2}}="-1",
(0,90)*{\cdots}="emp",
\ar@{->} "7";"3"
\ar@{->} "r+k";"k"
\ar@{->} "8";"4"
\ar@{->} "6";"2"
\ar@{->} "3";"-4"
\ar@{->} "k";"-j_k"
\ar@{->} "4";"-2"
\ar@{->} "2";"-1"
\ar@{->} "k";"7"
\ar@{->} "k";"8"
\ar@{->} "2";"8"
\ar@{->} "3";"k"
\ar@{->} "4";"k"
\ar@{->} "4";"2"
\ar@{->} "-j_k";"3"
\ar@{->} "-j_k";"4"
\ar@{->} "-1";"4"
\end{xy}
\]
Applying the mutation $\mu_{k+1}\mu_k$ to this quiver, the arrows between $(\varphi_{\mathbb{V}})_{k-\lfloor \frac{r}{2} \rfloor }$ and $(\varphi_{\mathbb{V}})_{-s}$ $(1\leq s\leq r)$
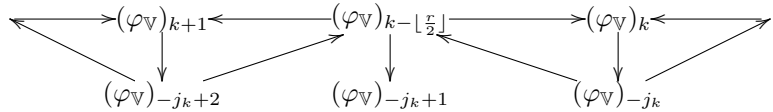
\[
\begin{xy}
(90,90) *{\ \ }="3",
(90,80)*{}="-4",
(70,90) *{(\varphi_{\mathbb{V}})_k}="k",
(70,80)*{(\varphi_{\mathbb{V}})_{-j_k}}="-j_k",
(40,90) *{(\varphi_{\mathbb{V}})_{k-\lfloor \frac{r}{2} \rfloor }}="4",
(40,80)*{(\varphi_{\mathbb{V}})_{-j_k+1}}="-2",
(10,90) *{(\varphi_{\mathbb{V}})_{k+1}}="2",
(10,80)*{(\varphi_{\mathbb{V}})_{-j_k+2}}="-1",
(-10,90)*{\ }="emp",
\ar@{->} "k";"-j_k"
\ar@{->} "4";"-2"
\ar@{->} "2";"-1"
\ar@{->} "3";"k"
\ar@{->} "4";"k"
\ar@{->} "4";"2"
\ar@{->} "-j_k";"3"
\ar@{->} "-j_k";"4"
\ar@{->} "-1";"4"
\ar@{->} "-1";"emp"
\ar@{->} "emp";"2"
\end{xy}
\]
are transformed to
\[
\begin{xy}
(90,90) *{\ \ }="3",
(90,80)*{}="-4",
(70,90) *{(\varphi_{\mu_k(\mathbb{V})})_k}="k",
(70,80)*{(\varphi_{\mathbb{V}})_{-j_k}}="-j_k",
(40,90) *{(\varphi_{\mathbb{V}})_{k-\lfloor \frac{r}{2} \rfloor }}="4",
(40,80)*{(\varphi_{\mathbb{V}})_{-j_k+1}}="-2",
(10,90) *{(\varphi_{\mu_{k+1}\mu_k(\mathbb{V})})_{k+1}}="2",
(10,80)*{(\varphi_{\mathbb{V}})_{-j_k+2}}="-1",
(-10,90)*{\ }="emp",
\ar@{->} "-j_k";"k"
\ar@{->} "4";"-2"
\ar@{->} "-1";"2"
\ar@{->} "k";"3"
\ar@{->} "k";"4"
\ar@{->} "2";"4"
\ar@{->} "2";"emp"
\end{xy}
\]
by Lemma \ref{mutgamlem}. Similarly, the arrows between $(\varphi_{(\mu_{k+1}\mu_k\mathbb{V})})_{k-\lfloor \frac{r}{2} \rfloor }=(\varphi_{\mathbb{V}})_{k-\lfloor \frac{r}{2} \rfloor }$ and $(\varphi_{\mu_{k+1}\mu_k\mathbb{V}})_{s}$ $(1\leq s\leq 2r)$ in $\mu_{k+1}\mu_k(\Gamma_{\textbf{i}})$ are
\[
\begin{xy}
(90,100)*{(\varphi_{\mathbb{V}})_{r+k-\lfloor \frac{r}{2} \rfloor -1}}="7",
(70,90) *{(\varphi_{(\mu_k\mathbb{V})})_k}="k",
(40,100)*{(\varphi_{\mathbb{V}})_{r+k-\lfloor \frac{r}{2} \rfloor}}="8",
(40,90) *{(\varphi_{\mathbb{V}})_{k-\lfloor \frac{r}{2} \rfloor }}="4",
(10,90) *{(\varphi_{(\mu_{k+1}\mu_{k}\mathbb{V})})_{k+1}}="2",
(-10,100)*{(\varphi_{\mathbb{V}})_{r+k-\lfloor \frac{r}{2} \rfloor+1}}="9",
\ar@{->} "k";"4"
\ar@{->} "4";"7"
\ar@{->} "4";"8"
\ar@{->} "2";"4"
\ar@{->} "4";"9"
\end{xy}
\]
Thus, by the exchange relation $(\ref{exrel})$,
\begin{eqnarray*}
& &\hspace{-20pt} (\varphi^G_{(\mu[0]\mathbb{V})})_{k-\lfloor \frac{r}{2} \rfloor } \\
&=&
\frac{(\varphi^G_{(\mu_k\mathbb{V})})_k (\varphi^G_{(\mu_{k+1}\mu_{k}\mathbb{V})})_{k+1}
+(\varphi^G_{\mathbb{V}})_{-j_k+1}
(\varphi^G_{\mathbb{V}})_{r+k-\lfloor \frac{r}{2} \rfloor-1} (\varphi^G_{\mathbb{V}})_{r+k-\lfloor \frac{r}{2} \rfloor}(\varphi^G_{\mathbb{V}})_{r+k-\lfloor \frac{r}{2} \rfloor+1}}{(\varphi^G_{\mathbb{V}})_{k-\lfloor \frac{r}{2} \rfloor }}.
\end{eqnarray*}

By (\ref{mukcoeff}) in Example \ref{initialex4}, we can write
\begin{equation}\label{muk}
 (\varphi^G_{(\mu_k\mathbb{V})})_k =
(\varphi_{(\mu_k V)_k})\circ x^G_{{\rm \bf{i}}}(\Phi_H;\phi(\textbf{Y})) =
(\Phi_H(\textbf{Y}))^{\Lambda_{j_k-1}}(\Phi_H(\textbf{Y}))^{\Lambda_{j_k+1}}\cdot
(\varphi_{(\mu_k V)_k})\circ x^G_{{\rm \bf{i}}}(1;\phi(\textbf{Y})). 
\end{equation}
 Similarly, we have 
\begin{equation}\label{muk+1}
(\varphi^G_{(\mu_{k+1}\mu_k\mathbb{V})})_{k+1}=(\varphi^G_{(\mu_{k+1}\mathbb{V})})_{k+1}
=(\Phi_H(\textbf{Y}))^{\Lambda_{j_k-3}}(\Phi_H(\textbf{Y}))^{\Lambda_{j_k-1}}\cdot
(\varphi_{(\mu_{k+1} V)_{k+1}})\circ x^G_{{\rm \bf{i}}}(1;\phi(\textbf{Y})).
\end{equation}
Therefore, using (\ref{genbasic}), (\ref{gtol}), (\ref{muk}), (\ref{muk+1}) and Theorem \ref{GLSthm} (ii), we obtain
\begin{eqnarray}
& &\hspace{-20pt} (\varphi^G_{(\mu[0]\mathbb{V})})_{k-\lfloor \frac{r}{2} \rfloor }(a;\textbf{Y})=\frac{(\Phi_H(\textbf{Y}))^{2\Lambda_{j_k-1}+\Lambda_{j_k-3}+\Lambda_{j_k+1}}}{(\Phi_H(\textbf{Y}))^{\Lambda_{j_k-1}}} \nonumber \\
&\times & 
\frac{(\varphi_{(\mu_k V)_k})(\varphi_{(\mu_{k+1}\mu_{k}V)_{k+1}})
+(\varphi_{V_{r+k-\lfloor \frac{r}{2} \rfloor-1}})(\varphi_{V_{r+k-\lfloor \frac{r}{2} \rfloor}})(\varphi_{V_{r+k-\lfloor \frac{r}{2} \rfloor+1}})}{(\varphi_{V_{k-\lfloor \frac{r}{2} \rfloor }})}
\circ x^G_{{\rm \bf{i}}} (1;\phi(\textbf{Y})) \nonumber \\
&=&(\Phi_H(\textbf{Y}))^{\Lambda_{j_k-1}+\Lambda_{j_k-3}+\Lambda_{j_k+1}}\times
(\varphi_{(\mu[0]V)_{k-\lfloor \frac{r}{2} \rfloor }})
\circ x^G_{{\rm \bf{i}}} (1;\phi(\textbf{Y})) \label{mainprop-pr1} \\
&=&a^{\Lambda_{j_k-1}+\Lambda_{j_k-3}+\Lambda_{j_k+1}} Y_{1,j_k-1}Y_{2,j_k-1}Y_{1,j_k-3}Y_{2,j_k-3}Y_{1,j_k+1}Y_{2,j_k+1}\times \nonumber \\
& &(\varphi_{(\mu[0]V)_{k-\lfloor \frac{r}{2} \rfloor }})
\circ x^G_{{\rm \bf{i}}} (1;\phi(\textbf{Y}))\nonumber.
\end{eqnarray}

The module $(\mu[0]V)_{k-\lfloor \frac{r}{2} \rfloor }$ is described as (\ref{mutlem1-pr2}). Using Proposition \ref{dualsemi} and (\ref{iatoia}), let us calculate $(\varphi_{(\mu[0]V)_{k-\lfloor \frac{r}{2} \rfloor }})
\circ x^G_{{\rm \bf{i}}} (1;\phi(\textbf{Y}))$, that is, let us find 
\[{\rm \bf{a}}=(a_{1,j_{r}},\cdots,a_{1,j_{2}},a_{1,j_{1}},
a_{2,j_r},\cdots, a_{2,j_2},a_{2,j_1})\in (\mathbb{Z}_{\geq0})^{2r}\]
 satisfying $\mathcal{F}_{{\rm \bf{i}}^{\rm \bf{a}},(\mu[0]V)_{k-\lfloor \frac{r}{2} \rfloor }}\neq \phi$ ( or equivalently, $\mathcal{F}_{{\rm \bf{i}},{\rm \bf{a}},(\mu[0]V)_{k-\lfloor \frac{r}{2} \rfloor }}\neq \phi$). If $\mathcal{F}_{{\rm \bf{i}}^{\rm \bf{a}},(\mu[0]V)_{k-\lfloor \frac{r}{2} \rfloor }}\neq \phi$, by counting the number of the bases in (\ref{mutlem1-pr2}), since the dimension at $j_k+1$ is $3$, we have $a_{1,j_k+1}+a_{2,j_k+1}=3$. Considering similarly,
\begin{eqnarray}\label{count}
& &a_{1,j_k+1}+a_{2,j_k+1}=a_{1,j_k-1}+a_{2,j_k-1}=a_{1,j_k-3}+a_{2,j_k-3}=3,\nonumber\\
& &a_{1,j_k}+a_{2,j_k}=a_{1,j_k-2}+a_{2,j_k-2}=2,\\
& &a_{1,j_k+3}+a_{2,j_k+3}=a_{1,j_k+2}+a_{2,j_k+2}=a_{1,j_k-4}+a_{2,j_k-4}=a_{1,j_k-5}+a_{2,j_k-5}=1\nonumber.
\end{eqnarray}
Since the module $(\mu[0]V)_{k-\lfloor \frac{r}{2} \rfloor }$ does not have the simple submodules isomorphic to $S_{j_r},S_{j_{r-1}},\cdots,S_{j_{\lfloor \frac{r+1}{2} \rfloor+1}}$, we have $a_{1,j_r}=a_{1,j_{r-1}}=\cdots=a_{1,j_{\lfloor \frac{r+1}{2} \rfloor+1}}=0$, which yields $a_{2,j_k-4}=1$, $a_{2,j_k-2}=2$, $a_{2,j_k}=2$ and $a_{2,j_k+2}=1$. We can also check that $a_{1,j_k-3}=a_{1,j_k-1}=a_{1,j_k+1}=1$. Thus, $\mathcal{F}_{{\rm \bf{i}}^{\rm \bf{a}},(\mu[0]V)_{k-\lfloor \frac{r}{2} \rfloor }}\neq \phi$ if and only if
 \begin{multline}\label{ia1} {\rm \bf{i}}^{\rm \bf{a}}=(j_k-3,j_k-1,j_k+1,j_k-4,j_k-2,j_k-2,j_k,j_k,j_k+2,\\
j_k-5,j_k-3,j_k-3,j_k-1,j_k-1,j_k+1,j_k+1,j_k+3). \end{multline}
Then we can check that $\mathcal{F}_{{\rm \bf{i}},{\rm \bf{a}},(\mu[0]V)_{k-\lfloor \frac{r}{2} \rfloor }}$ is a point. Here, we use the notation as in $(\ref{iatoia})$. By the above argument and (\ref{mbasea}), (\ref{mbase0}), (\ref{mbase01}), we have
\begin{eqnarray*}
& & \hspace{-25pt}(\varphi^G_{(\mu[0]\mathbb{V})})_{k-\lfloor \frac{r}{2} \rfloor }(a;\textbf{Y})=a^{\Lambda_{j_k-1}+\Lambda_{j_k-3}+\Lambda_{j_k+1}}Y_{1,j_k-1}Y_{2,j_k-1}Y_{1,j_k-3}Y_{2,j_k-3}Y_{1,j_k+1}Y_{2,j_k+1}
 \\
& &\times \Phi_{1,j_k-3}(\textbf{Y})\Phi_{1,j_k-1}(\textbf{Y})\Phi_{1,j_k+1}(\textbf{Y})
\Phi_{2,j_k-4}(\textbf{Y})\Phi_{2,j_k-2}^2(\textbf{Y})\Phi_{2,j_k}^2(\textbf{Y})
\Phi_{2,j_k}^2(\textbf{Y})\\
& &\times \Phi_{2,j_k+2}(\textbf{Y})\Phi_{2,j_k-5}(\textbf{Y})
\Phi^2_{2,j_k-3}(\textbf{Y})\Phi^2_{2,j_k-1}(\textbf{Y})\Phi^2_{2,j_k+1}(\textbf{Y})
\Phi_{2,j_k+3}(\textbf{Y})\\
& &= a^{\Lambda_{j_k-1}+\Lambda_{j_k-3}+\Lambda_{j_k+1}}Y_{2,j_k-1} ,
\end{eqnarray*}
which implies the claim (a) for $l=0$.

Next, let us consider the claim (b) for $l=0$. By Lemma \ref{mutgamlem}, the arrows between $(\varphi_{(\mu[0]\mathbb{V})})_{k+1}=(\varphi_{(\mu_{k+1}\mathbb{V})})_{k+1}$ and $(\varphi_{(\mu[0]\mathbb{V})})_{s}$ $(s\in[-r,-1]\cup[1,2r])$ in $\mu[0](\Gamma_{\textbf{i}})$.  are
\[
\begin{xy}
(90,100)*{(\varphi_{\mathbb{V}})_{r+k-\lfloor \frac{r}{2} \rfloor -1}}="7",
(30,100)*{(\varphi_{\mathbb{V}})_{r+k+1}}="8",
(0,93) *{(\varphi_{\mathbb{V}})_{k-\lfloor \frac{r}{2} \rfloor+1}}="4",
(60,93) *{(\varphi_{(\mu[0]\mathbb{V})})_{k-\lfloor \frac{r}{2} \rfloor}}="3",
(30,93) *{(\varphi_{(\mu_{k+1}\mathbb{V})})_{k+1}}="2",
(30,86) *{(\varphi_{\mathbb{V}})_{-j_k+2}}="2d",
(60,86) *{(\varphi_{\mathbb{V}})_{-j_k+1}}="3d",
\ar@{->} "2";"4"
\ar@{->} "2";"7"
\ar@{->} "2";"8"
\ar@{->} "3";"2"
\ar@{->} "2d";"2"
\ar@{->} "2";"3d"
\end{xy}
\]
Thus, 
\begin{eqnarray*}
& &\hspace{-20pt} (\varphi^G_{(\mu_{k+1}\mu[0]\mathbb{V})})_{k+1}(a;\textbf{Y}) \\
&=&
\frac{(\varphi^G_{(\mu[0]\mathbb{V})})_{k-\lfloor \frac{r}{2} \rfloor}(\varphi^G_{\mathbb{V}})_{-j_k+2}+(\varphi^G_{\mathbb{V}})_{r+k-\lfloor \frac{r}{2} \rfloor -1}(\varphi^G_{\mathbb{V}})_{r+k+1}(\varphi^G_{\mathbb{V}})_{k-\lfloor \frac{r}{2} \rfloor+1}(\varphi^G_{\mathbb{V}})_{-j_k+1}}{(\varphi^G_{(\mu_{k+1}\mathbb{V})})_{k+1}}.
\end{eqnarray*}
Using (\ref{genbasic}), (\ref{gtol}), (\ref{muk+1}), (\ref{mainprop-pr1})
 and Theorem \ref{GLSthm} (ii), we obtain
\begin{eqnarray}
& &\hspace{-20pt} (\varphi^G_{(\mu_{k+1}\mu[0]\mathbb{V})})_{k+1}(a;\textbf{Y})=(\Phi_H(a;\textbf{Y}))^{\Lambda_{j_k-2}+\Lambda_{j_k+1}} \nonumber\\
&\times&
\frac{(\varphi_{(\mu[0]V)_{k-\lfloor \frac{r}{2} \rfloor}})+(\varphi_{V_{r+k-\lfloor \frac{r}{2} \rfloor -1}})(\varphi_{V_{r+k+1}})(\varphi_{V_{k-\lfloor \frac{r}{2} \rfloor+1}})}{(\varphi_{(\mu_{k+1}V)_{k+1}})}\circ x^G_{{\rm \bf{i}}} (1;\phi(\textbf{Y}))\nonumber\\
&=&a^{\Lambda_{j_k-2}+\Lambda_{j_k+1}} Y_{1,j_k-2}Y_{2,j_k-2}Y_{1,j_k+1}Y_{2,j_k+1}
\times
(\varphi_{(\mu_{k+1}\mu[0]V)_{k+1}})\circ x^G_{{\rm \bf{i}}} (1;\phi(\textbf{Y})).\qquad\quad \label{k1-1}
\end{eqnarray}

Applying a similar argument as in (\ref{ia1}) to the module $(\mu_{k+1}\mu[0]V)_{k+1}$ in (\ref{mutlem1-pr3}), for ${\rm \bf{a}}\in(\mathbb{Z}_{\geq 0})^{2r}$, we find that $\mathcal{F}_{{\rm \bf{i}}^{\rm \bf{a}}, (\mu_{k+1}\mu[0]V)_{k+1}}\neq \phi$ if and only if
 \begin{eqnarray*} 
{\rm \bf{i}}^{\rm \bf{a}}&=&(j_k+1,j_k-2,j_k,j_k+2,
j_k-3,j_k-1,j_k+1,j_k+3),\\
 & & (j_k-2,j_k+1,j_k,j_k+2,
j_k-3,j_k-1,j_k+1,j_k+3),\\
 & & (j_k-2,j_k-3,j_k+1,j_k,j_k+2,j_k-1,j_k+1,j_k+3). 
\end{eqnarray*}

Therefore, it follows from (\ref{phi12})
 and Proposition \ref{dualsemi} that
\begin{eqnarray}
& &\hspace{-20pt}(\varphi_{(\mu_{k+1}\mu[0]V)_{k+1}})\circ x^G_{{\rm \bf{i}}} (1;\phi(\textbf{Y}))\nonumber\\
&=& \Phi_{1,j_k+1}\Phi_{2,j_k-2}\Phi_{2,j_k}\Phi_{2,j_k+2}
\Phi_{2,j_k-3}\Phi_{2,j_k-1}\Phi_{2,j_k+1}\Phi_{2,j_k+3}\nonumber\\
&+& \Phi_{1,j_k-2}\Phi_{1,j_k+1}\Phi_{2,j_k}\Phi_{2,j_k+2}
\Phi_{2,j_k-3}\Phi_{2,j_k-1}\Phi_{2,j_k+1}\Phi_{2,j_k+3}\nonumber\\
&+& \Phi_{1,j_k-2}\Phi_{1,j_k-3}\Phi_{1,j_k+1}\Phi_{2,j_k}\Phi_{2,j_k+2}
\Phi_{2,j_k-1}\Phi_{2,j_k+1}\Phi_{2,j_k+3}\nonumber\\
&=& \Phi_{1,j_k+1}\Phi_{2,j_k-2}\Phi_{2,j_k}\Phi_{2,j_k+2}
\Phi_{2,j_k-3}\Phi_{2,j_k-1}\Phi_{2,j_k+1}\Phi_{2,j_k+3}
(1+A^{-1}_{1,j_k-2}+A^{-1}_{1,j_k-2}A^{-1}_{1,j_k-3})\nonumber\\
&=& \frac{Y_{2,j_k-1}}{Y_{2,j_k-2}Y_{1,j_k+1}Y_{2,j_k+1}}(1+A^{-1}_{1,j_k-2}+A^{-1}_{1,j_k-2}A^{-1}_{1,j_k-3}). \label{k1-2}
\end{eqnarray}

Substituting (\ref{k1-2}) for (\ref{k1-1}), we obtain
\[(\varphi^G_{(\mu_{k+1}\mu[0]\mathbb{V})})_{k+1}(a;\textbf{Y})
=a^{\Lambda_{j_k-2}+\Lambda_{j_k+1}}Y_{1,j_k-2}Y_{2,j_k-1}(1+A^{-1}_{1,j_k-2}+A^{-1}_{1,j_k-2}A^{-1}_{1,j_k-3}),
\]
which means the claim (b) for $l=0$.

Next, assuming that the claims (a),(b) for $0,1,\cdots,l$, let us prove the claims for $l+1$. Using Lemma \ref{mut0}, we have $(\varphi_{(\mu_{k+l+2}\mu_{k+l+1}\mu[l]\mathbb{V})})_{k+l+2}=(\varphi_{(\mu_{k+l+2}\mathbb{V})})_{k+l+2}$. By Lemma \ref{mutgamlem}, we see that the arrows between $(\varphi_{\mathbb{V}})_{k-\lfloor \frac{r}{2} \rfloor+l+1}$ and $(\varphi_{\mathbb{V}})_{-s}$ $(s\in[1,r])$ in $\mu_{k+l+2}\mu_{k+l+1}\mu[l]\Gamma_{{\rm \bf{i}}}$ are as follows:
\[
\begin{xy}
(70,82)*{(\varphi_{\mathbb{V}})_{-(j_k-2l-2)}}="-j_k",
(40,90) *{(\varphi_{\mathbb{V}})_{k-\lfloor \frac{r}{2} \rfloor +l+1}}="4",
(40,82)*{(\varphi_{\mathbb{V}})_{-(j_k-2l-3)}}="-2",
\ar@{->} "4";"-2"
\ar@{->} "-j_k";"4"
\end{xy}
\]
It follows from the exchange sequence (\ref{exseq-af1}), Lemma \ref{titj} and Theorem \ref{GLSthm} that the arrows from $(\varphi_{(\mu_{k+l+2}\mu_{k+l+1}\mu[l]\mathbb{V})})_{s}$ $(1\leq s\leq 2r)$ to $(\varphi_{\mathbb{V}})_{k-\lfloor \frac{r}{2} \rfloor+l+1}$ are
\[ (\varphi_{(\mu_{k+l+2}\mathbb{V})})_{k+l+2}\rightarrow (\varphi_{\mathbb{V}})_{k-\lfloor \frac{r}{2} \rfloor+l+1},\quad
(\varphi_{(\mu[l]\mathbb{V})})_{k-\lfloor \frac{r}{2} \rfloor+l} \rightarrow (\varphi_{\mathbb{V}})_{k-\lfloor \frac{r}{2} \rfloor+l+1}.\]
Similarly, since there exist non-factorizable homomorphisms in the direct summands of $(\mu_{k+l+2}\mu_{k+l+1}\mu[l]V)$ from $(\mu_{k+l+1}\mu[l]V)_{k+l+1}$, $V_{r+k-\lfloor \frac{r}{2} \rfloor +l+2}$, $V_{r+k-\lfloor \frac{r}{2} \rfloor +l+1}$ and $V_{r+k-\lfloor \frac{r}{2} \rfloor +l}$ to
$V_{k-\lfloor \frac{r}{2} \rfloor+l+1}=S_{j_k-2l-3}$, we see that the arrows in $(\mu_{k+l+2}\mu_{k+l+1}\mu[l]\Gamma_{{\rm \bf{i}}})$ from $(\varphi_{\mathbb{V}})_{k-\lfloor \frac{r}{2} \rfloor+l+1}$ to $(\varphi_{(\mu_{k+l+2}\mu_{k+l+1}\mu[l]\mathbb{V})})_{s}$ $(1\leq s\leq 2r)$ are
\[
(\varphi_{\mathbb{V}})_{k-\lfloor \frac{r}{2} \rfloor+l+1}\rightarrow
(\varphi_{(\mu_{k+l+1}\mu[l]\mathbb{V})})_{k+l+1},
\quad
(\varphi_{\mathbb{V}})_{k-\lfloor \frac{r}{2} \rfloor+l+1}\rightarrow
(\varphi_{\mathbb{V}})_{r+k-\lfloor \frac{r}{2} \rfloor +l+2},
\]
\[
(\varphi_{\mathbb{V}})_{k-\lfloor \frac{r}{2} \rfloor+l+1}\rightarrow
(\varphi_{\mathbb{V}})_{r+k-\lfloor \frac{r}{2} \rfloor +l+1},\quad
(\varphi_{\mathbb{V}})_{k-\lfloor \frac{r}{2} \rfloor+l+1}\rightarrow
(\varphi_{\mathbb{V}})_{r+k-\lfloor \frac{r}{2} \rfloor +l}.
\]
Hence, by the induction hypothesis of the claim (b), we obtain the following by the same way as in (\ref{mainprop-pr1}):
\begin{eqnarray}
(\varphi^G_{(\mu[l+1]\mathbb{V})})_{k-\lfloor \frac{r}{2} \rfloor+l+1}
&=&(\Phi_H(a;\textbf{Y}))^{(\sum^{l+3}_{s=0}\Lambda_{j_k-2s+1})+(\sum^{l}_{s=0}\Lambda_{j_k-2s-2})}\nonumber \\
&\times&\varphi_{(\mu[l+1]V)_{k-\lfloor \frac{r}{2} \rfloor+l+1}}\circ x^G_{{\rm \bf{i}}} (1;\phi(\textbf{Y})). \label{pr-l+1a}
\end{eqnarray}

The module $(\mu[l+1]V)_{k-\lfloor \frac{r}{2} \rfloor+l+1}$ is described as $(\ref{mutlem1-pr5})$. Using Proposition \ref{dualsemi}, let us calculate $\varphi_{(\mu[l+1]V)_{k-\lfloor \frac{r}{2} \rfloor+l+1}}\circ x^G_{{\rm \bf{i}}} (1;\phi(\textbf{Y}))$. 

For ${\rm \bf{a}}=(a_{1,j_r},\cdots,a_{1,j_1},a_{2,j_r},\cdots,a_{1,j_1})\in(\mathbb{Z}_{\geq 0})^{2r}$, if the variety $\mathcal{F}_{{\rm \bf{i}},{\rm \bf{a}},(\mu[l+1]V)_{k-\lfloor \frac{r}{2} \rfloor+l+1}}$ is non-empty, we have
\begin{eqnarray}
& &a_{1,j_k+1}=a_{2,j_{k}+2}=a_{2,j_k+3}=1,\ \ a_{1,j_k-1}=1,\ a_{2,j_{k}}=2,\ a_{2,j_k+1}=2,\ a_{2,j_k-1}=3, \nonumber\\
& &a_{1,j_{k+l}-5}=a_{2,j_{k+l}-6}=a_{2,j_{k+l}-7}=1,\nonumber\\
& &a_{1,j_{k+l}-3}=1,\ a_{2,j_{k+l}-4}=2,\ a_{2,j_{k+l}-5}=2,\ a_{2,j_{k+l}-3}=3,\label{forl+1}\\ 
& &a_{1,j_k-2t-3}+a_{2,j_k-2t-3}=5\ \ (0\leq t\leq l-1),\nonumber\\
& &a_{1,j_k-2t-2}+a_{2,j_k-2t-2}=3\ \ (0\leq t\leq l),\nonumber
\end{eqnarray}
by the same argument in the proof of the claim (a) for $l=0$. Denoting the bases in $(\ref{mutlem1-pr5})$ by
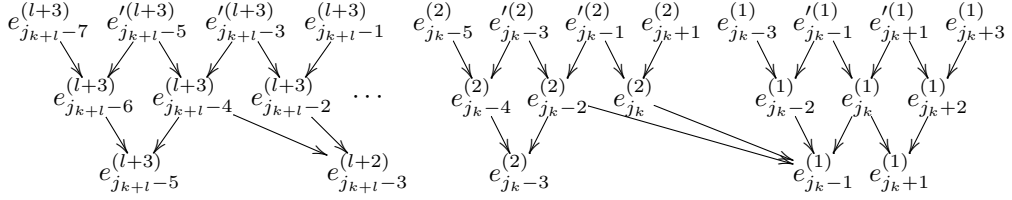
\[
\begin{xy}
(75,123)*{e^{(1)}_{j_k-3}}="a6",(85,123)*{e'^{(1)}_{j_k-1}}="a7",(95,123)*{e'^{(1)}_{j_k+1}}="a8",
(105,123)*{e^{(1)}_{j_k+3}}="a9",
(80,113)*{e^{(1)}_{j_k-2}}="a3",(90,113)*{e^{(1)}_{j_k}}="a4",(100,113)*{e^{(1)}_{j_k+2}}="a5",
(85,103)*{e^{(1)}_{j_k-1}}="a1",(95,103)*{e^{(1)}_{j_k+1}}="a2",
(35,123)*{e^{(2)}_{j_k-5}}="b6",(45,123)*{e'^{(2)}_{j_k-3}}="b7",(55,123)*{e'^{(2)}_{j_k-1}}="b8",
(65,123)*{e^{(2)}_{j_k+1}}="b9",
(40,113)*{e^{(2)}_{j_k-4}}="b3",(50,113)*{e^{(2)}_{j_k-2}}="b4",(60,113)*{e^{(2)}_{j_k}}="b5",
(45,103)*{e^{(2)}_{j_k-3}}="b1",
(-17,123)*{e^{(l+3)}_{j_{k+l}-7}}="c6",(-4,123)*{e'^{(l+3)}_{j_{k+l}-5}}="c7",(9,123)*{e'^{(l+3)}_{j_{k+l}-3}}="c8",
(22,123)*{e^{(l+3)}_{j_{k+l}-1}}="c9",
(-11,113)*{e^{(l+3)}_{j_{k+l}-6}}="c3",(2,113)*{e^{(l+3)}_{j_{k+l}-4}}="c4",(15,113)*{e^{(l+3)}_{j_{k+l}-2}}="c5",
(-5,103)*{e^{(l+3)}_{j_{k+l}-5}}="c1",(25,103)*{e^{(l+2)}_{j_{k+l}-3}}="z1",(25,113)*{\cdots}="dot",
\ar@{->} "a3";"a1"
\ar@{->} "a4";"a1"
\ar@{->} "a4";"a2"
\ar@{->} "a5";"a2"
\ar@{->} "a6";"a3"
\ar@{->} "a7";"a3"
\ar@{->} "a7";"a4"
\ar@{->} "a8";"a4"
\ar@{->} "a8";"a5"
\ar@{->} "a9";"a5"
\ar@{->} "b3";"b1"
\ar@{->} "b4";"b1"
\ar@{->} "b4";"a1"
\ar@{->} "b5";"a1"
\ar@{->} "b6";"b3"
\ar@{->} "b7";"b3"
\ar@{->} "b7";"b4"
\ar@{->} "b8";"b4"
\ar@{->} "b8";"b5"
\ar@{->} "b9";"b5"
\ar@{->} "c3";"c1"
\ar@{->} "c4";"c1"
\ar@{->} "c4";"z1"
\ar@{->} "c5";"z1"
\ar@{->} "c6";"c3"
\ar@{->} "c7";"c3"
\ar@{->} "c7";"c4"
\ar@{->} "c8";"c4"
\ar@{->} "c8";"c5"
\ar@{->} "c9";"c5"
\end{xy}
\]
we see that all $1$-dimensional simple submodule of $(\mu[l+1]V)_{k-\lfloor \frac{r}{2} \rfloor+l+1}$ are $\mathbb{C}e^{(t)}_{j_k-2t+1}$ $(1\leq t\leq l+3)$, $\mathbb{C}e^{(1)}_{j_k+1}$ and $\mathbb{C}(e^{(t+1)}_{j_k-2t-2}-e^{(t+2)}_{j_k-2t-2}+e^{(t+3)}_{j_k-2t-2})$ $(0\leq t\leq l)$, which are isomorphic to $S_{j_k-2t+1}$, $S_{j_k+1}$ and $S_{j_k-2t-2}$, respectively. We can also see that if $0\leq t\leq l-1$, then all $1$-dimensional simple submodule isomorphic to $S_{j_k-2t-3}$ of the quotient module 
\[
\frac{(\mu[l+1]V)_{k-\lfloor \frac{r}{2} \rfloor+l+1}}
{(\mathbb{C}(e^{(t+1)}_{j_k-2t-2}-e^{(t+2)}_{j_k-2t-2}+e^{(t+3)}_{j_k-2t-2})\oplus
\mathbb{C}(e^{(t+2)}_{j_k-2t-4}-e^{(t+3)}_{j_k-2t-4}+e^{(t+4)}_{j_k-2t-4}))}
\]
are $\mathbb{C}e^{(t+2)}_{j_k-2t-3}$ and
$\mathbb{C}(e^{(t+1)}_{j_k-2t-3}-e'^{(t+2)}_{j_k-2t-3}+e'^{(t+3)}_{j_k-2t-3}-e^{(t+4)}_{j_k-2t-3})$. Thus, $\mathcal{F}_{{\rm \bf{i}},{\rm \bf{a}},(\mu[l+1]V)_{k-\lfloor \frac{r}{2} \rfloor+l+1}}\neq \phi$ if and only if ${\rm \bf{a}}$ satisfies the following in addition to (\ref{forl+1}): For each $t\in[0,l-1]$,
\begin{eqnarray*}
& &\hspace{-20pt}(a_{1,j_k-2t-2},a_{1,j_k-2t-3},a_{1,j_k-2t-4})=(1,2,1),(0,1,0),(0,1,1),(1,1,0)\ {\rm or}\ (1,1,1)\ \\
& &\hspace{-20pt} {\rm and\ all\ other}\ a_{1,i},\ a_{2,i}=0.
\end{eqnarray*}
Let us calculate the monomial $M$ corresponding to $(a_{1,j_k-2t-2},a_{1,j_k-2t-3},a_{1,j_k-2t-4})=(0,1,0)$ for all $t\in[0,l-1]$, which means that $a_{1,j_k-2s+1}=1$ $(0\leq s\leq l+3)$ and $a_{1,j_k-2s+2}=0$ $(0\leq s\leq l+4)$. Thus, it is calculated as
\begin{eqnarray}
M&=& \Phi_{1,j_k+1}\nonumber \times \prod^{l+2}_{s=0} (\Phi_{1,j_k-2s-1} \Phi_{2,j_k-2s-2} \Phi_{2,j_k-2s-3} \\
& &\Phi_{2,j_k-2s} \Phi_{2,j_k-2s-1}
\Phi_{2,j_k-2s+2} \Phi_{2,j_k-2s+1} \Phi_{2,j_k-2s+3}) \nonumber \\
&=& \frac{Y_{2,j_k}Y_{2,j_k+2}}{Y_{1,j_k+1}Y^2_{2,j_k+1}} \prod^{l+2}_{s=0} \frac{Y_{2,j_k-2s+1}}{Y_{1,j_k-2s-1}Y_{2,j_k-2s-1}Y_{2,j_k-2s+2}}. \label{primary}
\end{eqnarray}
For $p\geq 1$ and $({\rm \bf{b}},{\rm \bf{c}})\in R^p_{l}$ $({\rm \bf{b}}=\{b_i\}^p_{i=1},\ {\rm \bf{c}}=\{c_i\}^p_{i=1})$, the monomial corresponding to $a_{1,j_k-2t-3}=2,\ a_{1,j_k-2t-2}=a_{1,j_k-2t-4}=1$ for $t\in [b_1,c_1-1]\cup \cdots \cup[b_p,c_p-1]$, and $a_{1,j_k-2t-3}=1$ for $t\in [1,l-1]\setminus ([b_1,c_1-1]\cup \cdots \cup[b_p,c_p-1])$, and $a_{1,j_k-2t-2}=0$ for $t\in[0,l]\setminus [{\rm \bf{b}},{\rm \bf{c}}]$ is $M\times
A[b_1,c_1;j_k]\cdots A[b_p,c_p;j_k]$ by (\ref{phi12}). Using (\ref{phi12}) again, we see that the partial sum of $\varphi_{(\mu[l+1]V)_{k-\lfloor \frac{r}{2} \rfloor+l+1}}\circ x^G_{{\rm \bf{i}}} (1;\phi(\textbf{Y}))$ corresponding to $a_{1,j_k-2t-3}=2,\ a_{1,j_k-2t-2}=a_{1,j_k-2t-4}=1$ for $t\in [b_1,c_1-1]\cup \cdots \cup[b_p,c_p-1]$ and $a_{1,j_k-2t-2}=0$ or $1$ for $t\in[0,l]\setminus[{\rm \bf{b}},{\rm \bf{c}}]$ is
\[
M \times A[b_1,c_1;j_k]\cdots A[b_p,c_p;j_k] \prod_{t\in [0,l]\setminus ([b_1,c_1]\cup \cdots\cup[b_p,c_p])} (1+A^{-1}_{1,j_k-2t-2}).
\]
On the other hand, by (\ref{mbasea}) and (\ref{primary}),
\begin{multline*}
(\Phi_H(a;\textbf{Y}))^{(\sum^{l+3}_{s=0}\Lambda_{j_k-2s+1})+(\sum^{l}_{s=0}\Lambda_{j_k-2s-2})}
\times M\\
=a^{(\sum^{l+3}_{s=0}\Lambda_{j_k-2s+1})+(\sum^{l}_{s=0}\Lambda_{j_k-2s-2})} \prod^{l}_{s=0} Y_{1,j_k-2s-2} \prod^{l+1}_{s=0} Y_{2,j_k-2s-1}.\qquad \qquad
\end{multline*}
Hence, by (\ref{pr-l+1a}), we have
\[
(\varphi^G_{(\mu[l+1]\mathbb{V})})_{k-\lfloor \frac{r}{2} \rfloor+l+1}
=a^{(\sum^{l+3}_{s=0}\Lambda_{j_k-2s+1})+(\sum^{l}_{s=0}\Lambda_{j_k-2s-2})}
\prod^{l}_{s=0} Y_{1,j_k-2s-2} \prod^{l+1}_{s=0} Y_{2,j_k-2s-1}\]
\[\times \sum_{p\geq 0,\ ({\rm \bf{b}},{\rm \bf{c}})\in R^p_{l}}
A[b_1,c_1;j_k]\cdots A[b_p,c_p;j_k] \prod_{t\in [0,l]\setminus ([b_1,c_1]\cup \cdots\cup[b_p,c_p])} (1+A^{-1}_{1,j_k-2t-2}),
\]
which implies the claim (a) for $l+1$.

Finally, let us prove the claim (b) for $l+1$. By the direct calculation, the arrows between $(\varphi_{(\mu_{k+l+2}\mu_{k+l+1}\mu[l]\mathbb{V})})_{k+l+2}=(\varphi_{(\mu_{k+l+2}\mathbb{V})})_{k+l+2}$ and $(\varphi_{\mathbb{V}})_{-s}$ $(s\in[1,r])$ in $\mu[l+1]\Gamma_{{\rm \bf{i}}}$ are as follows:
\[
\begin{xy}
(70,82)*{(\varphi_{\mathbb{V}})_{-(j_k-2l-3)}}="-j_k",
(40,90) *{(\varphi_{(\mu_{k+l+2}\mathbb{V})})_{k+l+2}}="4",
(40,82)*{(\varphi_{\mathbb{V}})_{-(j_k-2l-4)}}="-2",
\ar@{->} "-2";"4"
\ar@{->} "4";"-j_k"
\end{xy}
\]
By the exchange sequence $(\ref{exseq-2})$, Lemma \ref{titj} and Theorem \ref{GLSthm} imply that the arrow from $(\varphi_{(\mu[l+1]\mathbb{V})})_s$ $(1\leq s\leq 2r)$ to $(\varphi_{(\mu_{k+l+2}\mathbb{V})})_{k+l+2}$ are
\[ (\varphi_{(\mu[l+1]\mathbb{V})})_{k-\lfloor \frac{r}{2} \rfloor+l+1}
\rightarrow (\varphi_{(\mu_{k+l+2}\mathbb{V})})_{k+l+2}. \]
The arrows from $(\varphi_{(\mu_{k+l+2}\mathbb{V})})_{k+l+2}$ to $(\varphi_{(\mu[l+1]\mathbb{V})})_s$ $(1\leq s\leq2r)$ are
\[
(\varphi_{(\mu_{k+l+2}\mathbb{V})})_{k+l+2}\rightarrow
(\varphi_{(\mu_{k+l+1}\mu[l]\mathbb{V})})_{k+l+1},\ \ 
(\varphi_{(\mu_{k+l+2}\mathbb{V})})_{k+l+2}\rightarrow
(\varphi_{\mathbb{V}})_{r+k+l-\lfloor \frac{r}{2} \rfloor},
\]
\[
(\varphi_{(\mu_{k+l+2}\mathbb{V})})_{k+l+2}\rightarrow
(\varphi_{\mathbb{V}})_{r+k+l+2},\ \ 
(\varphi_{(\mu_{k+l+2}\mathbb{V})})_{k+l+2}\rightarrow
(\varphi_{\mathbb{V}})_{k+l-\lfloor \frac{r}{2} \rfloor+2}.
\]
By the same way as in (\ref{pr-l+1a}), we have
\begin{eqnarray}
(\varphi^G_{(\mu_{k+l+2}\mu[l+1]\mathbb{V})})_{k+l+2}
&=&(\Phi_H(a;\textbf{Y}))^{(\sum^{l+1}_{s=0}\Lambda_{j_k-2s+1}+\Lambda_{j_k-2s-2})} \nonumber\\
&\times&\varphi_{(\mu_{k+l+2}\mu[l+1]V)_{k+l+2}}\circ x^G_{{\rm \bf{i}}} (1;\phi(\textbf{Y}))\label{pr-l+1b}. 
\end{eqnarray}

The module $(\mu_{k+l+2}\mu[l+1]V)_{k+l+2}$ is described as (\ref{mutlem1-pr6}), and it has the simple submodules $S'_{j_k-2t}$ isomorphic to $S_{j_k-2t}$ $(1\leq t\leq l+2)$. The quotient modules $(\mu_{k+l+2}\mu[l+1]V)_{k+l+2}/(S'_{j_k-2t}\oplus S'_{j_k-2t-2})$ and $(\mu_{k+l+2}\mu[l+1]V)_{k+l+2}/(S'_{j_k-2l-4})$ have the simple submodules isomorphic $S_{j_k-2t-1}$ and $S_{j_k-2l-5}$ respectively. Therefore, for ${\rm \bf{a}}=(a_{1,j_r},\cdots,a_{1,j_1},a_{2,j_r},\cdots,a_{1,j_1})\in(\mathbb{Z}_{\geq 0})^{2r}$, the variety $\mathcal{F}_{{\rm \bf{i}},{\rm \bf{a}},(\mu_{k+l+2}\mu[l+1]V)_{k+l+2}}$ is non-empty if and only if 
\begin{eqnarray*}
& &a_{1,j_k+1}=a_{2,j_{k}+2}=a_{2,j_k+3}=1,\ 
a_{1,j_k-1}=1,\ a_{2,j_{k}}=a_{2,j_k+1}=2,\ a_{2,j_k-1}=3, \\
& &a_{1,j_k-2t-3}+a_{2,j_k-2t-3}=5\ \ (0\leq t\leq l-2),\\
& &a_{1,j_k-2t-2}+a_{2,j_k-2t-2}=3\ \ (0\leq t\leq l-1),\\
& &a_{1,j_k-2l-1}+a_{2,j_k-2l-1}=4,\quad 
a_{1,j_k-2l-3}+a_{2,j_k-2l-3}=2, \\
& &a_{1,j_k-2l-5}+a_{2,j_k-2l-5}=1, \quad
a_{1,j_k-2l-2}+a_{2,j_k-2l-2}=2,\\
& &a_{1,j_k-2l-4}+a_{2,j_k-2l-4}=1,
\end{eqnarray*}
and for each $t\in[0,l-1]$,
\begin{eqnarray*}
& &\hspace{-20pt}(a_{1,j_k-2t-2},a_{1,j_k-2t-3},a_{1,j_k-2t-4})=(1,2,1),(0,1,0),(0,1,1),(1,1,0)\ {\rm or}\ (1,1,1)\\
& &\hspace{-20pt}(a_{1,j_k-2l-2},a_{1,j_k-2l-3},a_{1,j_k-2l-4})=(1,1,1),(0,0,0),(0,0,1),(1,0,0)\ {\rm or}\ (1,0,1)\\
& &\hspace{-20pt}0\leq a_{1,j_{k}-2l-5} \leq a_{1,j_{k}-2l-4}\leq 1, \quad
{\rm all\ other}\ a_{1,i},\ a_{2,i}=0.
\end{eqnarray*}
Let us calculate the monomial $M'$ corresponding to 
$(a_{1,j_k-2t-2},a_{1,j_k-2t-3},a_{1,j_k-2t-4})=(0,1,0)$ for all $t\in[0,l-1]$ and $a_{1,j_k-2l-3}=a_{1,j_{k}-2l-4}=a_{1,j_{k}-2l-5}=0$. It is calculated as
\begin{eqnarray}
M'&=& \prod^{l+1}_{s=0} (\Phi_{1,j_k-2s+1} \Phi_{2,j_k-2s-2} \Phi_{2,j_k-2s-3} \nonumber\\
& &\Phi_{2,j_k-2s} \Phi_{2,j_k-2s-1}
\Phi_{2,j_k-2s+2} \Phi_{2,j_k-2s+1} \Phi_{2,j_k-2s+3}) \nonumber \\
&=&\prod^{l+1}_{s=0} \frac{Y_{2,j_k-2s-1}}{Y_{1,j_k-2s+1}Y_{2,j_k-2s+1}Y_{2,j_k-2s-2}}. \label{primary2}
\end{eqnarray}
For $p\geq 1$ and $({\rm \bf{b}},{\rm \bf{c}})\in R^p_{l+1}$ $({\rm \bf{b}}=\{b_i\}^p_{i=1},\ {\rm \bf{c}}=\{c_i\}^p_{i=1})$ such that $c_p<l+1$, the monomial corresponding to $a_{1,j_k-2t-3}=2,\ a_{1,j_k-2t-2}=a_{1,j_k-2t-4}=1$ for $t\in [b_1,c_1-1]\cup \cdots \cup[b_p,c_p-1]$, $a_{1,j_k-2t-3}=1$ for $t\in [1,l-1]\setminus ([b_1,c_1-1]\cup \cdots \cup[b_p,c_p-1])$, $a_{1,j_k-2l-3}=a_{1,j_k-2l-5}=0$ and $a_{1,j_k-2t-2}=0$ for $t\in[0,l+1]\setminus [{\rm \bf{b}},{\rm \bf{c}}]$ is $M'\times
A[b_1,c_1;j_k]\cdots A[b_p,c_p;j_k]$ by (\ref{phi12}). Using (\ref{phi12}) again, we see that the partial sum of $\varphi_{(\mu_{k+l+2}\mu[l+1]V)_{k+l+2}}\circ x^G_{{\rm \bf{i}}} (1;\phi(\textbf{Y}))$ corresponding to $a_{1,j_k-2t-3}=2,\ a_{1,j_k-2t-2}=a_{1,j_k-2t-4}=1$ $(t\in [b_1,c_1-1]\cup \cdots \cup[b_p,c_p-1])$ and $a_{1,j_k-2t-2}=1$ or $0$ for $t\in[0,l+1]\setminus[{\rm \bf{b}},{\rm \bf{c}}]$ and $0\leq a_{1,j_{k}-2l-5} \leq a_{1,j_{k}-2l-4}\leq 1$
is
\begin{eqnarray*}
& &M' \times A[b_1,c_1;j_k]\cdots A[b_p,c_p;j_k] \prod_{t\in [0,l+1]\setminus ([b_1,c_1]\cup \cdots\cup[b_p,c_p])} (1+A^{-1}_{1,j_k-2t-2})\\
&\times& (1+A^{-1}_{1,j_k-2l-4}+A^{-1}_{1,j_k-2l-4}A^{-1}_{1,j_k-2l-5}).
\end{eqnarray*}

Similarly, for $p\geq 1$ and $({\rm \bf{b}},{\rm \bf{c}})\in R^p_{l+1}$ $({\rm \bf{b}}=\{b_i\}^p_{i=1},\ {\rm \bf{c}}=\{c_i\}^p_{i=1})$ such that $c_p=l+1$, the monomial corresponding to $a_{1,j_k-2t-3}=2,\ a_{1,j_k-2t-2}=a_{1,j_k-2t-4}=1$ for $t\in [b_1,c_1-1]\cup \cdots \cup[b_p,l]$, $a_{1,j_k-2t-3}=1$ for $t\in [1,l]\setminus ([b_1,c_1-1]\cup \cdots \cup[b_p,l-1])$, $a_{1,j_k-2l-3}=a_{1,j_k-2l-4}=1$, $a_{1,j_k-2l-5}=0$ and $a_{1,j_k-2t-2}=0$ for $t\in[0,l+1]\setminus ([b_1,c_1]\cup \cdots \cup[b_p,l+1])$ is $M'\times A[b_1,c_1;j_k]\cdots A[b_p,c_p;j_k]$. We see that the partial sum of $\varphi_{(\mu_{k+l+2}\mu[l+1]V)_{k+l+2}}\circ x^G_{{\rm \bf{i}}} (1;\phi(\textbf{Y}))$ corresponding to $a_{1,j_k-2t-3}=2,\ a_{1,j_k-2t-2}=a_{1,j_k-2t-4}=1$ $(t\in [b_1,c_1-1]\cup \cdots \cup[b_p,l-1])$ and $a_{1,j_k-2l-2}=a_{1,j_k-2l-3}=a_{1,j_k-2l-4}=1$ is $M' \times (1+A^{-1}_{1,j_k-2l-5})\times A[b_1,c_1;j_k]\cdots A[b_p,c_p;j_k] \prod_{t\in [0,l+1]\setminus ([b_1,c_1]\cup \cdots\cup[b_p,l+1])} (1+A^{-1}_{1,j_k-2t-2})$. On the other hand, by (\ref{mbasea}) and (\ref{primary2}), 
\begin{eqnarray*}
& &\hspace{-20pt}(\Phi_H(a;\textbf{Y}))^{(\sum^{l+1}_{s=0}\Lambda_{j_k-2s+1}+\Lambda_{j_k-2s-2})}
\times M'\\
&=&a^{(\sum^{l+1}_{s=0}\Lambda_{j_k-2s+1}+\Lambda_{j_k-2s-2})} \prod^{l+1}_{s=0} Y_{1,j_k-2s-2} \prod^{l+1}_{s=0} Y_{2,j_k-2s-1}.\qquad \qquad
\end{eqnarray*}

Hence, by (\ref{pr-l+1b}), we have
\begin{multline*}
(\varphi^G_{(\mu_{k+l+2}\mu[l+1]\mathbb{V})})_{k+l+2}
=a^{(\sum^{l+1}_{s=0}\Lambda_{j_k-2s+1}+\Lambda_{j_k-2s-2})}
\prod^{l+1}_{s=0} Y_{1,j_k-2s-2} \prod^{l+1}_{s=0} Y_{2,j_k-2s-1}\\
\times \sum_{p\geq 0,\ ({\rm \bf{b}},{\rm \bf{c}})\in R^p_{l+1}}
A[b_1,c_1;j_k]\cdots A[b_p,c_p;j_k](1-\delta_{c_p,l+1}+A^{-1+\delta_{c_p,l+1}}_{1,j_k-2l-4}(1+A^{-1}_{1,j_k-2l-5}))
\\
\times\prod_{t\in [0,l+1]\setminus ([{\rm \bf{b}},{\rm \bf{c}}])} (1+A^{-1}_{1,j_k-2t-2}),
\qquad \qquad \qquad \qquad \qquad \qquad \qquad \qquad \qquad \qquad \qquad
\end{multline*}
which implies the claim (b) for $l+1$. \qed

\vspace{2mm}

The proofs for Proposition \ref{spprop2} and \ref{spprop3} are similar to the one for Proposition \ref{spprop1}. 

\vspace{1mm}
\nd
{\sl [Proof of Theorem \ref{thm1}.]}

Let us prove (1) since (2) and (3) are proven in the same way as (1).

(a) For $p\geq0$ and $({\rm \bf{b}},{\rm \bf{c}})\in R^p_{l-1}$ $({\rm \bf{b}}=\{b_i\}^p_{i=1},{\rm \bf{c}}=\{c_i\}^p_{i=1})$, let $\mu_1:B\left(\left(\sum^{j_k-1}_{s=j_k-2l-1}\Lambda_s\right)-\al[{\rm \bf{b}},{\rm \bf{c}};j_k]\right)\rightarrow \mathcal{Y}$ be the monomial realization which maps the highest weight vector in $B\left(\left(\sum^{j_k-1}_{s=j_k-2l-1}\Lambda_s\right)-\al[{\rm \bf{b}},{\rm \bf{c}};j_k]\right)$ to the monomial $H_1[{\rm \bf{b}},{\rm \bf{c}}]:=H_1\cdot A[b_1,c_1;j_k]\cdots A[b_p,c_p;j_k]$, where $\mathcal{Y}$ is defined in \ref{monoreacry}. By Proposition \ref{spprop1}, we need show that
$H_1[{\rm \bf{b}},{\rm \bf{c}}]
\prod_{q\in[0,l-1]\setminus[{\rm \bf{b}},{\rm \bf{c}}]}(1+A_{1,j_k-2q-2}^{-1})$ coincides with
\begin{equation}\label{desire}
\sum_{b\in B\left(\left(\sum^{j_k-1}_{s=j_k-2l-1}\Lambda_s\right)-\al[{\rm \bf{b}},{\rm \bf{c}};j_k]\right)_{\prod_{q\in[0,l-1]\setminus[{\rm \bf{b}},{\rm \bf{c}}]}s_{j_k-2q-2}}}\mu_1(b).
\end{equation}

First, let us show that each factor in the monomial $H_1[{\rm \bf{b}},{\rm \bf{c}}]$ has non-negative degree. For $1\leq i\leq p$, we can easily see that
\begin{eqnarray}
A[b_i,c_i;j_k]&=&\left(\prod^{c_i-1}_{s=b_i}A^{-1}_{1,j_k-2s-2}A^{-1}_{1,j_k-2s-3} \right)
A^{-1}_{1,j_k-2c_i-2}\nonumber\\
&=&\left(\prod^{c_i-1}_{s=b_i}\frac{Y_{1,j_k-2s-1}Y_{2,j_k-2s-4}}{Y_{1,j_k-2s-2}Y_{2,j_k-2s-3}} \right)
\frac{Y_{1,j_k-2c_i-3}Y_{1,j_k-2c_i-1}}{Y_{1,j_k-2c_i-2}Y_{2,j_k-2c_i-2}}\nonumber\\
&=& \label{eachhigh1}
\frac{\prod^{c_i+1}_{s=b_i}Y_{1,j_k-2s-1}\prod^{c_i-1}_{s=b_i+1}Y_{2,j_k-2s-2}}
{\prod^{c_i}_{s=b_i}Y_{1,j_k-2s-2}\prod^{c_i}_{s=b_i+1}Y_{2,j_k-2s-1}}.
\end{eqnarray}
Hence, each factor in the monomial $H_1[{\rm \bf{b}},{\rm \bf{c}}]$ has non-negative degree. Since the monomial $A_{1,i}$ has the weight $\al_i$ (see (\ref{wtph}),(\ref{asidefsp})), the monomial $H_1[{\rm \bf{b}},{\rm \bf{c}}]$ has the weight $\left(\sum^{j_k-1}_{s=j_k-2l-1}\Lambda_s\right)-\al[{\rm \bf{b}},{\rm \bf{c}};j_k]$. Furthermore, by (\ref{eachhigh1}), we can verify that for $q\in[1,l-1]\setminus[{\rm \bf{b}},{\rm \bf{c}}]$, in the monomial $H_1[{\rm \bf{b}},{\rm \bf{c}}]$, the factor $Y_{1,j_k-2q-2}$ has the degree $1$, and the factor $Y_{2,j_k-2q-2}$ does not appear. Thus, the definition of Kashiwara operators in \ref{monoreacry} implies that
\[ \tilde{f}_{j_k-2q-2}H_1[{\rm \bf{b}},{\rm \bf{c}}]=
H_1[{\rm \bf{b}},{\rm \bf{c}}]\cdot A^{-1}_{1,j_k-2q-2},
\]
and $\tilde{f}^2_{j_k-2q-2}H_1[{\rm \bf{b}},{\rm \bf{c}}]=0$.
More generally, by the definition of the monomials $A_{1,i}$ $(i\in I)$ in (\ref{asidefsp}), for $q_1,\cdots,q_m\in[1,l-1]\setminus[{\rm \bf{b}},{\rm \bf{c}}]$ $(m\in \mathbb{Z}_{\geq0})$, if $q\in[1,l-1]\setminus[{\rm \bf{b}},{\rm \bf{c}}]$ and $q\neq q_1,\cdots,q_m$, then in the monomial $H_1[{\rm \bf{b}},{\rm \bf{c}}]\prod^m_{s=1}A^{-1}_{1,j_k-2q_s-2}$, the factor $Y_{1,j_k-2q-2}$ has the degree $1$, and factors $Y^{\pm1}_{2,j_k-2q-2}$ do not appear. Hence,
\[
\tilde{f}_{j_k-2q-2}(H_1[{\rm \bf{b}},{\rm \bf{c}}]\prod^m_{s=1}A^{-1}_{1,j_k-2q_s-2})=
(H_1[{\rm \bf{b}},{\rm \bf{c}}]\prod^m_{s=1}A^{-1}_{1,j_k-2q_s-2})
\cdot A^{-1}_{1,j_k-2q-2},
\]
and $\tilde{f}^2_{j_k-2q-2}(H_1[{\rm \bf{b}},{\rm \bf{c}}]\prod^m_{s=1}A^{-1}_{1,j_k-2q_s-2})=0$.

Let ${\rm id}_{\mathcal{Y}}$ be the identity map on the set $\mathcal{Y}$. By the above argument, we obtain
\[
H_1[{\rm \bf{b}},{\rm \bf{c}}]
\prod_{q\in[0,l-1]\setminus[{\rm \bf{b}},{\rm \bf{c}}]}(1+A_{1,j_k-2q-2}^{-1})
=\left(\prod_{q\in[0,l-1]\setminus[{\rm \bf{b}},{\rm \bf{c}}]}({\rm id}_{\mathcal{Y}}+\tilde{f}_{1,j_k-2q-2})\right)H_1[{\rm \bf{b}},{\rm \bf{c}}],\]
and from Theorem \ref{kashidem}, it coincides with (\ref{desire}).

\vspace{1mm}

(b) For $p\geq0$ and $({\rm \bf{b}},{\rm \bf{c}})\in R^p_{l}$ $({\rm \bf{b}}=\{b_i\}^p_{i=1},{\rm \bf{c}}=\{c_i\}^p_{i=1})$, let 
\[\mu_2:B\left(\left(\sum^{j_k-1}_{s=j_k-2l-2}\Lambda_s\right)-\al[{\rm \bf{b}},{\rm \bf{c}};j_k]\right)\rightarrow \mathcal{Y}\]
 be the monomial realization which maps the highest weight vector to the monomial $H_2[{\rm \bf{b}},{\rm \bf{c}}]:=H_2\cdot A[b_1,c_1;j_k]\cdots A[b_p,c_p;j_k]$. By Proposition \ref{spprop1}, we need show that
\begin{equation}\label{bmainmono}
H_2[{\rm \bf{b}},{\rm \bf{c}}]\cdot
(1-\delta_{c_p,l}+A^{-1+\delta_{c_p,l}}_{1,j_k-2l-2}(1+A^{-1}_{1,j_k-2l-3}))
\prod_{q\in[0,l-1]\setminus[{\rm \bf{b}},{\rm \bf{c}}]}(1+A_{1,j_k-2q-2}^{-1})
\end{equation}
coincides with
\begin{equation}\label{desireb}
\sum_{b\in B\left(\left(\sum^{j_k-1}_{s=j_k-2l-2}\Lambda_s\right)-\al[{\rm \bf{b}},{\rm \bf{c}};j_k]\right)_{s_{j_k-2l-3}s^{1-\delta_{c_p,l}}_{j_k-2l-2}\prod_{q\in[0,l-1]\setminus[{\rm \bf{b}},{\rm \bf{c}}]}s_{j_k-2q-2}}}\mu_2(b).
\end{equation}
In the same way as (a), we see that each factor in the monomial $H_2[{\rm \bf{b}},{\rm \bf{c}}]$ has non-negative power, and it has the weight $\left(\sum^{j_k-1}_{s=j_k-2l-2}\Lambda_s\right)-\al[{\rm \bf{b}},{\rm \bf{c}};j_k]$. For $q_1,\cdots,q_m\in[0,l-1]\setminus[{\rm \bf{b}},{\rm \bf{c}}]$ $(m\in \mathbb{Z}_{\geq0})$, if $q\in[0,l-1]\setminus[{\rm \bf{b}},{\rm \bf{c}}]$ and $q\neq q_1,\cdots,q_m$, then in the monomial $H_2[{\rm \bf{b}},{\rm \bf{c}}]\prod^m_{s=1}A^{-1}_{1,j_k-2q_s-2}$, the factor $Y_{1,j_k-2q-2}$ has the degree $1$, and factors $Y^{\pm1}_{2,j_k-2q-2}$ do not appear. Hence,
\[
\tilde{f}_{j_k-2q-2}(H_2[{\rm \bf{b}},{\rm \bf{c}}]\prod^m_{s=1}A^{-1}_{1,j_k-2q_s-2})=
(H_2[{\rm \bf{b}},{\rm \bf{c}}]\prod^m_{s=1}A^{-1}_{1,j_k-2q_s-2})
\cdot A^{-1}_{1,j_k-2q-2},
\]
and $\tilde{f}^2_{j_k-2q-2}(H_2[{\rm \bf{b}},{\rm \bf{c}}]\prod^m_{s=1}A^{-1}_{1,j_k-2q_s-2})=0$. Moreover, if $c_p<l$, we have 
\begin{equation}\label{ftilh2}
\tilde{f}_{j_k-2l-2}(H_2[{\rm \bf{b}},{\rm \bf{c}}]\prod^m_{s=1}A^{-1}_{1,j_k-2q_s-2})=
(H_2[{\rm \bf{b}},{\rm \bf{c}}]\prod^m_{s=1}A^{-1}_{1,j_k-2q_s-2})
\cdot A^{-1}_{1,j_k-2l-2}
\end{equation}
and $\tilde{f}^2_{j_k-2l-2}(H_2[{\rm \bf{b}},{\rm \bf{c}}]\prod^m_{s=1}A^{-1}_{1,j_k-2q_s-2})=0$. It follows from the explicit forms of $H_2[{\rm \bf{b}},{\rm \bf{c}}]$ and $A^{-1}_{1,j_k-2q_s-2}$ that in the monomial (\ref{ftilh2}), the factor $Y_{1,j_k-2l-3}$ has the degree $1$, and factors $Y^{\pm1}_{2,j_k-2l-3}$ do not appear. Hence, the definition of Kashiwara operators in \ref{monoreacry} implies that
\begin{eqnarray*}
& &\hspace{-20pt}\tilde{f}_{j_k-2l-3}((H_2[{\rm \bf{b}},{\rm \bf{c}}]\prod^m_{s=1}A^{-1}_{1,j_k-2q_s-2})A^{-1}_{1,j_k-2l-2})\\
&=&((H_2[{\rm \bf{b}},{\rm \bf{c}}]\prod^m_{s=1}A^{-1}_{1,j_k-2q_s-2})A^{-1}_{1,j_k-2l-2})
\cdot A^{-1}_{1,j_k-2l-3},
\end{eqnarray*} 
and $\tilde{f}^2_{j_k-2l-3}((H_2[{\rm \bf{b}},{\rm \bf{c}}]\prod^m_{s=1}A^{-1}_{1,j_k-2q_s-2})A^{-1}_{1,j_k-2l-2})=0$. Similarly, if $c_p=l$, we obtain $\tilde{f}_{j_k-2l-3}(H_2[{\rm \bf{b}},{\rm \bf{c}}]\prod^m_{s=1}A^{-1}_{1,j_k-2q_s-2})=(H_2[{\rm \bf{b}},{\rm \bf{c}}]\prod^m_{s=1}A^{-1}_{1,j_k-2q_s-2})A^{-1}_{1,j_k-2l-3}$, and $\tilde{f}^2_{j_k-2l-3}(H_2[{\rm \bf{b}},{\rm \bf{c}}]\prod^m_{s=1}A^{-1}_{1,j_k-2q_s-2})=0$. By the above argument, the sum in (\ref{bmainmono}) is the same as
\[ ((1-\delta_{c_p,l}){\rm id}_{\mathcal{Y}}+\tilde{f}^{-1+\delta_{c_p,l}}_{j_k-2l-2}({\rm id}_{\mathcal{Y}}+\tilde{f}_{j_k-2l-3}))\left(\prod_{q\in[0,l-1]\setminus[{\rm \bf{b}},{\rm \bf{c}}]}({\rm id}_{\mathcal{Y}}+\tilde{f}_{1,j_k-2q-2})\right)H_2[{\rm \bf{b}},{\rm \bf{c}}],\]
and from Theorem \ref{kashidem}, it coincides with (\ref{desireb}). \qed







\end{document}